\documentclass[12pt]{article}
\author{S\'everine Fiedler-Le Touz\'e}
\usepackage{amsfonts, amssymb, epsfig}
\newtheorem{theorem}{Theorem}

\newtheorem{lemma}{Lemma}
\newtheorem{proposition}{Proposition}
\title{Korchagin's third conjecture}
\begin{document}
\maketitle
\begin{abstract}
We consider the $M$-curve of degree nine with three nests $1 \langle \alpha_i \rangle, i = 1, 2, 3$ in $\mathbb{R}P^2$. After systematic constructions, Korchagin conjectured that at least two of the $\alpha_i$ must be odd. It was later proved that there is always one odd $\alpha_i$. 
We say that the curve has a jump in a non-empty oval $O$ if there exist four ovals $A$, $B$, $C$, $D$, with $A$ interior to some other non-empty oval $O'$, $D$ exterior, $B, C$ interior to $O$, such that $B$ and $C$ are separated inside of $O$ by any line passing through $A$ and $D$. 
In this paper, we prove the conjecture for the curves without jump, and we find restrictions on the complex orientations and rigid isotopy types admissible for the curves even, even, odd with jump. 
\end{abstract}

\section{Admissible complex types for the ninth degree $M$-curves with three nests}
\subsection{Introduction}
In \cite{ko5}-\cite{ko6}, Korchagin stated three conjectures about ninth degree $M$-curves, two of them are completely proved (or partly disproved), see\cite{fi1}, \cite{fi3}, \cite{or5}. The third one is still open. 
Let $C_9$ be a ninth degree $M$-curve with three nests, its real scheme is:
$\langle \mathcal{J} \amalg 1 \langle \alpha_1 \rangle
\amalg 1 \langle \alpha_2 \rangle \amalg 1 \langle \alpha_3 \rangle \amalg \beta\rangle$, with $\alpha_1 + \alpha_2 + \alpha_3 + \beta = 25$. 
\begin{quote}
{\bf Conjecture} At least two of the $\alpha_i$ are odd.
\end{quote}
It is already known that at least one of the $\alpha_i$ is odd, see \cite{fi}. We say that the curve has a jump in a non-empty oval $O$ if there exist four ovals $A$, $B$, $C$, $D$, with $A$ interior to some other non-empty oval $O'$, $D$ exterior, $B, C$ interior to $O$, such that $B$ and $C$ are separated inside of $O$ by any line passing through $A$ and $D$, see Figure~\ref{saut}. 
The main result of the present paper is the following:

\begin{figure}[htbp]
\centering
\includegraphics{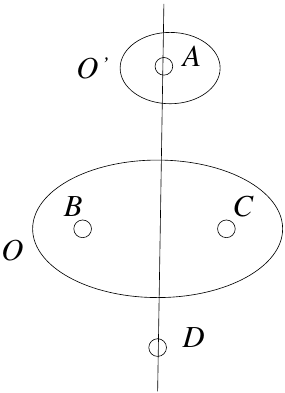}
\caption{\label{saut} Jump in $O$}
\end{figure}

\begin{theorem}
Let $C_9$ be an $M$-curve of degree $9$ with three nests and no jump.
Then at least two of the $\alpha_i$ are odd. 
\end{theorem}

\begin{figure}[htbp]
\centering
\includegraphics{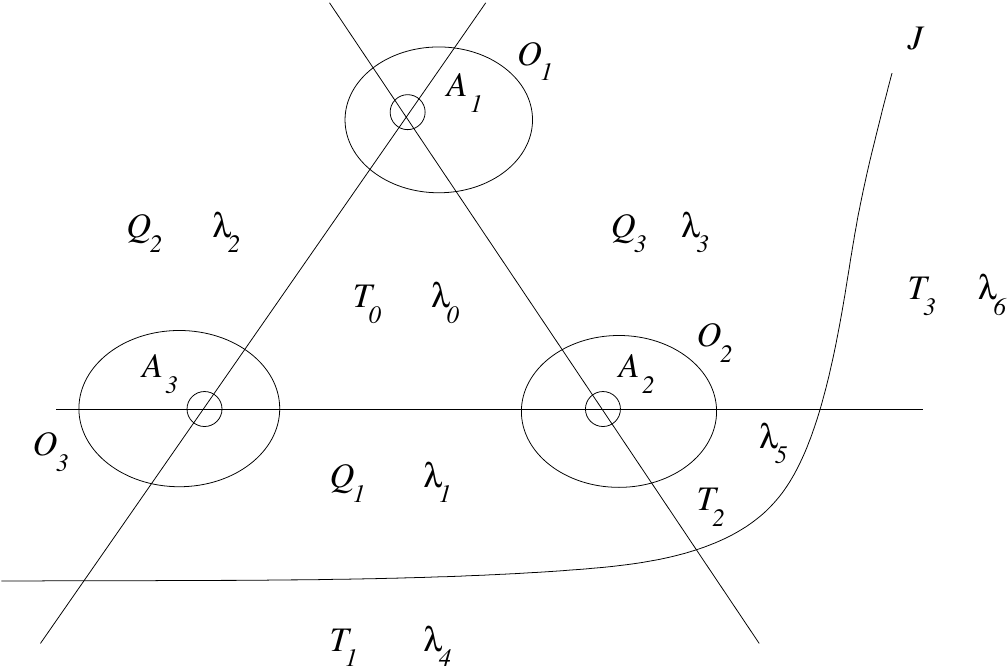}
\caption{\label{conj1} $M$-curve of degree $9$ with three nests}
\end{figure}
\begin{figure}[htbp]
\centering
\includegraphics{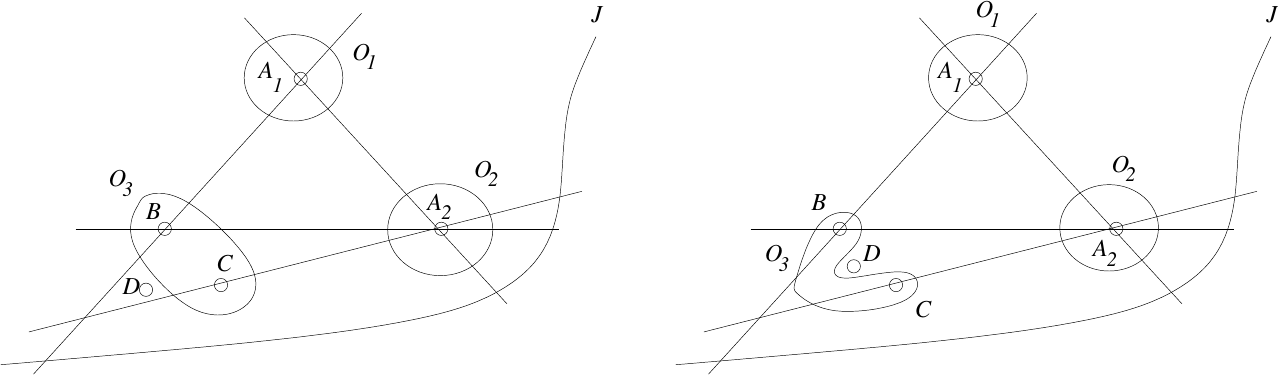}
\caption{\label{conj2} Jump in $O_3$, $D$ front, $D$ back}
\end{figure}

In this section, we determine as much information as possible about the complex orientations and rigid isotopy types realizable by the ninth degree $M$-curves with three nests (Proposition~1 and Tables~\ref{orev1}-\ref{orev16}).
To this purpose, we apply lemmas and propositions stated in \cite{fi}, in particular
Lemmas~19 and 20, which are handy reformulations, for the specific case of the ninth degree $M$-curves with three nests, of Orevkov's complex orientation formulas from \cite{or6}. 
Let us first describe these curves, using the results from \cite{fi}.
For the reader's convenience, we repeat here some definitions and notations introduced in that paper. 
For each empty oval $P$ of $C_9$, we choose a point lying in it, denoted also by $P$, it will be clear from the context whether we speak of the point or of the oval.
Let $O_1, O_2, O_3$ be the three non-empty ovals.
For any pair of ovals $A_i, A_j$, interior to $O_i$ and $O_j$ respectively, one denotes by $[A_iA_j]$ and calls {\em affine segment\/} the segment of the line $A_iA_j$ that doesn't cut $\mathcal{J}$, see Figure~\ref{conj1}. Let $[A_iA_j]'$ be the other segment. Any pencil $\mathcal{F}_A$, based in an oval $A$ interior to $O_j$ or $O_k$ and sweeping out $O_i$ from the first oval interior to $O_i$ to the last one, meets the same sequence (Fiedler chain) of empty ovals, each of them being either interior to $O_i$ or exterior. 
The occurence of exterior ovals means that $O_i$ has some jump. If $O_i$ has no jump, the sequence is called {\em chain of $O_i$\/}. If $\alpha_i$ is odd, the chain of $O_i$ is {\em positive\/} or {\em negative\/}, depending on whether it contributes by $+1$ or $-1$ to $\Lambda_+ - \Lambda_-$. 
The word {\em chain\/}, when used without further precision, has a very restrictive meaning. A pencil of lines (conics) sweeping out a curve of degree $m$ is {\em maximal\/} if it has alternatively $m$ ($2m$) and $m-2$ ($2m-2$) real intersection points with the curve. A {\em chain\/} is a Fiedler chain of ovals arising from a maximal pencil, in such a way that the consecutive ovals are connected by vanishing cycles. In the pseudoholomorphic setting, the curve may be perturbed so as to replace these cycles by real crossings.
It was proved in \cite{fi} that at most one non-empty oval may have a jump, by convention it is $O_3$. Moreover the sequence of $O_3$ (or jump sequence) splits into exactly three {\em jump subsequences\/} (int, ext, int), we say that there is only one jump. 
If these subsequences have parities (odd, odd, odd), we speak of an {\em odd jump\/}. The two inner subsequences are the {\em chains of $O_3$\/}. 
Let $O_i$ be any one of the non-empty ovals, the lines passing through pairs of ovals interior to $O_i$ cut either all of the affine segments $[A_jA_k]$ or all of the non-affine segments $[A_jA_k]'$, where $A_j$ ($A_k$) ranges over all of the ovals interior to $O_j$ ($O_k$). In the first case, we say that $O_i$ is {\em separating\/}, otherwise $O_i$ is {\em non-separating\/}. If $O_3$ has a jump, $O_3$ is non-separating. Let $B, D, C$ be three ovals distributed in the three successive jump subsequences. Up to a swap of $B$ and $C$, the ovals $A_1$, $A_2$, $C$, $D$, $B$ lie in convex position in $\mathbb{R}P^2 \setminus \mathcal{J}$ for any choice of $A_1, A_2$, see \cite{fi1} and \cite{fi}.
Let $O_3$ have a jump and $D$ be an exterior oval of the jump sequence. Then $O_3$ cuts all of the segments $[AD]$ or all of the segments $[AD]'$, when $A$ ranges over all ovals interior to $O_1$ and $O_2$. In the first case, we say that $D$ is {\em front\/}, in the second case, $D$ is {\em back\/}, see Figure~\ref{conj2}. If each exterior oval $D$ in the exterior jump subsequence is front (back), we say that $O_3$ is {\em crossing\/} ({\em non-crossing\/}), the exterior subsequence is a chain, the three jump subsequences will be called {\em jump chains\/}. Note that there may exist both front and back ovals, in this case, the exterior subsequence splits into two chains formed respectively by front and back ovals, see Proposition~3 in \cite{fi}.

\begin{figure}[htbp]
\centering
\includegraphics{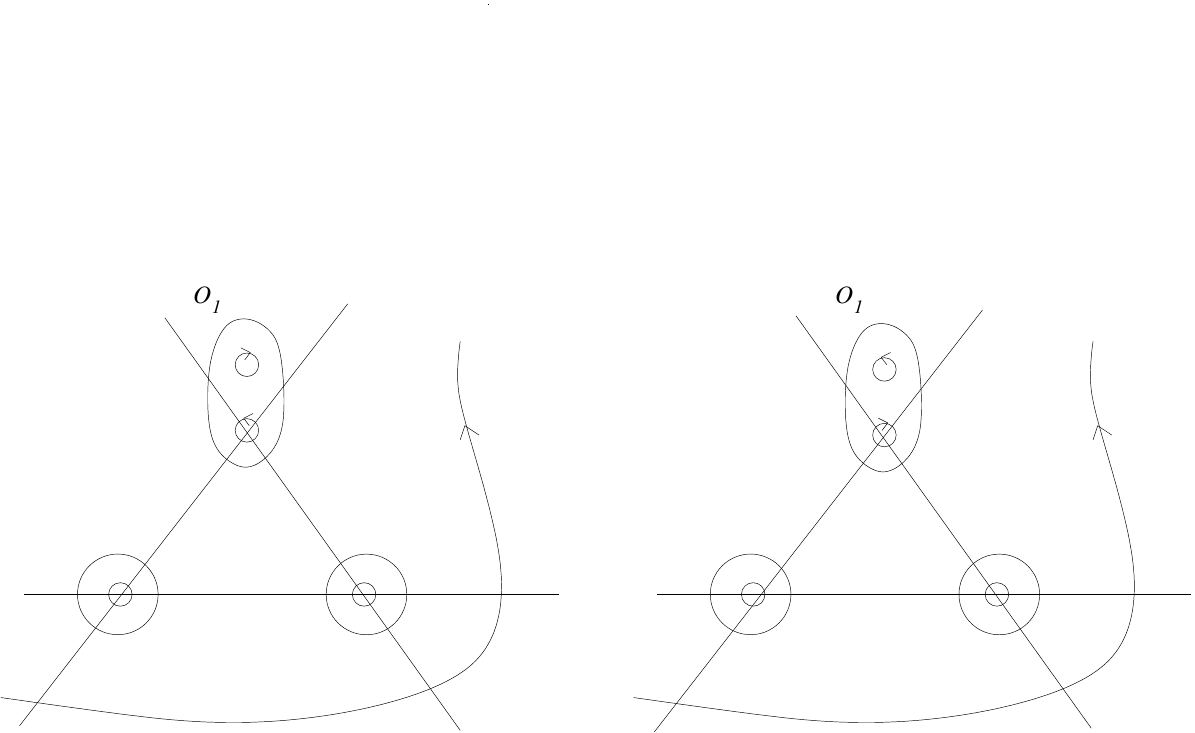}
\caption{\label{updown} $O_1$ up, $O_1$ down}
\end{figure}

Let $A_1, A_2, A_3$ be three chosen ovals, distributed inside of the three non-empty ovals $O_1, O_2, O_3$. The six ovals $O_i, A_i$ will be called the {\em principal ovals\/} of $C_9$. The three lines $A_1A_2$, $A_1A_3$, $A_2A_3$, and $\mathcal{J}$ divide the plane in four quadrangles and three triangles. The contributions of the ovals in $T_0$, $Q_1$, $Q_2$, $Q_3$, $T_1$, $T_2$, $T_3$ to $\Lambda_+ - \Lambda_-$ are denoted by $\lambda_i$, $i = 0 \dots 6$, see Figure~\ref{conj1}.
Let $A_i$ be extremity of the chain of $O_i$ if $O_i$ has no jump; if $O_3$ has a jump, $A_3$ is extremity of the jump sequence.
The interior of $O_i$ is divided in four parts by the lines $A_iA_j$, $A_iA_k$ (two in quadrangles, two in triangles), if $O_i$ has no jump, the chain of $O_i$ lies entirely in one of these four parts, if $O_3$ has a jump, the two chains of $O_3$ lie in the same part, contained in a quadrangles $Q_1$ or $Q_2$ (by convention, we choose $A_3$ such that the quadrangle is $Q_1$).
If $O_i$ is separating, the chain of $O_i$ lies in $T_i$ or in $T_0$. Let $O_i$ be separating and even. We say that $O_i$ is {\em up\/} if, when $A_i$ is chosen negative (positive), the chain is in $T_i$ ($T_0$), otherwise $O_i$ is {\em down\/}, see Figure~\ref{updown}. 

Let us assign to each nest a code that gives the sign of the outer oval $O_i$, and the number $\alpha_i^+ - \alpha_i^- \in \{ 0, \pm 1, \pm 2 \}$. 
Note that the values $\pm 2$ correspond to the nests with odd jump. The codes are defined as explained in the example hereafter: $+$ stands for an even nest (without odd jump) with positive outer oval, $(+, -)$ stands for an odd nest with positive outer oval and negative inner chain, $(+, -, -)$ stands for a nest with odd jump, such that the outer oval $O_3$ is positive and the two chains of $O_3$ are negative.
The code $(\pm, \mp)$ stands for a nest that is either $(+, -)$ or $(-, +)$.
The {\em short complex scheme\/} of $C_9$ is the set of three codes realized by the three nests. 
The {\em complex type\/} of $C_9$ is the short complex scheme, enhanced with the information for the odd nests, whether they are separating or non-separating; for the even nests, whether they are non-separating, up or down: $(+, -, n)$ and $(+, -, s)$ correspond respectively to a non-separating and a separating odd nest with positive outer oval $O_i$ and negative interior chain; $(+, n)$, $(+, u)$, $(+, d)$ correspond respectively to a non-separating, up, and down even nest with positive outer oval $O_i$.

In the next subsection, we will prove:

\begin{proposition}
The admissible complex types for $C_9$ with at least one even nest are: listed in Table~\ref{orev6} (even, even, odd without jump); listed in Table~\ref{orev10} (even, even, odd with jump); listed in Table~\ref{orev16} (even, odd, odd without jump); 
$(\pm, \mp, n)$, $(\pm, \mp, n)$, $(+, -, -)$ and $(\pm, \mp, n)$, $(\pm, \mp, n)$, $(-, +, +)$ (even, odd, odd with jump).

The admissible short complex schemes for $C_9$ with three odd nests are listed in Table~\ref{orev18}. 
\end{proposition}

The admissible complex types contradicting the conjecture are listed in Table~\ref{orev6} (case without jump), and Table~\ref{orev10} (case with jump).
Theorem~1 is proved in section~2.

\begin{figure}[htbp]
\centering
\includegraphics{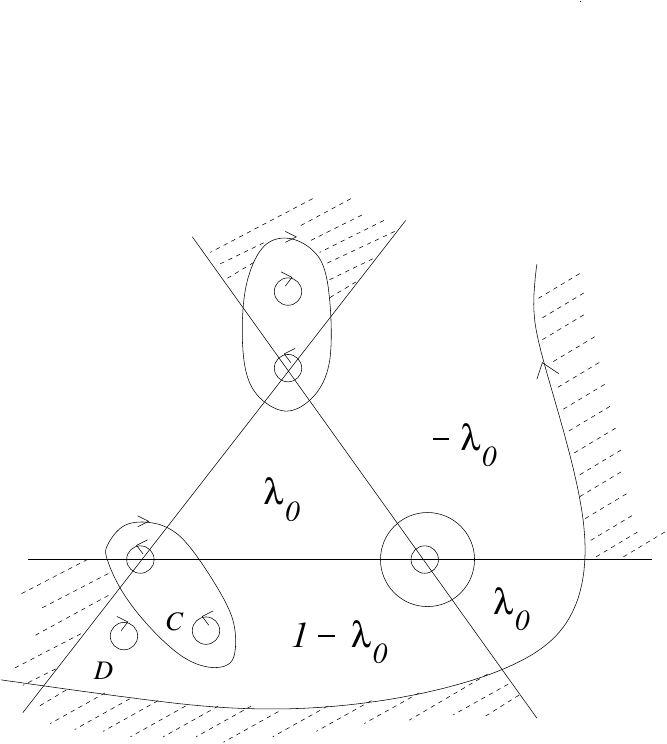}
\caption{\label{upcross} $C_9$ with complex type $(+, u), (\pm, \mp), (+, -, -)$, $\lambda_0 = 0$ or $1$ }
\end{figure}

\begin{figure}[htbp]
\centering
\includegraphics{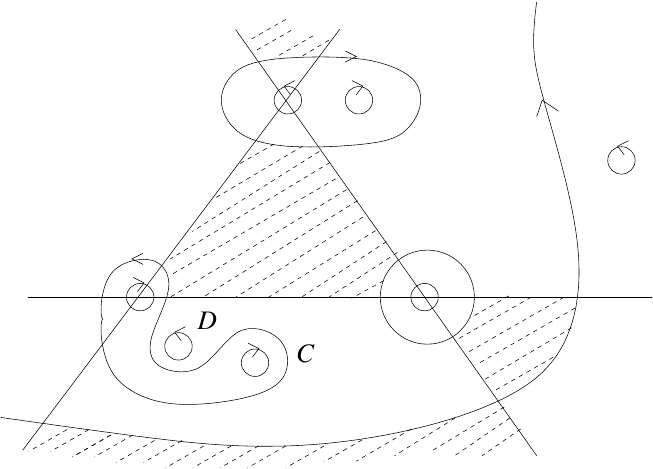}
\caption{\label{rightnocross} $C_9$ with complex type $(+, n), (\pm, \mp, n)$, $(-, +, +)$ }
\end{figure}

Let $O_1$ be even and non-separating. We say that $O_1$ is {\em left\/} ({\em right\/}) if, when $A_1$ was chosen with orientation opposite to that of $O_1$, the chain of $O_1$ lies in $Q_2$ ($Q_3)$.
In section~3, we will prove:

\begin{theorem}
Let $C_9$ be a ninth degree $M$-curve even, even, odd with jump.
Then $C_9$ is as shown in Figure~\ref{upcross} or Figure~\ref{rightnocross}. 
The curve $C_9$ has complex type $(+, u), (\pm, \mp), (+, -, -)$, with $O_3$ crossing and $\lambda_0 = 0$ or $1$ (Figure~\ref{upcross}); or $C_9$ has complex type $(+, n), (\pm, \mp, n)$, $(-, +, +)$ with $O_3$ non-crossing, $O_1$ right and $\lambda_0 = 0$ (Figure~\ref{rightnocross}). 
\end{theorem} 

Let us add a word of explanation to the Figures.
The principal ovals, the base lines and $\mathcal{J}$ divide the plane in zones.
We have indicated the contributions of the non-principal ovals to $\Lambda_+ - \Lambda_-$ in each zone (a single positive oval means $+1$). In a zone left blank, the contribution is $0$, the zones that we know to be empty are hatched.

\subsection{Proof of Proposition~1}
All along the proof, we refer to Lemmas and Propositions from \cite{fi}. 
For the reader's convenience, we give here a quick summary of the main results. 
If $O_3$ is crossing, $T_3$ is empty; if $O_3$ is non-crossing, $T_0 \cup T_1 \cup T_2$ is emtpy; if the exterior jump sequence has both front and back ovals,
$T_0 \cup T_1 \cup T_2 \cup T_3$ is empty (Lemma~9).
Let $\Lambda = \lambda_0 - \lambda_4 - \lambda_5 - \lambda_6$. One has $\Lambda = \Pi_+ - \Pi_- - 4$ (Lemma~10). 
Let $C_9$ have a jump, then $\Pi_+ - \Pi_- = 3$ or $4$.
If $\Pi_+ - \Pi_- = 4$, then $C_9$ is even, odd, odd with odd jump.
If $\Pi_+ - \Pi_- = 3$, then $O_3$ is either crossing and positive, or non-crossing and negative (Lemma~18). 
Note that if $\Pi_+ - \Pi_- = 3$, the curve $C_9$ may be even, even, odd with odd jump, or odd, odd, odd with jump sequence either (odd, odd, even) or (odd, even, even).  
In \cite{fi}, we defined parameters $\pi_l$, $\pi'_l$, $N_l$, $M_l$, $G_l$, $F_l$ associated to the code of each nest $O_l$ ($F_l$ is the only one that depends on the data $n, s, d, u$). We defined also parameters $E_i, i = 0, 1, 2, 3$ depending on the three codes. Lemmas~19 and 20 are handy reformulations of Orevkov's formulas: 
If $T_i$ contains some exterior ovals, then $E_i = 0$ (Lemma~19). If $O_i$ is separating, then $F_i - G_j - G_k = 0$ (Lemma~20). 
If $T_0 \cup T_1 \cup T_2 \cup T_3$ is empty, then the short complex scheme of $C_9$ is $(\pm, \mp), (\pm, \mp), (+, -, -)$ or $(\pm, \mp), (\pm, \mp), (-, +, +)$ (Lemma~21). 
Finally, let us mention two propositions from \cite{fi}:
If $T_i$ contains some exterior ovals, then $\vert \lambda_{i+3} \vert \leq 3$; if $\vert \lambda_{i+3} \vert = 3$, then 
$\lambda_{i+3} = 3$ and $\Lambda = -2$ (Proposition~1). If $T_0$ contains only
exterior ovals, then $\vert \lambda_0 \vert \leq 2$ (Proposition~2).

Consider first the case even, even, odd. 
The twelve admissible complex schemes without jump are displayed in Table~\ref{orev1}, along with: the numbers $E_i$, the indices $j$ of the zones $T_j$ that may contain some exterior ovals (row $\mathcal{Z}$), and $\Lambda$. Note that by Lemma~21, the first five complex schemes of Table~\ref{orev1} cannot be realized with all of the $O_i$ non-separating. 
Adding to this the informations in Tables~\ref{orev2}-\ref{orev4} allows to get the list of admissible complex types in Table~\ref{orev5}. 
In the tables, each row with $n = 0, 1, 2, 3$ nests $(\pm, \mp)$ corresponds actually to $n+1$ complex schemes or types. Choose for any separating even nest $O_i$ the base oval $A_i$ in such a way that $(A_i, O_i)$ form a positive pair, and $A_i$ is extremal in the chain of inner ovals. 
The complex types in Table~\ref{orev5} are such that $\Lambda = \Pi_+ - \Pi_- - 4 = -3$ if $O_3$ is $(\pm, \mp)$, and $-5$ if $O_3$ is $(+, +)$ or $(-, -)$.
This, combined with the data of $\mathcal{Z}$, allows to find out the admissible values for the triangular parameters $\lambda_i, i = 0, 4, 5, 6$.
The six types with $\mathcal{Z} = (1)$ (first four rows of Table~\ref{orev5}) are such that the triangle $T_1$ contains only exterior ovals, as $O_1$ is either non-separating or $(+, d)$. For these types, $\lambda_4 = 3, 4, 5, 6$, this contradicts Proposition~1 (from \cite{fi}).
The type with $\mathcal{Z} = (3)$ is such that the triangle $T_3$ contains only exterior ovals. But $\lambda_6 = 5$, this contradicts also Proposition~1.
The types with $\mathcal{Z} = (0)$ are such that the triangle $T_0$ contains only exterior ovals. But $\lambda_0 \leq - 5$, this contradicts Proposition~2.
The types $(-, d), (-, d), (\pm, \mp, n)$ and  $(-, d), (-, d), (\pm, \mp, s)$ verify $\lambda_0 - \lambda_6 = -5$. For $O_3$ non-separating, both $T_0$ and $T_3$ contain only exterior ovals, so that one may apply Propositions~1 and 2.
One must have $\vert \lambda_6 \vert \leq 2$, $\vert \lambda_0 \vert \leq 2$, this is a contradiction. When $O_3$ is separating, one may choose the oval $A_3$ so that $T_3$ contains only exterior ovals, to apply Proposition~1.
Then, one chooses the other extremal oval $A'_3$, so that $T'_0$ contains only exterior ovals, to apply Proposition~2.The values of the $\lambda_i$ are the same for both choices, so one gets the same contradiction as for a non-separating $O_3$. In the sequel, keep in mind that we only need to take care of the even nests to achieve the hypothesis needed for Proposition~1 or 2. 
The types left are gathered in Table~\ref{orev6}, note that each of the three rows corresponds actually to four types, depending on whether the nest $O_3$ is $(+, -)$ or $(-, +)$, separating or non-separating. We get the admissible values of $\lambda_0$ and $\lambda_6$ applying Propositions~1 and 2. As $T_0$ contains only exterior ovals, $\vert \lambda_0 \vert \leq 2$. If the even ovals are both non-separating, $\vert \lambda_6 \vert = \vert \lambda_0 + 3 \vert \leq 2$, hence $\lambda_0 = -1$ or $-2$; if $O_1$ is separating, $\vert \lambda_6 \vert = \vert \lambda_0 + 4 \vert \leq 2$, hence $\lambda_0 = -2$.  

The $M$-curves with odd jump are either even, even, odd with $\Pi_+ - \Pi_- = 3$
or even, odd, odd with $\Pi_+ - \Pi_- = 4$ (Lemma~18).
The list of admissible short complex schemes with odd jump is shown in Table~\ref{orev7}. Consider the case even, even, odd. If $O_3$ is non-crossing, the ovals $O_1$ and $O_2$ are non-separating (by Lemma~9). If moreover $O_1$ is negative, then there are no exterior triangular ovals. This is a contradiction as one should have $\lambda_6 = -\Lambda = 1$.
Adding this to the informations in Tables~\ref{orev8}-\ref{orev9} allows to get the list of admissible complex types even, even, odd with jump displayed in Table~\ref{orev10}, and to see that the ovals $O_1$, $O_2$ of the even, odd, odd curves with jump must be non-separating. In Table~\ref{orev10}, we have also indicated the parameters $\lambda_1$, $\lambda_2$, $\lambda_3$, they are
obtained from the identities $\lambda_0 + \lambda_1 - \lambda_4 = 0$, 
$\lambda_0 + \lambda_2 - \lambda_5 = 0$, $\lambda_0 + \lambda_3 - \lambda_6 = 0$. (These identities follow from Fiedler's theorem applied with the pencils of lines $\mathcal{F}_{A_i}$ sweeping out $T_0 \cup Q_i \cup T_i$, see first part of 
Lemma~ 10). 

Table~\ref{orev11} shows the admissible short complex schemes even, odd, odd without jump. As there is at most one separating even nest, the contribution of the interior ovals to $\Lambda$ is $0, +1$ or $-1$. If $O_1$ is separating, choose the base oval $A_1$ so as to realize each time the required hypothesis to apply Propositions~1 or 2. We get thus:
if $\mathcal{Z} = \emptyset$, then $\vert \Lambda \vert \leq 1$; 
if $\mathcal{Z} = (0)$, then $\vert \lambda_0 \vert \leq 2$ hence 
$\vert \Lambda \vert \leq 3$;
if $\mathcal{Z} = (i), i = 1, 2, 3$, then $\vert \lambda_{i+3} \vert \leq 2$ hence
$\vert \Lambda \vert \leq 3$, or $\lambda_{i+3} = 3$ and $\Lambda = -2$. 
Using the informations from Tables~\ref{orev12}-\ref{orev15}, we get the list of admissible complex types even, odd, odd without jump in Table~\ref{orev16}.

Table~\ref{orev17} shows the admissible short complex schemes odd, odd, odd.
By Propositions~1 and 2:  $\vert \lambda_0 \vert \leq 2$ and 
$\vert \lambda_{i+3} \vert \leq 2$. The only admissible schemes left are shown in Table~\ref{orev18}. For $O_1 = (\pm, \mp, s)$, one has $F_1 - G_2 - G_3 = 0$ for both choices of $O_3$, there is no contradiction (see Lemma~20). For $O_3 = (-, -, s)$, one has $F_3 - G_1 - G_2 = -1$, contradiction.
Any one of the nests $(\pm, \mp)$ may be separating, whereas the nest $(-, -)$ must be non-separating.
The curve may have a jump only if all three nests are $(\pm, \mp)$
(Lemma~18).

\begin{table}
\centering
\begin{tabular}{c c c c c c c c c }
$O_1$ & $O_2$ & $O_3$ & $E_0$ & $E_1$ & $E_2$ & $E_3$ & $\mathcal{Z}$ & $\Lambda$\\
$+$ & $+$ & $(+, +)$ & $-4$ & $-3$ & $-3$ & $-2$ & $\emptyset$ & $-5$\\
$+$ & $+$ & $(+, -)$ & $-2$ & $-1$ & $-1$ & $-2$ & $\emptyset$ & $-3$\\
$+$ & $+$ & $(-, +)$ & $-2$ & $-1$ & $-1$ & $-2$ & $\emptyset$ & $-3$\\
$+$ & $+$ & $(-, -)$ & $-2$ & $-1$ & $-1$ & $-4$ & $\emptyset$ & $-5$\\
$+$ & $-$ & $(+, +)$ & $-3$ & $-2$ & $-4$ & $-1$ & $\emptyset$ & $-5$\\
$+$ & $-$ & $(+, -)$ & $-1$ & $0$ & $-2$ & $-1$ & $(1)$ & $-3$\\
$+$ & $-$ & $(-, +)$ & $-1$ & $0$ & $-2$ & $-1$ & $(1)$ & $-3$\\
$+$ & $-$ & $(-, -)$ & $-1$ & $0$ & $-2$ & $-3$ & $(1)$ & $-5$\\
$-$ & $-$ & $(+, +)$ & $-2$ & $-3$ & $-3$ & $0$ & $(3)$ & $-5$\\
$-$ & $-$ & $(+, -)$ & $0$ & $-1$ & $-1$ & $0$ & $(0, 3)$ & $-3$\\
$-$ & $-$ & $(-, +)$ & $0$ & $-1$ & $-1$ & $0$ & $(0, 3)$ & $-3$\\
$-$ & $-$ & $(-, -)$ & $0$ & $-1$ & $-1$ & $-2$ & $(0)$ & $-5$\\
\end{tabular}
\caption{\label{orev1} Admissible short complex schemes even, even, odd, without jump}
\end{table}

\begin{table}
\centering
\begin{tabular}{c c c c}
$O_1$ & $O_2$ & $O_3$ & $F_1 - G_2 - G_3$\\
$(+, u)$ & $+$ & $(+, +)$ & $-4$\\
$(+, d)$ & $+$ & $(+, +)$ & $-3$\\
$(+, u)$ & $+$ & $(\pm, \mp)$ & $-2$\\
$(+, d)$ & $+$ & $(\pm, \mp)$ & $-1$\\
$(+, u)$ & $+$ & $(-, -)$ & $-2$\\
$(+, d)$ & $+$ & $(-, -)$ & $-1$\\
$(+, u)$ & $-$ & $(+, +)$ & $-3$\\
$(+, d)$ & $-$ & $(+, +)$ & $-2$\\
$(+, u)$ & $-$ & $(\pm, \mp)$ & $-1$\\
$(+, d)$ & $-$ & $(\pm, \mp)$ & $0$\\
$(+, u)$ & $-$ & $(-, -)$ & $-1$\\
$(+, d)$ & $-$ & $(-, -)$ & $0$\\
$(-, u)$ & $-$ & $(+, +)$ & $-3$\\
$(-, d)$ & $-$ & $(+, +)$ & $-2$\\
$(-, u)$ & $-$ & $(\pm, \mp)$ & $-1$\\
$(-, d)$ & $-$ & $(\pm, \mp)$ & $0$\\
$(-, u)$ & $-$ & $(-, -)$ & $-1$\\
$(-, d)$ & $-$ & $(-, -)$ & $0$\\
\end{tabular}
\caption{\label{orev2} Even, even, odd, without jump $O_1$ separating}
\end{table}

\begin{table}
\centering
\begin{tabular}{c c c c}
$O_1$ & $O_2$ & $O_3$ & $F_2 - G_1 - G_3$\\
$+$ & $(-, u)$ & $(+, +)$ & $-4$\\
$+$ & $(-, d)$ & $(+, +)$ & $-3$\\
$+$ & $(-, u)$ & $(\pm, \mp)$ & $-2$\\
$+$ & $(-, d)$ & $(\pm, \mp)$ & $-1$\\
$+$ & $(-, u)$ & $(-, -)$ & $-2$\\
$+$ & $(-, d)$ & $(-, -)$ & $-1$\\
\end{tabular}
\caption{\label{orev3} Even, even, odd, without jump, $O_2$ separating}
\end{table}

\begin{table}
\centering
\begin{tabular}{c c c c}
$O_1$ & $O_2$ & $O_3$ & $F_3 - G_1 - G_2$\\
$+$ & $+$ & $(+, +, s)$ & $-3$\\
$+$ & $+$ & $(\pm, \mp, s)$ & $-2$\\
$+$ & $+$ & $(-, -, s)$ & $-3$\\
$+$ & $-$ & $(+, +, s)$ & $-2$\\
$+$ & $-$ & $(\pm, \mp, s)$ & $-1$\\
$+$ & $-$ & $(-, -, s)$ & $-2$\\
$-$ & $-$ & $(+, +, s)$ & $-1$\\
$-$ & $-$ & $(\pm, \mp, s)$ & $0$\\
$-$ & $-$ & $(-, -, s)$ & $-1$\\
\end{tabular}
\caption{\label{orev4} Even, even, odd, without jump, $O_3$ separating}
\end{table}

\begin{table}
\centering
\begin{tabular}{c c c c c c c c }
$O_1$ & $O_2$ & $O_3$ & $\mathcal{Z}$ & $\lambda_0$ & $\lambda_4$ & $\lambda_5$ & $\lambda_6$\\
$(+, n)$ & $(-, n)$ & $(\pm, \mp, n)$ & $(1)$ & $0$ & $3$ & $0$ & $0$\\
$(+, d)$ & $(-, n)$ & $(\pm, \mp, n)$ & $(1)$ & $1$ & $4$ & $0$ & $0$\\
$(+, n)$ & $(-, n)$ & $(-, -, n)$ & $(1)$ & $0$ & $5$ & $0$ & $0$\\
$(+, d)$ & $(-, n)$ & $(-, -, n)$ & $(1)$ & $1$ & $6$ & $0$ & $0$\\
$(-, n)$ & $(-, n)$ & $(+, +, n)$ & $(3)$ & $0$ & $0$ & $0$ & $5$\\
$(-, n)$ & $(-, n)$ & $(\pm, \mp, n)$ & $(0, 3)$ & $\lambda_0$ & $0$ & $0$ & $\lambda_0 +3$\\
$(-, d)$ & $(-, n)$ & $(\pm, \mp, n)$ & $(0, 3)$ & $\lambda_0$ & $-1$ & $0$ & $\lambda_0 + 4$\\
$(-, n)$ & $(-, n)$ & $(\pm, \mp, s)$ & $(0, 3)$ & $\lambda_0$ & $0$ & $0$ & $\lambda_0 +3$\\
$(-, d)$ & $(-, n)$ & $(\pm, \mp, s)$ & $(0, 3)$ & $\lambda_0$ & $-1$ & $0$ & $\lambda_0 + 4$\\
$(-, d)$ & $(-, d)$ & $(\pm, \mp, n)$ & $(0, 3)$ & $\lambda_0$ & $-1$ & $-1$ & $\lambda_0 + 5$\\
$(-, d)$ & $(-, d)$ & $(\pm, \mp, s)$ & $(0, 3)$ & $\lambda_0$ & $-1$ & $-1$ & $\lambda_0 + 5$\\
$(-, n)$ & $(-, n)$ & $(-, -, n)$ & $(0)$ & $-5$ & $0$ & $0$ & $0$\\
$(-, d)$ & $(-, n)$ & $(-, -, n)$ & $(0)$ & $-6$ & $-1$ & $0$ & $0$\\
$(-, d)$ & $(-, d)$ & $(-, -, n)$ & $(0)$ & $-7$ & $-1$ & $-1$ & $0$\\
\end{tabular}
\caption{\label{orev5} Admissible complex types even, even, odd without jump}
\end{table}

\begin{table}
\centering
\begin{tabular}{c c c c c c c c c c }
$O_1$ & $O_2$ & $O_3$ & $\lambda_0$ & $\lambda_1$ & $\lambda_2$ & $\lambda_3$ & $\lambda_4$ & $\lambda_5$ & $\lambda_6$\\
$(-, n)$ & $(-, n)$ & $(\pm, \mp)$ & $-1$ & $1$ & $1$ & $3$ & $0$ & $0$ & $2$\\
$(-, n)$ & $(-, n)$ & $(\pm, \mp)$ & $-2$ & $2$ & $2$ & $3$ & $0$ & $0$ & $1$\\
$(-, d)$ & $(-, n)$ & $(\pm, \mp)$ & $-2$ & $1$ & $2$ & $4$ & $-1$ & $0$ & $2$\\
\end{tabular}
\caption{\label{orev6} Even, even, odd, without jump, last cases left}
\end{table}

\begin{table}
\centering
\begin{tabular}{c c c c c c c c c }
$O_1$ & $O_2$ & $O_3$ & $E_0$ & $E_1$ & $E_2$ & $E_3$ & $\mathcal{Z}$ & $\Lambda$\\
$+$ & $(+, -)$ & $(+, -, -)$ & $0$ & $1$ & $0$ & $-1$ & $(0, 2)$ & $-1$\\
$-$ & $(+, -)$ & $(+, -, -)$ & $1$ & $0$ & $1$ & $0$ & $(1)$ & $-1$\\
$+$ & $(-, +)$ & $(+, -, -)$ & $0$ & $1$ & $0$ & $-1$ & $(0, 2)$ & $-1$\\
$-$ & $(-, +)$ & $(+, -, -)$ & $1$ & $0$ & $1$ & $0$ & $(1)$ & $-1$\\
$(+, -)$ & $(+, -)$ & $(+, -,-)$ & $1$ & $1$ & $1$ & $0$ & $(3)$ & $0$\\
$(+, -)$ & $(-, +)$ & $(+, -,-)$ & $1$ & $1$ & $1$ & $0$ & $(3)$ & $0$\\
$(-, +)$ & $(-, +)$ & $(+, -,-)$ & $1$ & $1$ & $1$ & $0$ & $(3)$ & $0$\\
$+$ & $(+, -)$ & $(-, +, +)$ & $-1$ & $0$ & $-1$ & $0$ & $(3)$ & $0$\\
$-$ & $(+, -)$ & $(-, +, +)$ & $0$ & $-1$ & $0$ & $1$ & $\emptyset$ & $-1$\\
$+$ & $(-, +)$ & $(-, +, +)$ & $-1$ & $0$ & $-1$ & $0$ & $(3)$ & $-1$\\
$-$ & $(-, +)$ & $(-, +, +)$ & $0$ & $-1$ & $0$ & $1$ & $\emptyset$ & $-1$\\
$(+, -)$ & $(+, -)$ & $(-, +, +)$ & $0$ & $0$ & $0$ & $1$ & $(0, 1, 2)$ & $0$\\
$(+, -)$ & $(-, +)$ & $(-, +, +)$ & $0$ & $0$ & $0$ & $1$ & $(0, 1, 2)$ & $0$\\
$(-, +)$ & $(-, +)$ & $(-, +, +)$ & $0$ & $0$ & $0$ & $1$ & $(0, 1, 2)$ & $0$\\
\end{tabular}
\caption{\label{orev7} Admissible short complex schemes with odd jump}
\end{table}

\begin{table}
\centering
\begin{tabular}{c c c c}
$O_1$ & $O_2$ & $O_3$ & $F_1 - G_2 - G_3$\\
$(+, u)$ & $(\pm, \mp)$ & $(+, -, -)$ & $0$\\
$(+, d)$ & $(\pm, \mp)$ & $(+, -, -)$ & $1$\\
$(-, u)$ & $(\pm, \mp)$ & $(+, -, -)$ & $0$\\
$(-, d)$ & $(\pm, \mp)$ & $(+, -, -)$ & $1$\\
$(\pm, \mp, s)$ & $(\pm, \mp)$ & $(+, -,-)$ & $1$\\
\end{tabular}
\caption{\label{orev8} Odd jump, $O_1$ separating}
\end{table}

\begin{table}
\centering
\begin{tabular}{c c c c}
$O_1$ & $O_2$ & $O_3$ & $F_2 - G_1 - G_3$\\
$+$ & $(\pm, \mp, s)$ & $(+, -, -)$ & $0$\\
$-$ & $(\pm, \mp, s)$ & $(+, -, -)$ & $1$\\
\end{tabular}
\caption{\label{orev9} Odd jump, $O_2$ separating}
\end{table}

\begin{table}
\centering
\begin{tabular}{c c c c c c c c c c c}
$O_1$ & $O_2$ & $O_3$ & $\mathcal{Z}$ & $\lambda_0$ & $\lambda_4$ & $\lambda_5$ & $\lambda_6$ & $\lambda_1$ & $\lambda_2$ & $\lambda_3$\\
$(+, n)$ & $(\pm, \mp, n)$ & $(+, -, -)$ & $(0, 2)$ & $\lambda_0$ & $0$ & $\lambda_0 + 1$ & $0$ & $-\lambda_0$ & $1$ & $-\lambda_0$\\
$(+, u)$ & $(\pm, \mp, n)$ & $(+, -, -)$ & $(0, 2)$ & $\lambda_0$ & $1$ & $\lambda_0$ & $0$ & $1-\lambda_0$ & $0$ & $-\lambda_0$\\
$(+, n)$ & $(\pm, \mp, s)$ & $(+, -, -)$ & $(0, 2)$ & $\lambda_0$ & $0$ & $\lambda_0 + 1$ & $0$ & $-\lambda_0$ & $1$ & $-\lambda_0$\\
$(+, u)$ & $(\pm, \mp, s)$ & $(+, -, -)$ & $(0, 2)$ & $\lambda_0$ & $1$ & $\lambda_0$ & $0$ & $1-\lambda_0$ & $0$ & $-\lambda_0$\\
$(-, n)$ & $(\pm, \mp, n)$ & $(+, -, -)$ & $(1)$ & $0$ & $1$ & $0$ & $0$ & $1$ & $0$ & $0$\\
$(-, u)$ & $(\pm, \mp, n)$ & $(+, -, -)$ & $(1)$ & $-1$ & $0$ & $0$ & $0$ & $1$ & $1$ & $1$\\
$(+, n)$ & $(\pm, \mp, n)$ & $(-, +, +)$ & $(3)$ & $0$ & $0$ & $0$ & $1$ & $0$ & $0$ & $1$\\
\end{tabular}
\caption{\label{orev10} Admissible complex types even, even, odd with jump}
\end{table}

\begin{table}
\centering
\begin{tabular}{c c c c c c c c c}
$O_1$ & $O_2$ & $O_3$ & $E_0$ & $E_1$ & $E_2$ & $E_3$ & $\mathcal{Z}$ & $\Lambda$\\
$+$ & $(\pm, \mp)$ & $(+, +)$ & $-3$ & $-2$ & $-3$ & $-1$ & $\emptyset$ & $-4$\\
$+$ & $(\pm, \mp)$ & $(\pm, \mp)$ & $-1$ & $0$ & $-1$ & $-1$ & $(1)$ & $-2$\\
$+$ & $(\pm, \mp)$ & $(-, -)$ & $-1$ & $0$ & $-1$ & $-3$ & $(1)$ & $-4$\\
$+$ & $(+, +)$ & $(+, +)$ & $-5$ & $-4$ & $-3$ & $-3$ & $\emptyset$ & $-6$\\
$+$ & $(+, +)$ & $(-, -)$ & $-3$ & $-2$ & $-1$ & $-5$ & $\emptyset$ & $-6$\\
$+$ & $(-, -)$ & $(-, -)$ & $-1$ & $0$ & $-3$ & $-3$ & $(1)$ & $-6$\\
$-$ & $(\pm, \mp)$ & $(+, +)$ & $-2$ & $-3$ & $-2$ & $0$ & $(3)$ & $-4$\\
$-$ & $(\pm, \mp)$ & $(\pm, \mp)$ & $0$ & $-1$ & $0$ & $0$ & $(0, 2, 3)$ & $-2$\\
$-$ & $(\pm, \mp)$ & $(-, -)$ & $0$ & $-1$ & $0$ & $-2$ & $(0, 2)$ & $-4$\\
$-$ & $(+, +)$ & $(+, +)$ & $-4$ & $-5$ & $-2$ & $-2$ & $\emptyset$ & $-6$\\
$-$ & $(+, +)$ & $(-, -)$ & $-2$ & $-3$ & $0$ & $-4$ & $(2)$ & $-6$\\
$-$ & $(-, -)$ & $(-, -)$ & $0$ & $-1$ & $-2$ & $-2$ & $(0)$ & $-6$\\
\end{tabular}
\caption{\label{orev11} Admissible short complex schemes even, odd, odd without jump}
\end{table}

\begin{table}
\centering
\begin{tabular}{c c c c}
$O_1$ & $O_2$ & $O_3$ & $F_3 - G_1 - G_2$\\
$+$ & $(\pm, \mp)$ & $(+, +, s)$ & $-2$\\
$+$ & $(\pm, \mp)$ & $(\pm, \mp, s)$ & $-1$\\
$+$ & $(\pm, \mp)$ & $(-, -, s)$ & $-2$\\
$+$ & $(+, +)$ & $(+, +, s)$ & $-4$\\
$+$ & $(+, +)$ & $(-, -, s)$ & $-4$\\
$+$ & $(-, -)$ & $(-, -, s)$ & $-2$\\
$-$ & $(\pm, \mp)$ & $(+, +, s)$ & $-1$\\
$-$ & $(\pm, \mp)$ & $(\pm, \mp, s)$ & $0$\\
$-$ & $(\pm, \mp)$ & $(-, -, s)$ & $-1$\\
$-$ & $(+, +)$ & $(+, +, s)$ & $-3$\\
$-$ & $(+, +)$ & $(-, -, s)$ & $-3$\\
$-$ & $(-, -)$ & $(-, -, s)$ & $-1$\\
\end{tabular}
\caption{\label{orev12} Even, odd, odd, $O_3$ separating}
\end{table}

\begin{table}
\centering
\begin{tabular}{c c c c}
$O_1$ & $O_2$ & $O_3$ & $F_2 - G_1 - G_3$\\
$+$ & $(\pm, \mp, s)$ & $(+, +)$ & $-3$\\
$+$ & $(\pm, \mp, s)$ & $(-, -)$ & $-1$\\
$+$ & $(+, +, s)$ & $(-, -)$ & $-2$\\
$-$ & $(\pm, \mp, s)$ & $(+, +)$ & $-2$\\
$+$ & $(\pm, \mp, s)$ & $(-, -)$ & $-1$\\
$-$ & $(+, +, s)$ & $(-, -)$ & $-1$\\
$-$ & $(\pm, \mp, s)$ & $(-, -)$ & $0$\\
\end{tabular}
\caption{\label{orev13} Even, odd, odd, $O_2$ separating}
\end{table}

\begin{table}
\centering
\begin{tabular}{c c c c}
$O_1$ & $O_2$ & $O_3$ & $F_1 - G_2 - G_3$\\
$(+, u)$ & $(\pm, \mp)$ & $(+, +)$ & $-3$\\
$(+, u)$ & $(\pm, \mp)$ & $(\pm, \mp)$ & $-1$\\
$(+, u)$ & $(\pm, \mp)$ & $(-, -)$ & $-1$\\
$(+, u)$ & $(+, +)$ & $(+, +)$ & $-5$\\
$(+, u)$ & $(+, +)$ & $(-, -)$ & $-3$\\
$(+, u)$ & $(-, -)$ & $(-, -)$ & $-1$\\
$(-, u)$ & $(\pm, \mp)$ & $(+, +)$ & $-3$\\
$(-, u)$ & $(\pm, \mp)$ & $(\pm, \mp)$ & $-1$\\
$(-, u)$ & $(\pm, \mp)$ & $(-, -)$ & $-1$\\
$(-, u)$ & $(+, +)$ & $(+, +)$ & $-5$\\
$(-, u)$ & $(+, +)$ & $(-, -)$ & $-3$\\
$(-, u)$ & $(-, -)$ & $(-, -)$ & $-1$\\
\end{tabular}
\caption{\label{orev14} Even, odd, odd, $O_1$ up}
\end{table}

\begin{table}
\centering
\begin{tabular}{c c c c}
$O_1$ & $O_2$ & $O_3$ & $F_1 - G_2 - G_3$\\
$(+, d)$ & $(\pm, \mp)$ & $(+, +)$ & $-2$\\
$(+, d)$ & $(\pm, \mp)$ & $(\pm, \mp)$ & $0$\\
$(+, d)$ & $(\pm, \mp)$ & $(-, -)$ & $0$\\
$(+, d)$ & $(+, +)$ & $(+, +)$ & $-4$\\
$(+, d)$ & $(+, +)$ & $(-, -)$ & $-2$\\
$(+, d)$ & $(-, -)$ & $(-, -)$ & $0$\\
$(-, d)$ & $(\pm, \mp)$ & $(+, +)$ & $-2$\\
$(-, d)$ & $(\pm, \mp)$ & $(\pm, \mp)$ & $0$\\
$(-, d)$ & $(\pm, \mp)$ & $(-, -)$ & $0$\\
$(-, d)$ & $(+, +)$ & $(+, +)$ & $-4$\\
$(-, d)$ & $(+, +)$ & $(-, -)$ & $-2$\\
$(-, d)$ & $(-, -)$ & $(-, -)$ & $0$\\
\end{tabular}
\caption{\label{orev15} Even, odd, odd, $O_1$ down}
\end{table}

\begin{table}
\centering
\begin{tabular}{c c c c c }
$O_1$ & $O_2$ & $O_3$ & $\mathcal{Z}$ & $\Lambda$\\
$(+, n)$ & $(\pm, \mp, n)$ & $(\pm, \mp, n)$ & $(1)$ & $-2$\\
$(+, d)$ & $(\pm, \mp, n)$ & $(\pm, \mp, n)$ & $(1)$ & $-2$\\ 
$(-, n)$ & $(\pm, \mp, n)$ & $(\pm, \mp, n)$ & $(0, 2, 3)$ & $-2$\\ 
$(-, n)$ & $(\pm, \mp, n)$ & $(\pm, \mp, s)$ & $(0, 2, 3)$ & $-2$\\ 
$(-, n)$ & $(\pm, \mp, s)$ & $(\pm, \mp, s)$ & $(0, 2, 3)$ & $-2$\\ 
$(-, d)$ & $(\pm, \mp, n)$ & $(\pm, \mp, n)$ & $(0, 2, 3)$ & $-2$\\ 
$(-, d)$ & $(\pm, \mp, n)$ & $(\pm, \mp, s)$ & $(0, 2, 3)$ & $-2$\\ 
$(-, d)$ & $(\pm, \mp, s)$ & $(\pm, \mp, s)$ & $(0, 2, 3)$ & $-2$\\ 
$(-, n)$ & $(\pm, \mp, n)$ & $(-, -, n)$ & $(0, 2)$ & $-4$\\ 
$(-, n)$ & $(\pm, \mp, s)$ & $(-, -, n)$ & $(0, 2)$ & $-4$\\ 
$(-, d)$ & $(\pm, \mp, n)$ & $(-, -, n)$ & $(0, 2)$ & $-4$\\ 
$(-, d)$ & $(\pm, \mp, s)$ & $(-, -, n)$ & $(0, 2)$ & $-4$\\ 
\end{tabular}
\caption{\label{orev16} Admissible complex types even, odd, odd without jump}
\end{table}

\begin{table}
\centering
\begin{tabular}{c c c c c c c c c }
$O_1$ & $O_2$ & $O_3$ & $E_0$ & $E_1$ & $E_2$ & $E_3$ & $\mathcal{Z}$ & $\Lambda$\\
$(+, +)$ & $(+, +)$ & $(+, +)$ & $-6$ & $-4$ & $-4$ & $-4$ & $\emptyset$ & $-7$\\
$(+, +)$ & $(+, +)$ & $(\pm, \mp)$ & $-4$ & $-2$ & $-2$ & $-4$ & $\emptyset$ & $-5$\\
$(+, +)$ & $(+, +)$ & $(-, -)$ & $-4$ & $-2$ & $-2$ & $-6$ & $\emptyset$ & $-7$\\
$(+, +)$ & $(\pm, \mp)$ & $(\pm, \mp)$ & $-2$ & $0$ & $-2$ & $-2$ & $(1)$ & $-3$\\
$(+, +)$ & $(\pm, \mp)$ & $(-, -)$ & $-2$ & $0$ & $-2$ & $-4$ & $(1)$ & $-5$\\
$(+, +)$ & $(-, -)$ & $(-, -)$ & $-2$ & $0$ & $-4$ & $-4$ & $(1)$ & $-7$\\
$(\pm, \mp)$ & $(\pm, \mp)$ & $(\pm, \mp)$ & $0$ & $0$ & $0$ & $0$ & $(0, 1, 2, 3)$ & $-1$\\
$(\pm, \mp)$ & $(\pm, \mp)$ & $(-, -)$ & $0$ & $0$ & $0$ & $-2$ & $(0, 1, 2)$ & $-3$\\
$(\pm, \mp)$& $(-, -)$ & $(-, -)$ & $0$ & $0$ & $-2$ & $-2$ & $(0, 1)$ & $-5$\\
$(-, -)$ & $(-,-)$ & $(-, -)$ & $0$ & $-2$ & $-2$ & $-2$ & $(0)$ & $-7$\\
\end{tabular}
\caption{\label{orev17} Admissible short complex schemes odd, odd, odd}
\end{table}

\begin{table}
\centering
\begin{tabular}{c c c c c c c c c }
$O_1$ & $O_2$ & $O_3$ & $E_0$ & $E_1$ & $E_2$ & $E_3$ & $\mathcal{Z}$ & $\Lambda$\\
$(\pm, \mp)$ & $(\pm, \mp)$ & $(\pm, \mp)$ & $0$ & $0$ & $0$ & $0$ & $(0, 1, 2, 3)$ & $-1$\\
$(\pm, \mp)$ & $(\pm, \mp)$ & $(-, -)$ & $0$ & $0$ & $0$ & $-2$ & $(0, 1, 2)$ & $-3$\\
\end{tabular}
\caption{\label{orev18} Short complex schemes odd, odd, odd, cases left}
\end{table}

\section{Proof of the conjecture in the case without jump}
\subsection{Definitions and first steps}
Assume there exists $C_9$ contradicting the conjecture, $C_9$ realizes one of the three complex types in Table~\ref{orev6}. We will find a contradiction using Bezout's theorem with auxiliary rational quartics. The first step consists in studying a pencil of conics based at $A_1, A_2, A_3$ and a fourth point $D$, chosen in a quadrangular oval. It turns out that this pencil has a "double jump": it sweeps out successively five ovals $A, E, B, F, C$, such that $A, B, C$ are triangular, $E, F$ are quadrangular (see Lemma~3 ahead).

\begin{figure}[htbp]
\centering
\includegraphics{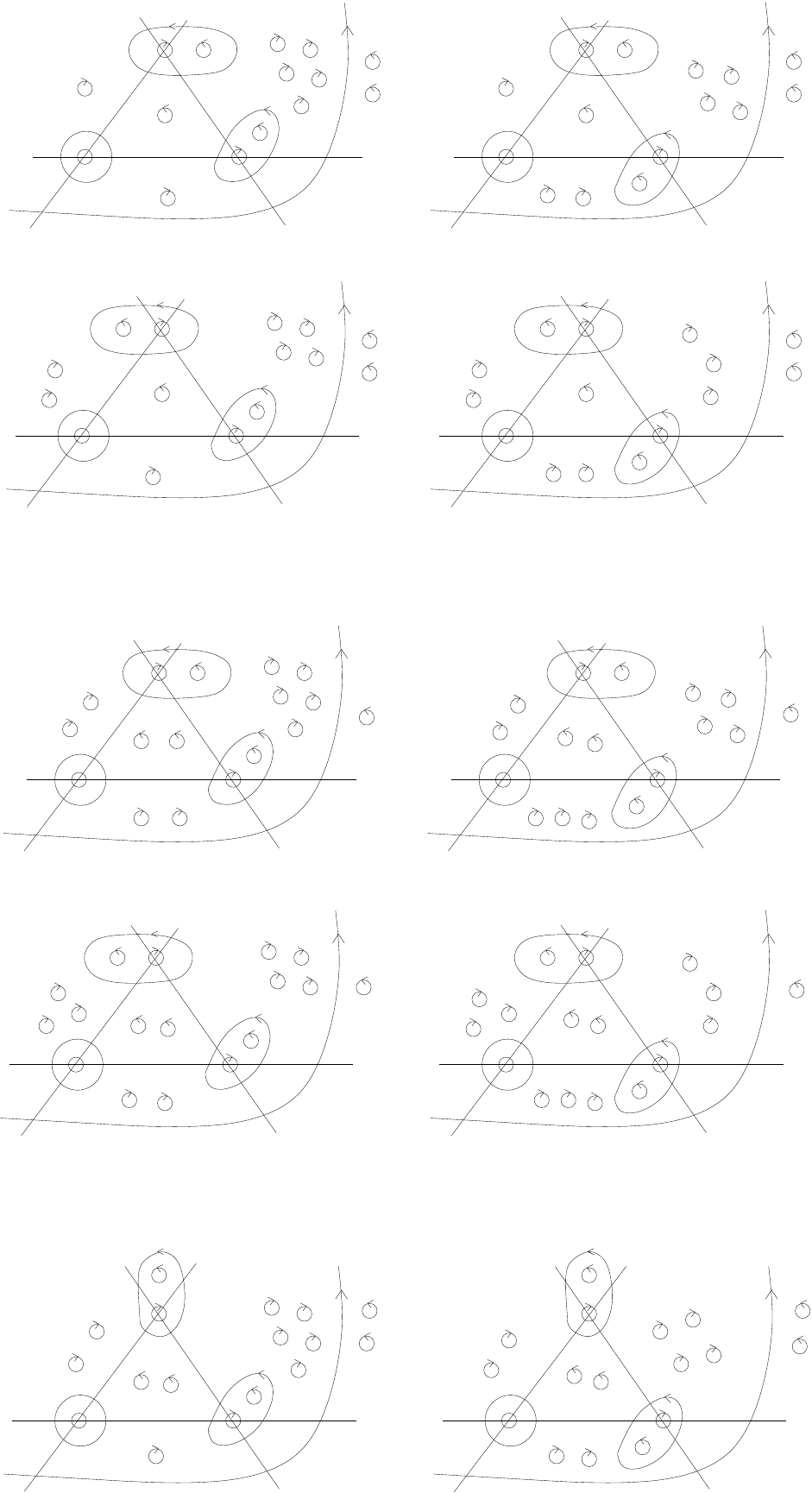}
\caption{\label{types} Admissible distributions of ovals for the curves even, even, odd without jump}
\end{figure}

\begin{figure}[htbp]
\centering
\includegraphics{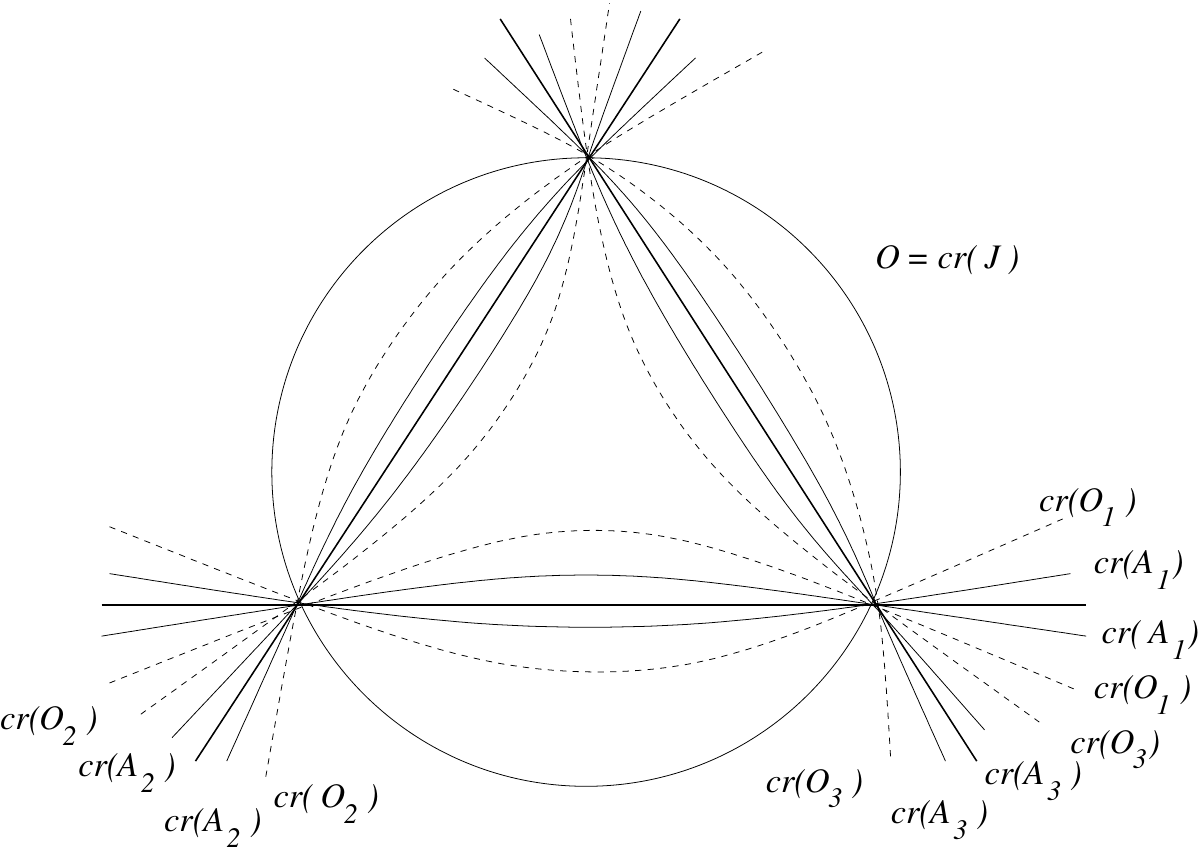}
\caption{\label{conj3} The main part of $C_{18}$ and the base lines}
\end{figure}

For $i = 1, 2$, let $A_i$ be the positive extreme oval in the chain of $O_i$, and let $A_3$ be any one of the extreme ovals in the chain of $O_3$, the pair of ovals $A_3$, $O_3$ have opposite orientations. In Figure~\ref{types}, we show all admissible distributions of ovals for $C_9$. 
The base lines, $\mathcal{J}$ and the ovals $O_i, i = 1, 2, 3$ divide the plane in zones. The ovals that are not represented are distributed in such a way that the difference {\em number of positive minus number of negative ones\/} is $0$ in each zone. Let $O_1$ be non-separating ($C_9$ has the first or the second type), we say that $O_1$ is {\em left\/} if the chain of $O_1$ is in $Q_2$, {\em right\/} if the chain of $O_1$ is in $Q_3$. The oval $O_2$ is non-separating for all three types, we say that $O_2$ is {\em left\/} if the chain of $O_2$ is in $Q_3$, {\em right\/} if the chain of $O_2$ is in $Q_1$.
Perform the Cremona transformation $cr$ based at the points $A_1, A_2, A_3$. The respective images of the lines $A_1A_2$, $A_2A_3$, $A_3A_1$ will be denoted by $A_3$,$A_1$, $A_2$ (note that after $cr$, the new base lines $A_1A_2$, $A_2A_3$, $A_3A_1$ divide the plane in four triangles $cr(T_0)$, $cr(T_i \cup Q_i)$, $i = 1, 2, 3$). For the points lying inside of non-principal ovals, we use the same notation as before $cr$. The curve $C_9$ is mapped onto a curve $C_{18}$ with ordinary $9$-fold singularities at the points $A_i, i = 1, 2, 3$, see Figure~\ref{conj3}.  
An oval of $C_9$ is {\em interior\/} if it lies inside of one of the three nests, otherwise it is {\em exterior\/}. An oval of $C_9$ is {\em triangular (quadrangular)\/} if it lies in $T_0 \cup T_1 \cup T_2 \cup T_3$ ($Q_1 \cup Q_2 \cup Q_3$). An oval of $C_{18}$ is {\em interior (exterior, triangular, quadrangular)\/} if its preimage is. The triangular ovals of $C_{18}$ are the ovals lying inside of $\mathcal{O} = cr(\mathcal{J})$. The {\em main part\/} of $C_{18}$ is the union of $\mathcal{O}$ and the images of the six principal ovals. 
A conic through $A_1, A_2, A_3$ and two further points is mapped by $cr$ onto the line passing through the images of the two points. 
A rational quartic with nodes at $A_1, A_2, A_3$ and passing through five further points is mapped onto the conic passing through the images of these five points.
We shall always choose the two (or the five) points inside of non-principal ovals.
Our goal is to find configurations of five ovals of $C_9$ leading to a contradiction:
after $cr$, a conic passing through the images of these five ovals intersects the curve $C_{18}$ at $38$ points. Its preimage would be a rational quartic intersecting: each of the ovals $A_i, O_i, i = 1, 2, 3$ at four points, each of the other five ovals at two points, and $\mathcal{J}$ at four points.  
A pencil of lines (conics) sweeping out a curve of degree $m$ is {\em maximal\/} if it has alternatively $m$ ($2m$) and $m-2$ ($2m-2$) real intersection points with the curve.  It is {\em totally real\/} if it has only real intersections with the curve. 
The main part of $C_{18}$ divides the plane in $25$ zones. The six zones corresponding to the interiors of the ovals $A_1, A_2, A_3$ contain no ovals. For $i = 1, 2, 3$, at most one of the four zones corresponding to the interior of $O_i$ may contain some ovals.
After $cr$, we will call {\em chain\/} a maximal Fiedler chain of ovals lying in the same zone, swept out by some maximal pencil of lines. The consecutive ovals are connected by vanishing cycles.


Consider two non-principal ovals: $G$ exterior in $Q_1$ and $D$ in $Q_2$, either exterior or extremal negative in the chain of $O_1$ (in the latter case, $O_1$ is left). Note that such ovals exist, see Figure~ \ref{types}.
After $cr$, denote by $\mathcal{F}_D$ the pencil of lines based at $D$ sweeping out $\mathcal{O}$ from the first triangular oval to the last one, and meeting successively $A_1, A_2, A_3$, this pencil is maximal. 
We will say that an oval is interior to $cr(O_i)$ if its preimage is interior to $O_i$.
If $\mathcal{F}_D$ sweeps out some oval interior to $cr(O_i)$, then it sweeps out  all of the ovals interior to $cr(O_i)$, consecutively; the ordering of these ovals given by $\mathcal{F}_D$ corresponds to the ordering of their preimages in the chain of $O_i$. The ovals interior to $cr(O_i)$ (other than $D$) build a chain, formed by the images of the non-principal ovals interior to $O_i$ (other than $D$).
This follows from Bezout's theorem, along with the fact that, in the pseudo-holomorphic setting, one may deform the curve, replacing the vanishing cycles connecting consecutive ovals in the chains by real double points, \cite{or1}-\cite{or7}. 

The triangular ovals forming positive (negative) pairs with $\mathcal{O}$ will be called {\em plus-ovals (minus-ovals)\/}. The plus-ovals (minus-ovals) are thus the negative (positive) ovals in $cr(T_0)$ and the positive (negative) ovals in $cr(T_1) \cup cr(T_2) \cup cr(T_3)$.
For all three complex types of Table~\ref{orev6}, $cr(T_2)$ is empty; for the first two complex types, $cr(T_1)$ is also empty.

\begin{lemma}
There cannot exist a configuration of three ovals $P, H, Q$ met successively by $\mathcal{F}_D$, with $P$ in $cr(T_0)$, $Q$ in $cr(T_3)$ and $H$ quadrangular. 
\end{lemma}

\begin{figure}[htbp]
\centering
\includegraphics{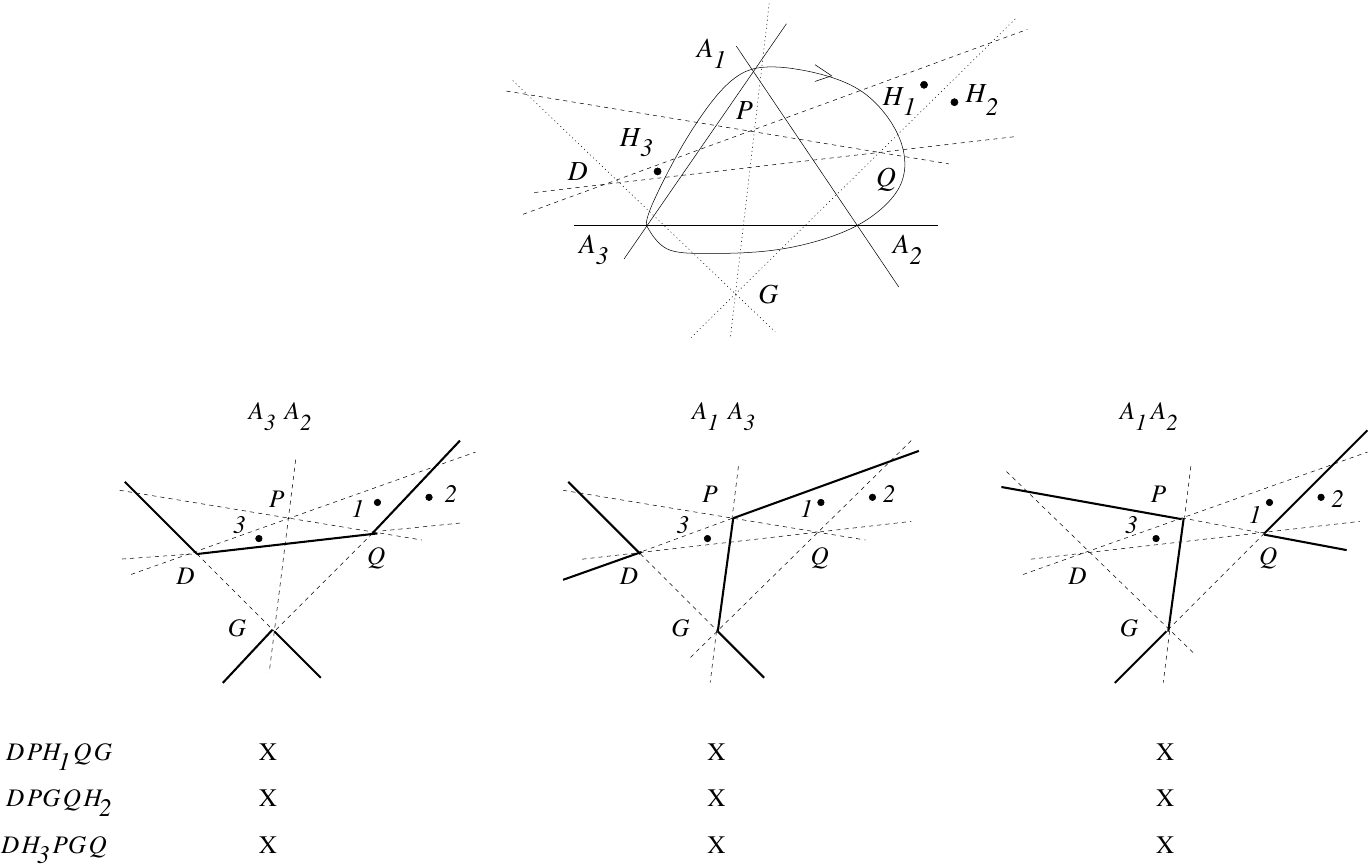}
\caption{\label{pqg} $P$ in $cr(T_0)$, $Q$ in $cr(T_3)$}
\end{figure}

{\em Proof:\/} The lines $PG$, $PQ$ and $QG$ divide the sector $DP$, $DQ$ containing $H$ in four triangles, one of them is entirely contained in $\mathcal{O}$, see Figure~\ref{pqg}. 
The point $H$ may be chosen in any one of the other three triangles, denote by $H_1$, $H_2$ and $H_3$ the three admissible choices of $H$,
the conic through $D, P, Q, G, H$ is $DPH_1QG$, $DPGQH_2$ or $DH_3PGQ$. 
For each of the six lines determined by two of the points $D, P, Q, G$, consider the three intersection points with the base lines, we have thus five particular points on each line. The ordering of these five points on the line is known only for $DP$, $DQ$ and $PQ$. In the upper part of Figure~\ref{pqg}, the position of the base lines with respect to the other three lines $DG$, $PG$, $QG$ has been chosen arbitrarily among all admissible positions. 
In the lower part, we have drawn in bold the boundaries of the zones containing each of the three base lines, each of these zones is a M\"obius band bounded by segments of lines: for $A_3A_2$, it is bounded by $DQ$, $DG$, $QG$, $P$ lies outside;
for $A_1A_3$, it is bounded by $DP$, $PG$, $DG$, $Q$ lies outside;
for $A_1A_2$, it is bounded by $PQ$, $PG$, $QG$, $D$ lies outside.
(In the three pictures, the positions of $H$ are indicated only with the indices $1, 2, 3$.)
The X assigned to a pair (conic, base line $A_iA_j$) means that the conic cuts $A_iA_j$ (twice) and each of the objects $cr(O_k)$, $cr(A_k)$ at four points. We say that the conic is {\em maximal\/} with respect to the base line $A_iA_j$. (We may say similarly that an arc of conic is maximal with respect to some base line.) The three conics are all maximal with respect to the three base lines, and have each four intersection points with $\mathcal{O}$. They cut $C_{18}$ at $38$ points, this is a contradiction. $\Box$ 

\begin{lemma}
There cannot exist a configuration of four ovals $P$, $Q$, $R$, $S$, with $P, R$ in $cr(T_3)$, $Q, S$ in $cr(T_0)$ and $Q$ not in $cr(O_3)$, such that $\mathcal{F}_D$ meets successively $P, Q, R, S$ or $S, P, Q, R$.
\end{lemma}

\begin{figure}[htbp]
\centering
\includegraphics{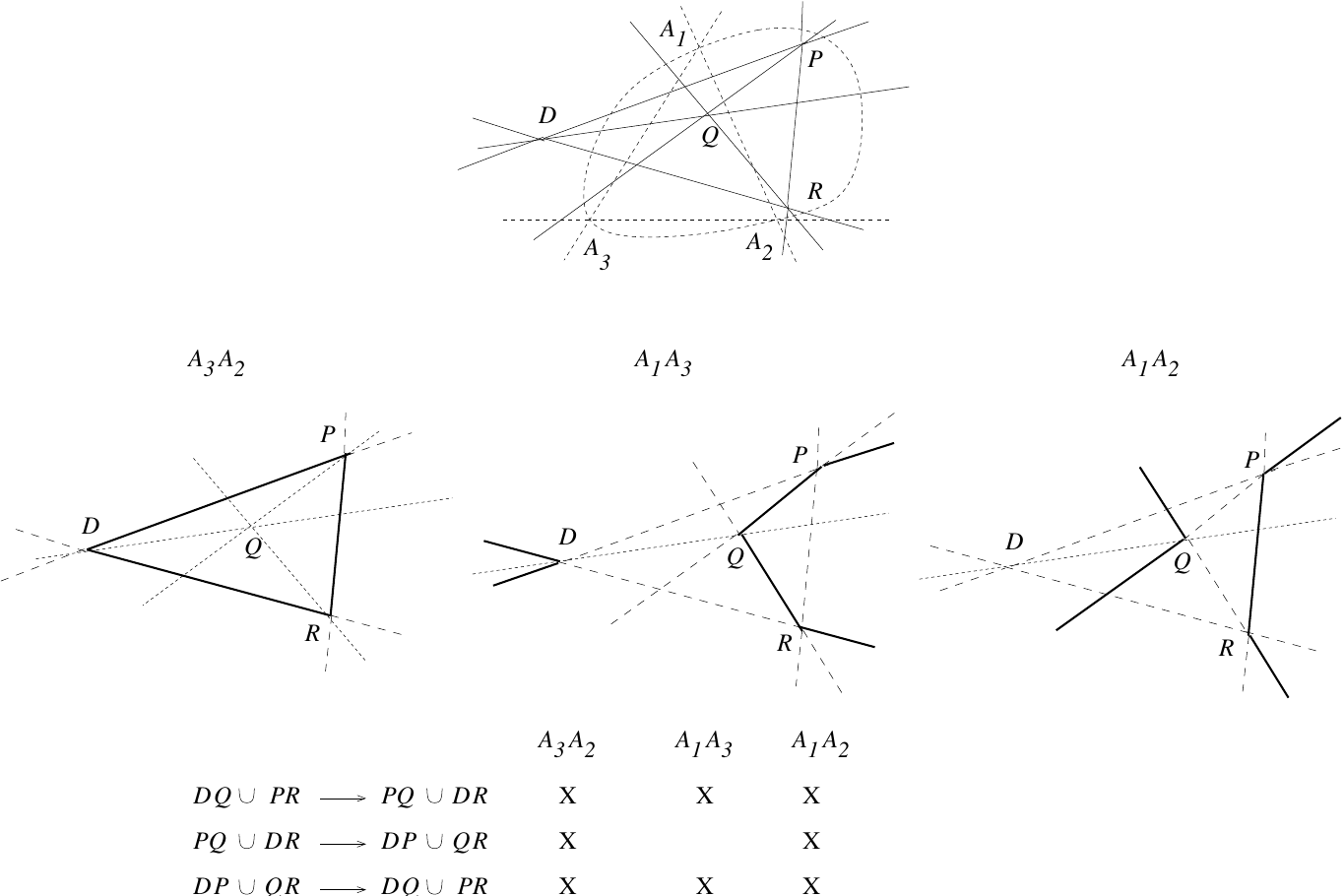}
\caption{\label{dpqrs} $Q$ in $cr(T_0)$ exterior to $cr(O_3)$, $P, R$ in $cr(T_3)$}
\end{figure}

{\em Proof:\/} Assume there exist four such ovals, note that they are either all exterior or one unique oval among $P, R, S$ is interior to $cr(O_3)$. Consider the pencil of conics $\mathcal{F}_{DPQR}$, all of its conics have in total eight intersection points with the set of base ovals. The pencil is divided in three portions by the double lines, in two of these portions, the conics cut $\mathcal{O}$ four times and are maximal with respect to the three base lines. These portions are totally real, see Figure~\ref{dpqrs}, where the $X$ assigned to a pair (portion, base line) means that the portion is maximal with respect to the base line.
The remaining empty ovals are swept out by the third portion $PQ \cup DR \to DP \cup QR$. 
Depending on whether $\mathcal{F}_D$ sweeps out $P, Q, R, S$ or $S, P, Q, R$, the conic of $\mathcal{F}_{DPQR}$ passing through $S$ is $DPQRS$ or $DSPQR$. 
Taking account of the positions of the points in the various zones, and of the positions of the base lines given in the lower part of Figure~\ref{dpqrs}, we get the following arguments.
In the first case, the arc $RSD$ cuts $A_3A_2$ (once), and doesn't cut $A_1A_2$. So it must cut $A_1A_3$ twice, in such a way that it is maximal with respect to this line. In the second case, the arc $DSP$ cuts $A_3A_2$ (once), and doesn't cut $A_1A_2$. So it must cut $A_1A_3$ twice, in such a way that it is maximal with respect to this line. In both cases, the conic is maximal with respect to the three base lines, and cuts $\mathcal{O}$ four times, this is a contradiction. $\Box$

For the first two complex types, the ovals swept out by $\mathcal{F}_D$ are of three {\em kinds\/}: quadrangular, triangular in $cr(T_0)$, triangular in $cr(T_3)$. For the third complex type, we have the same three  kinds, plus supplementarily a chain of ovals in $cr(T_1)$, interior to $cr(O_1)$. A triangular oval is exterior or interior to $cr(O_3)$. Let us call {\em sequence\/} a maximal Fiedler chain of ovals of the same kind, swept out consecutively by $\mathcal{F}_D$. Each sequence is formed of one or several chains.
A pair of odd chains of the same kind is {\em inessential\/} if the pencil sweeps out only even sequences between them.

\begin{figure}[htbp]
\centering
\includegraphics{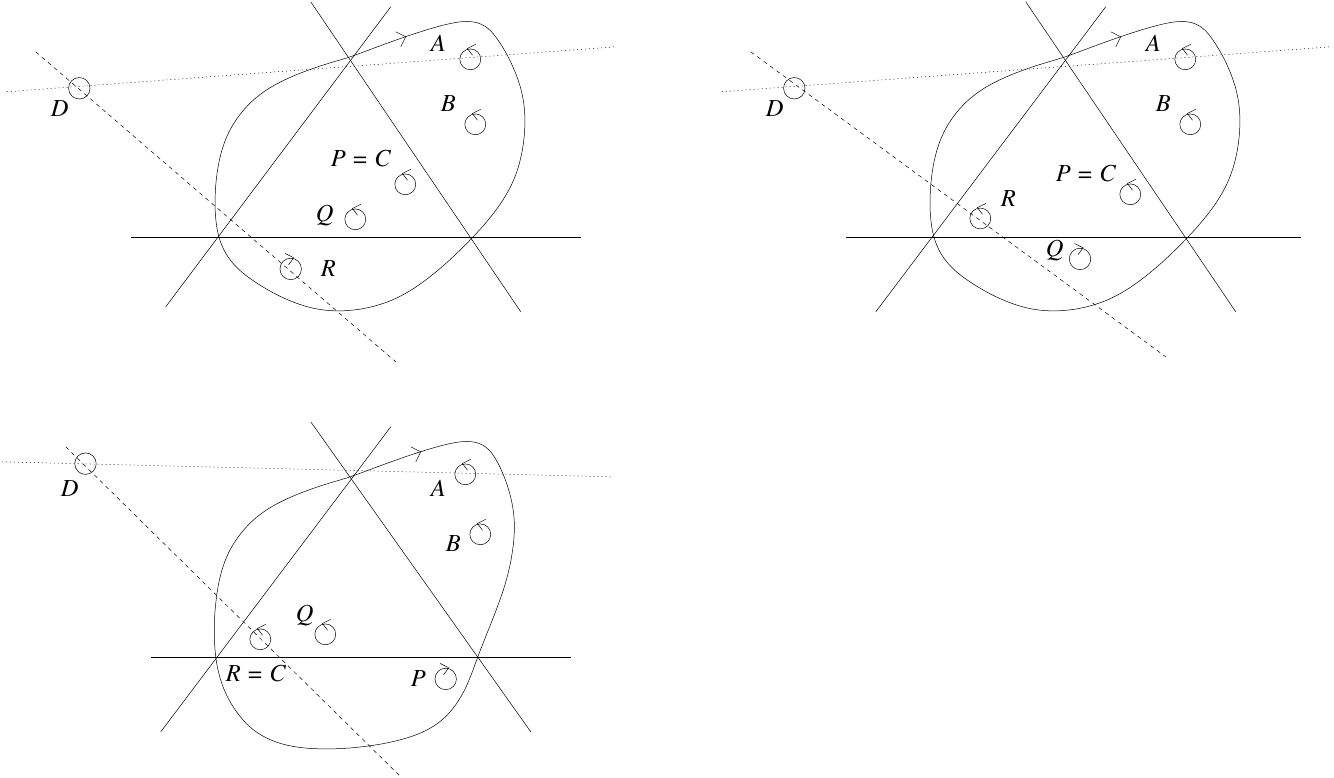}
\caption{\label{abpqr} Third complex type, ovals $A, B, P, Q, R$}
\end{figure}

\begin{lemma}
\begin{enumerate}
\item
There exist three exterior triangular plus-ovals $A, B, C$, swept out successively by $\mathcal{F}_D$, each of them extremity of an odd chain, such that: for the first and third type, $A, B$ are in $cr(T_3)$, $C$ is in $cr(T_0)$; for the second type, $A$ is in $cr(T_3)$, $B, C$ are in $cr(T_0)$.  
For any choice of $A, B, C$, there exist two quadrangular ovals $E, F$, each of them extremity of an odd chain, such that $\mathcal{F}_D$ sweeps out successively
$A, E, B, F, C$ with alternating orientations. 
\item
Let $C_9$ realize the third complex type, there exist five ovals $A, B, P, Q, R$, swept out successively by $\mathcal{F}_D$, extremities of odd chains, such that: $A, B$ are plus-ovals in $cr(T_3)$, two of the ovals $P, Q, R$ are plus-ovals in $cr(T_0)$, the third one is a minus-oval in $cr(T_1)$, interior to $cr(O_1)$. If the minus-oval is $Q$ or $R$, let $C = P$, if the minus-oval is $P$, let $C = R$, see Figure~\ref{abpqr}. For any choice of $A, B, P, Q, R$, there exist furthermore two quadrangular ovals $E, F$, each of them extremity of an odd chain, such that the pencil sweeps out successively: $A, E, B, F, P=C, Q, R$ or $A, E, B, P, Q, F, R=C$ with alternating orientations.
\item
Let $C_9$ realize the first type or the third type with $P = C$, one may choose $B$ and $C$ such that the ovals in each zone $cr(T_0)$ and $cr(T_3)$ swept out beween $B$ and $C$ are distributed in even chains and disjoint pairs of inessential odd chains. 
\end{enumerate}
\end{lemma}

{\em Proof:\/} 
Note for the third type that 2. implies 1.  
Let us "delete" all of the even sequences, and systematically disjoint inessential pairs of odd chains (there are several ways to do that.) In each odd sequence, delete all of the even chains, and then all of the ovals left but one extremity.
Now, the pencil $\mathcal{F}_D$ meets successively single ovals, such that any two consecutive ones are of different kinds, and each of them was originally extremity of an odd chain. Note that the deletion procedure has preserved the triangular parameters $\lambda_i, i = 0, 4, 5, 6$.
Assume that $\mathcal{F}_D$ meets successively $M$ in $cr(T_0)$ and $N$ in $cr(T_3)$. (Note that for the third type, the ovals in $cr(T_1)$ will be swept out after $N$.) By Lemma~1, $\mathcal{F}_D: M \to N$ meets no quadrangular ovals. So, this pencil must sweep out alternatively ovals in $cr(T_3)$ and $cr(T_0)$ (the first one being in $cr(T_3)$. By Lemma~2, there are no ovals at all between $M$ and $N$. 
With the non-zero data $\lambda_0$, $\lambda_6$ (and $\lambda_4$ for the third type), we deduce immediately the following property. For the first (second) type, there exist three plus-oval $A, B, C$ swept out successively, such that $A, B$ are in $cr(T_3)$, $C$ is in $cr(T_0)$ ($A$ is in $cr(T_3)$, $B, C$ are in $cr(T_0)$). One gets also the existence of $P, Q, R$ for the third type.
The existence of ovals $E$ and $F$ as required is also clear by Fiedler's theorem, 
we need yet to prove that they are quadrangular. 

\begin{figure}[htbp]
\centering
\includegraphics{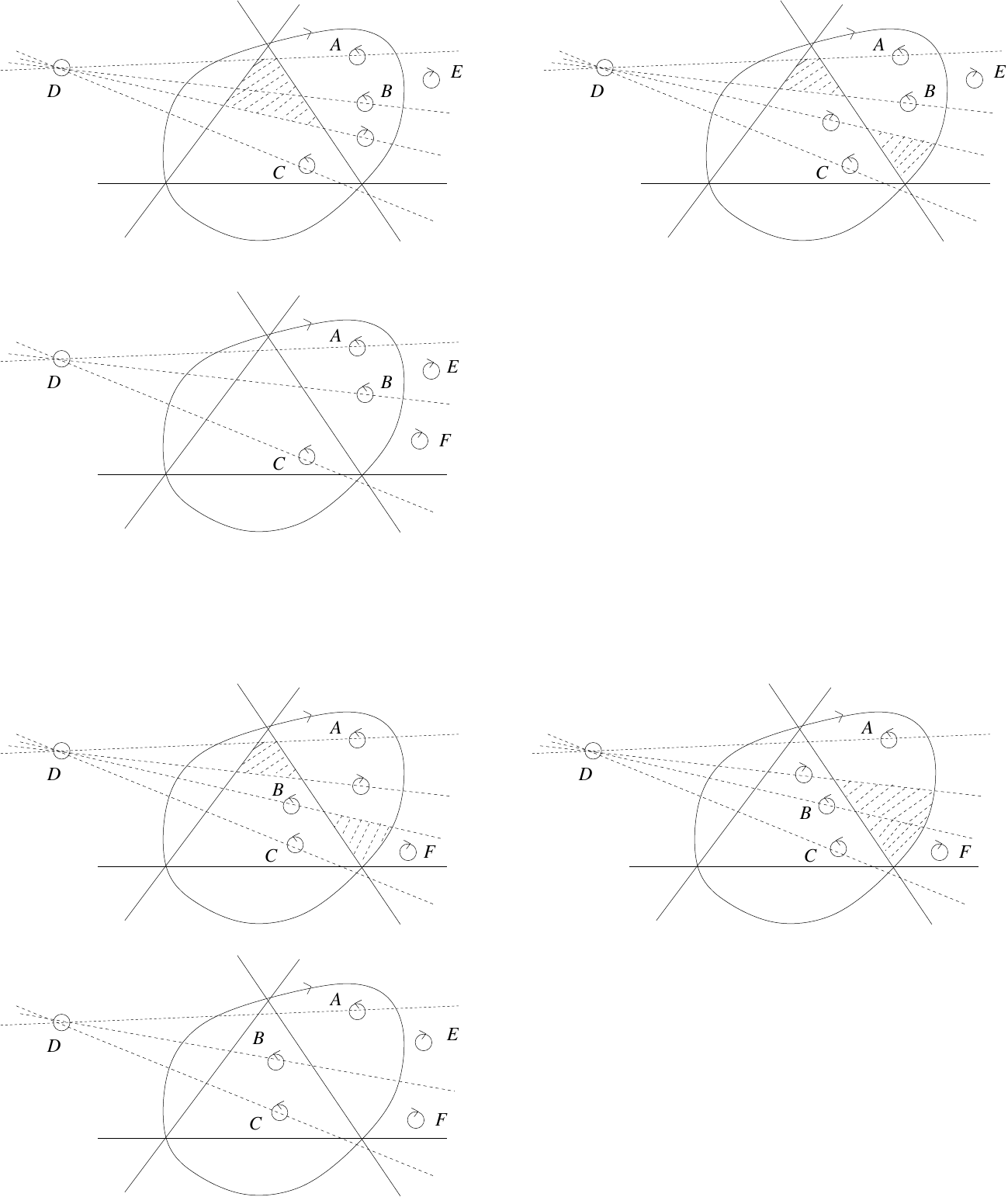}
\caption{\label{efquad} Quadrangular ovals $E$, $F$}
\end{figure}

After the deletion procedure, the following obvious property holds:
{\em Consider two ovals $M, N$ in the same zone $cr(T_i)$ ($i = 0$ or $3$), such that $\mathcal{F}_D$ meets no oval between them in $cr(T_{(3-i)})$. Then $M$ and $N$ have the same orientation and there is an odd number of ovals between them, alternatively quadrangular and in $cr(T_i)$ (the first and the last are quadrangular).\/} 
Let $C_9$ realize the first complex type, see upper part of Figure~\ref{efquad} (the quadrangular ovals have been placed arbitrarily in $cr(Q_3)$). By Lemma~2, there is no oval in $cr(T_0)$ between $A$ and $B$, hence $E$ is quadrangular. 
If $F$ is in $cr(T_3)$, then Lemma~2 implies that there is no oval in $cr(T_0)$ between $B$ and $F$, contradiction. If $F$ is in $cr(T_0)$, then Lemma~2 implies that there is no oval in $cr(T_3)$ between $F$ and $C$, contradiction again. 
So, $F$ is quadrangular, point 1. is proved.
One may choose $B$ and $C$ such that $B$ is the last plus-ovals before $F$ in $cr(T_3)$ and $C$ is the first plus-oval after $F$ in $cr(T_0)$. We have already proved that there is no minus-oval between $B$ and $C$ in $cr(T_0) \cup cr(T_3)$. Consider now again the whole set of ovals. The pencil $\mathcal{F}_D$ meets between $B$ and $C$ in $cr(T_0) \cup cr(T_3)$ only even chains and disjoint pairs of inessential odd chains, point 3. is proved.

Caution: the argument doesn't work if, instead of deleting ovals, we had chosen $B$ extremity of the last odd plus-chain in $cr(T_3)$ and $C$ extremity of the first odd plus-chain in $cr(T_0)$. The pencil $\mathcal{F}_D$ could sweep out for example: $B$, an even chain in $cr(Q_3)$, an odd chain in $cr(T_3)$ starting with a minus-oval $B'$ and $C$ (the odd chains containing $B$ and $B'$ form an inessential pair). In this case, we cannot find $F$ quadrangular as required.

Let $C_9$ have the second type, see lower part of Figure~\ref{efquad}. By Lemma~2, there is no oval in $cr(T_3)$ between $B$ and $C$, so $F$ is quadrangular.
If $E$ is in $cr(T_3)$, Lemma~2 implies that there is no oval in $cr(T_0)$ between
$A$ and $E$, contradiction. If $E$ is in $cr(T_0)$, Lemma~2 implies that there is no oval in $cr(T_3)$ between $E$ and $B$, contradiction again. Hence, $E$ is quadrangular.
Let $C_9$ have the third type. For $R = C$, there must be some oval $F$ between $Q$ and $R$ by Fiedler's theorem, and $F$ can be only quadrangular.
Apart from that, the same arguments as for the first type apply to prove that $E$ and $F$ are quadrangular. For $P = C$, the same arguments as for the first type apply to prove point 3.
$\Box$

\subsection{The first and the third complex type}
In the next Lemmas~4-7, $D$ and $\mathcal{F}_D$ are as in the previous section, $A, E, B, F, C$ are five ovals swept out successively by $\mathcal{F}_D$, such that
$A, B$ are exterior in $cr(T_3)$, $C$ is exterior in $cr(T_0)$, $E, F$ are quadrangular. (These conditions are weaker than those in Lemma 3.)
Denote by $\mathcal{S}$ the sector swept out by the piece of pencil $\mathcal{F}_D: B \to F \to C$.
Let $[DC]$ be the segment of line $DC$ that does not cut $A_1A_3$, and $[DC]'$ be the other segment. 
We shall study the admissible distributions of the ovals in the various zones.
Keep in mind that an oval in $cr(Q_i)$ may be exterior, interior to $cr(O_j)$ or interior to $cr(O_k)$. If $O_1$ is left (right), the ovals interior to $cr(O_1)$ are in
$cr(Q_2)$ ($cr(Q_3)$). If $O_1$ is separating, the ovals interior to $cr(O_1)$ are in $cr(T_1)$.

\begin{lemma}
The intersection point $AB \cap DC$ lies on $[DC]$.
\end{lemma}

\begin{figure}[htbp]
\centering
\includegraphics{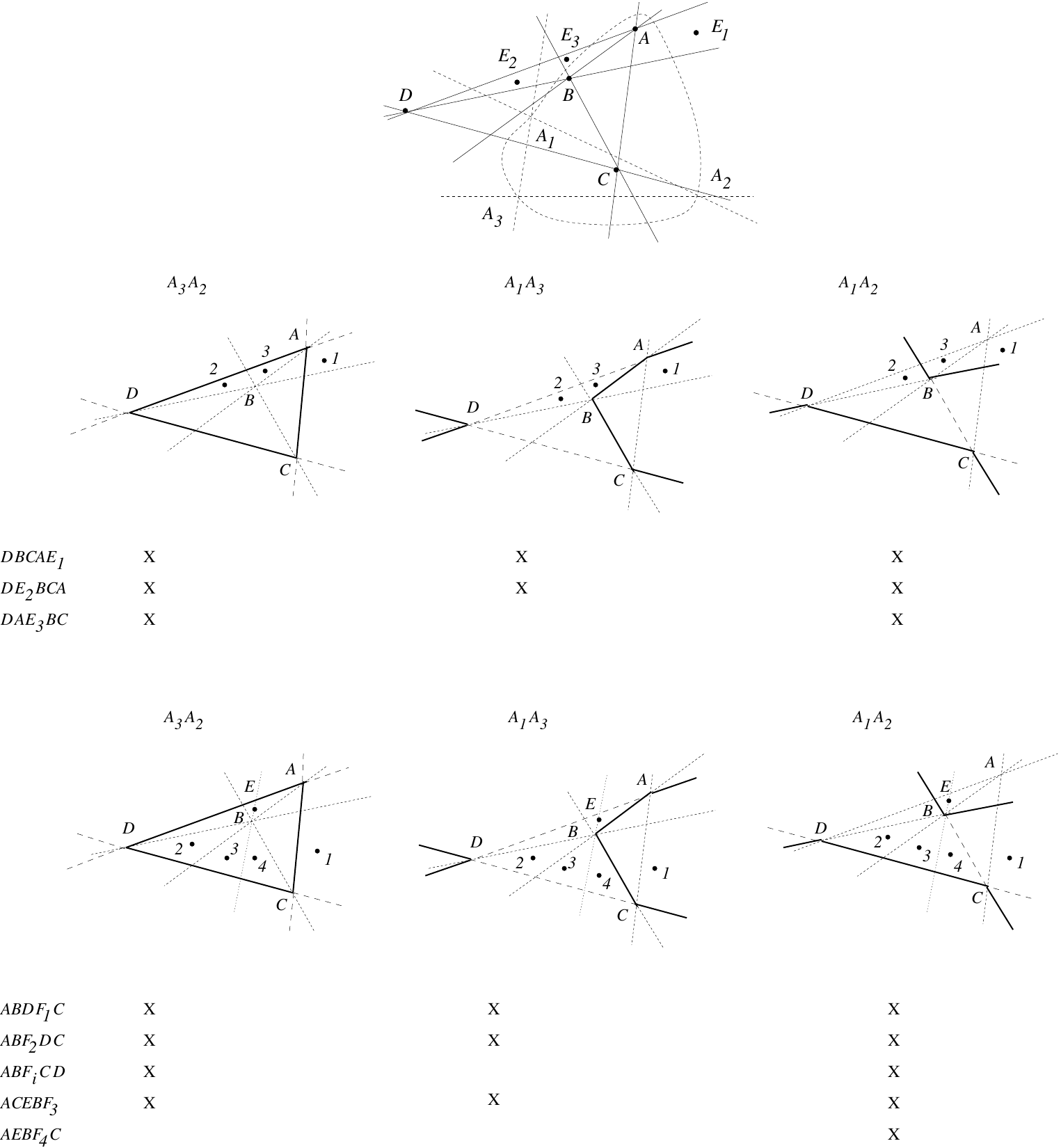}
\caption{\label{abfcd} $AB \cap DC$ on $[DC]'$, position of $E$ and $F$, conics}
\end{figure}

\begin{figure}[htbp]
\centering
\includegraphics{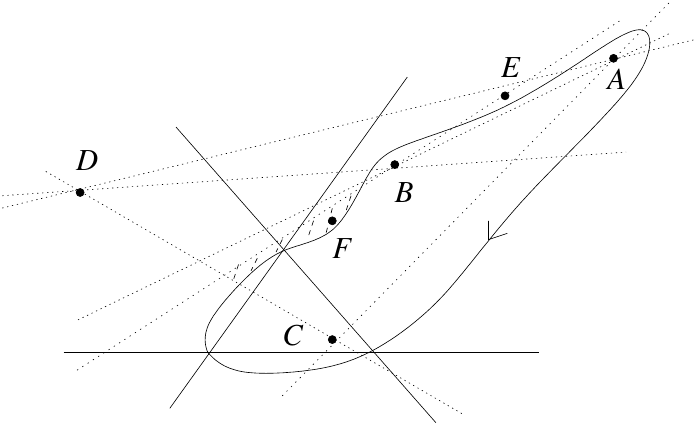}
\caption{\label{case2} $A, E, B, F, C, D$ lie in convex position, contradiction}
\end{figure}

{\em Proof:\/}
Assume that $AB \cap DC$ lies on $[DC]'$.
The lines $AC$, $AB$ and $BC$ divide the sector $DA$, $DB$ containing $E$ in four triangles, one of these triangles is entirely contained in $\mathcal{O}$, see Figure~\ref{abfcd}.
Denote by $E_1, E_2, E_3$ the three choices of $E$. The conic through $D, A, B, C, E$ is $DBCAE_1$, $DE_2BCA$ or $DAE_3BC$, they cut all $\mathcal{O}$ four times, and the first two are maximal with respect to the three base lines. Hence, $E = E_3$, the conic is $DAEBC$. Note that if $A_1A_3$ doesn't cut the triangle containing $E$, then this triangle lies in $cr(T_3 \cup Q_3)$, otherwise, $A_3A_1$ divides this triangle in two pieces, one is in $cr(T_3 \cup Q_3)$, the other in $cr(Q_1)$. As $DAEBC$ cannot be maximal with respect to $A_3A_1$, the oval $E$ is in $cr(Q_3)$, or interior to $cr(O_2)$ in $cr(Q_1)$.
The lines $AB$, $AC$, $BC$, divide the sector $DB$, $DC$ containing $F$ in four triangles, one of them is is entirely contained in $\mathcal{O}$, another is divided in two subtriangles by the line $BE$. Denote by $F_1, F_2, F_3, F_4$ the four admissible choices of $F$, corresponding to $ABDF_1C$, $ABF_2DC$, $ABF_iCD$, $i = 3, 4$, $ACEBF_3$, $AEBF_4C$. All of these conics cut $\mathcal{O}$ four times.
The conics $ABDF_1C$, $ABF_2DC$ and $ACEBF_3$ are maximal with respect to all three base lines, and the configuration of six points $A, E, B, F_4, C, D$ contradicts Lemma~15 from \cite{fi}, see Figure~\ref{case2}. $\Box$ 

\begin{lemma}
The oval $E$ lies on a conic $DAEBC$ and is in $cr(Q_3)$ (interior to $cr(O_1)$ or exterior), or interior to $cr(O_1)$ in $cr(Q_2)$ (in this latter case, $O_1$ is left).
\end{lemma}

\begin{figure}[htbp]
\centering
\includegraphics{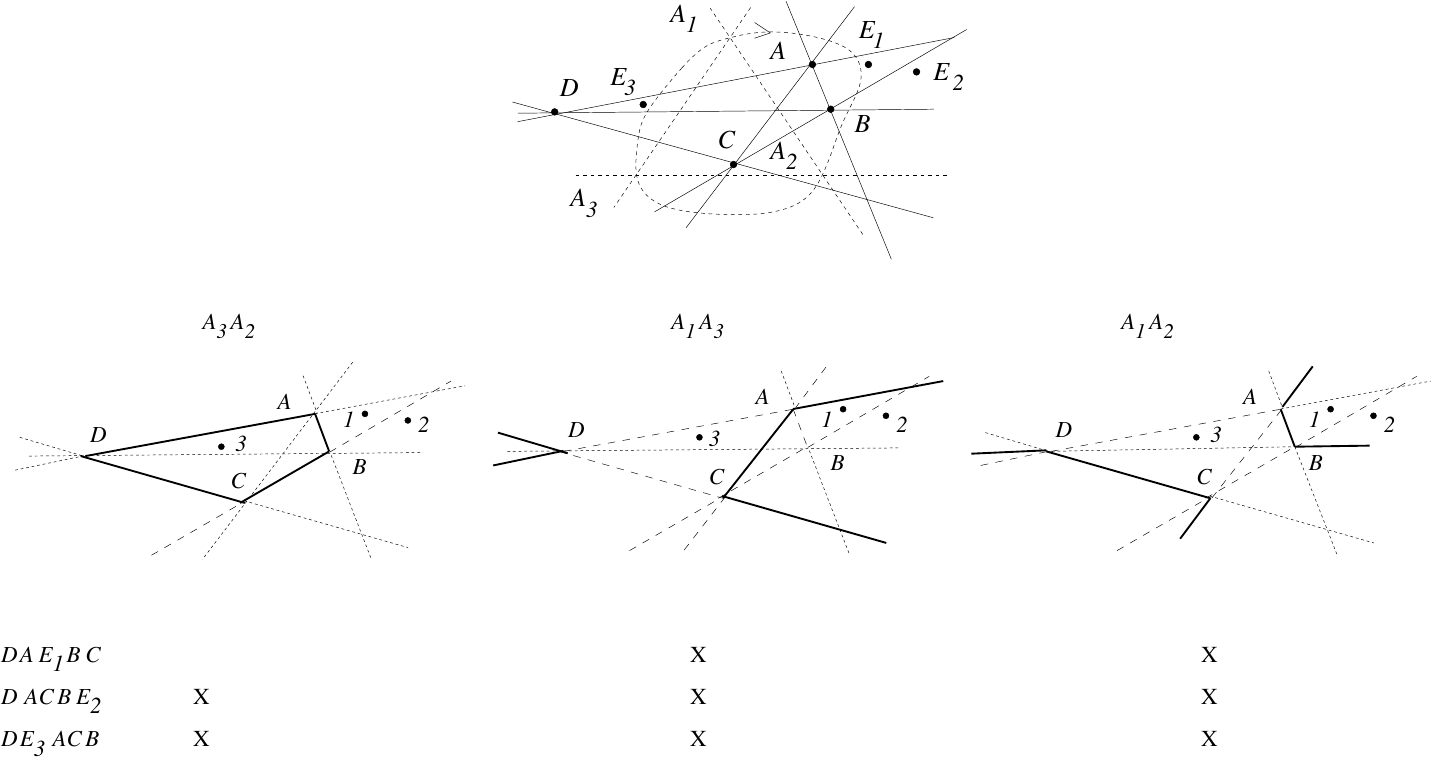}
\caption{\label{daebc} Position of $E$, conics}
\end{figure}

{\em Proof:\/}
The lines $AC$, $AB$ and $BC$ divide the sector $DA$, $DB$ containing $E$ in four triangles, one of these triangles is entirely contained in $\mathcal{O}$, see Figure~\ref{daebc}. The point $E$ may be chosen in any one of the other three triangles, denote by $E_1, E_2, E_3$ the three admissible choices of $E$, the conic through $D, A, B, C, E$ is correspondingly $DAE_1BC$, $DACBE_2$ or $DE_3ACB$,
they cut all $\mathcal{O}$ four times and the last two are maximal with respect to the three base lines.
The conic is $DAEBC$, $E$ lies in a triangle bounded by the lines $AD$, $BC$, $AB$. 
Note that $cr(O_2)$ cuts twice each of the segment $[DA]$, $[DB]$ with non-empty intersection with the M\"obius band corresponding to $A_1A_3$. Bezout's theorem between $C_9$ and the lines $DA$, $DB$ implies thus that the triangle containing $E$ doesn't intersect the interior of $cr(O_2)$.
Note also that if this triangle is not cut by $A_3A_2$, it lies entirely in $cr(T_3 \cup Q_3)$. If the triangle is cut by $A_3A_2$, one of the pieces obtained is in $cr(T_3 \cup Q_3)$, the other is in $cr(Q_2)$. Assume $E$ is in $cr(Q_2)$, as $DAEBC$ cannot be maximal with respect to $A_2A_3$, $E$ is interior to $cr(O_1)$. 
$\Box$

\begin{lemma}
The oval $F$ lies on a conic $AEBCF$.
\end{lemma}

\begin{figure}[htbp]
\centering
\includegraphics{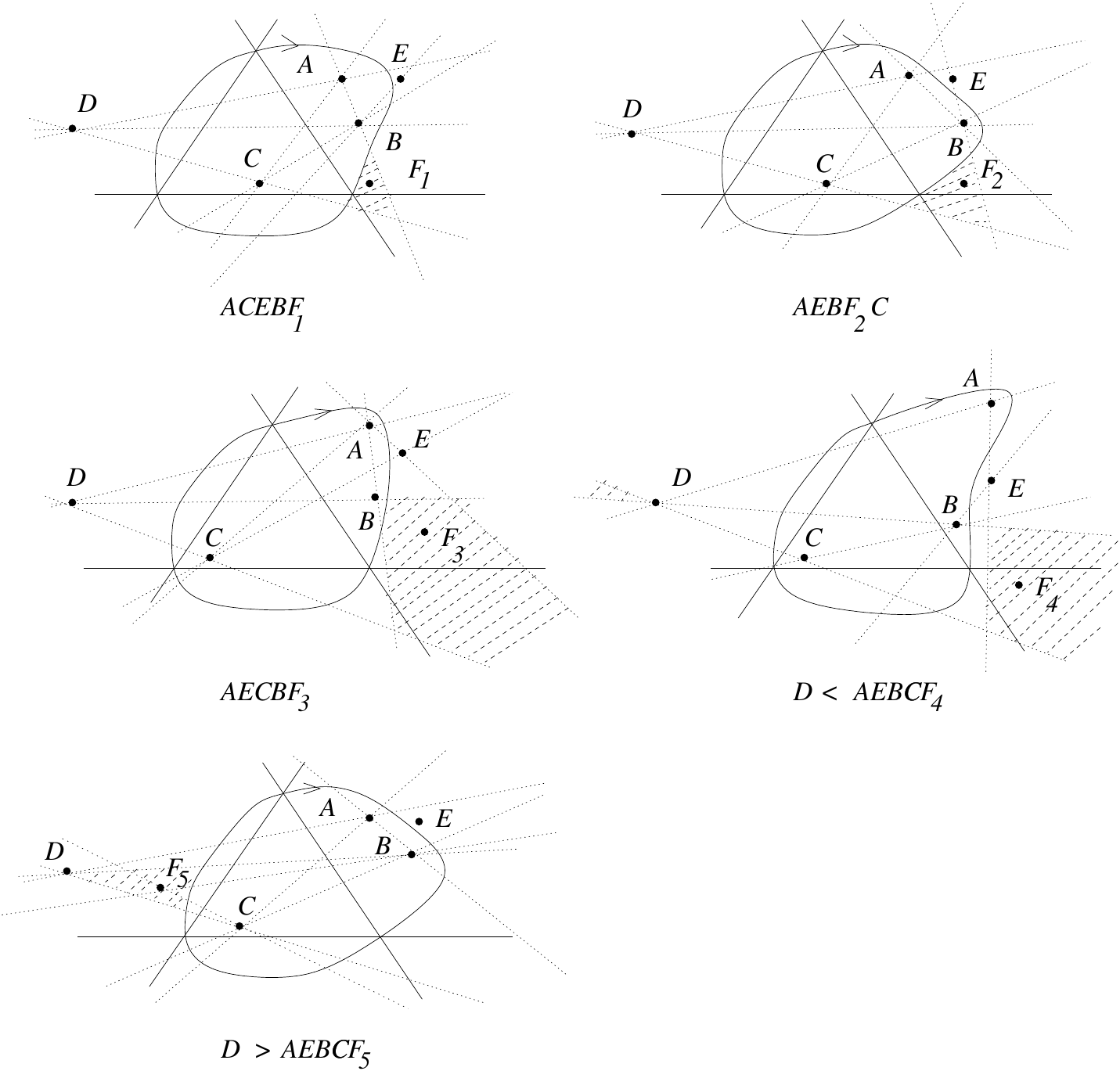}
\caption{\label{case1} The five admissible positions of $F$}
\end{figure}

\begin{figure}[htbp]
\centering
\includegraphics{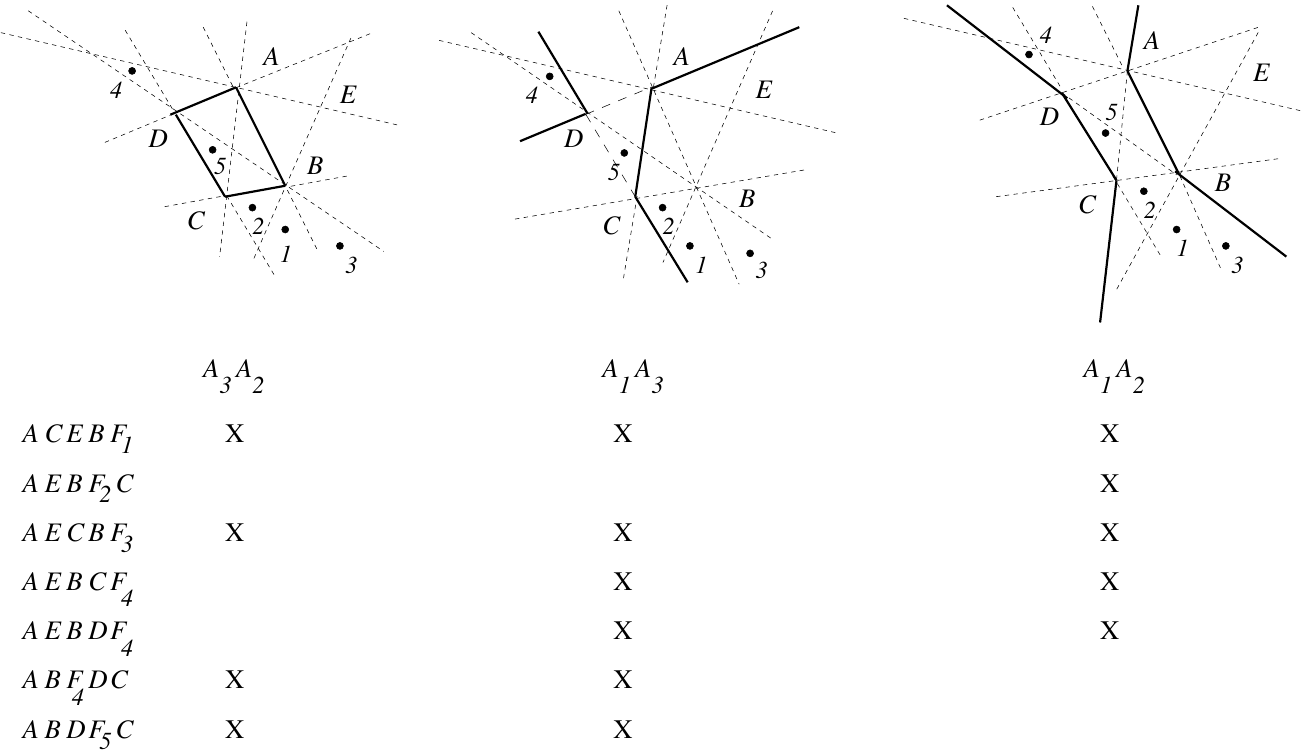}
\caption{\label{acebf} Position of the base lines and of $F$, conics}
\end{figure}

{\em Proof:\/}
The six lines determined by $A, B, C, E$ divide the sector $DB, DC$ containing $F$ in seven zones, two of them lie entirely inside of $\mathcal{O}$. Denote by $F_1, \dots F_5$ the remaining five admissible positions of $F$. Each of them is characterized by the conic through $A, B, C, E, F$, plus (cases $4, 5$) the position of $D$ inside ($<$) or outside ($>$) of this conic, see Figure~\ref{case1}, where $F$ lies in the hatched zones and the position of the base lines with respect to the dotted lines has been chosen arbitrarily among all possible positions. All of the conics intersect $\mathcal{O}$ at four points. The zones containing the base lines are M\"obius bands bounded by segments of lines, they are shown in Figure~\ref{acebf}, along with all five admissible positions of $F$ (indicated with help of the indices $1, \dots 5$). The conics $ACEBF_1$ and $AECBF_3$ are maximal with respect to all three base lines, contradiction. The configuration of six points $A, E, B, F_2, C, D$ contradicts Lemma~15 from \cite{fi}, see second picture in Figure~\ref{case1}. $\Box$

If an oval $K$ swept out between $B$ and $C$ is such that $D < AEBCK$ ($D > AEBCK$), we say that $K$ is in the zone $F_4$ ($F_5$), see the last two pictures in Figure~\ref{case1}. Note that the zone $F_4$ is a triangle bounded by the lines $DB$, $DC$, $AE$. In the case $K = F$, we write shortly $F = F_4$ ($F = F_5$).
From now on, we assume that if $O_1$ is left, then the oval $D$ of $C_9$ is the negative extremity of the chain of $O_1$.

\begin{lemma}
\begin{enumerate}
\item
Let $O_1$ be left. The images by $cr$ of the non-extreme ovals interior to $O_1$ are all together either in the zone $F_4$, or outside of the sector $\mathcal{S}$. The zone $F_4$ contains only these ovals, or is empty. If $F = F_5$, then $F$ may be: exterior or interior to $cr(O_3)$ in $cr(Q_2)$; interior to $cr(O_3)$ in $cr(Q_1)$.
\item
Let $O_1$ be separating. $F = F_4$ may be: exterior or interior to $cr(O_3)$ in $cr(Q_2)$; interior to $cr(O_3)$ in $cr(Q_1)$. $F = F_5$ may be: exterior or interior to $cr(O_3)$ in $cr(Q_2)$; interior to $cr(O_3)$ in $cr(Q_1)$.
\item
Let $O_1$ be right. The admissible positions of $F$ are the same as for separating $O_1$, plus one: $F = F_4$ interior to $cr(O_1)$ in $cr(Q_3)$
\end{enumerate}
\end{lemma}

{\em Proof:\/} 
Let $F = F_4$ and denote by $[FD]$ the segment $FD$ contained in the zone $F_4$, $[FD]$ doesn't cut $A_1A_3$, and is interior to both conics $AEBDF$ and $ABFDC$. Each of these conics cuts $\mathcal{O}$ four times, and is maximal with respect to two base lines, see Figure~\ref{acebf}. 
If $F$ is in $cr(Q_3)$, the segment $[FD]$ cuts $A_3A_2$. If moreover $F$ is exterior to $cr(O_1)$, the conic $AEBDF$ is maximal with respect to $A_3A_2$. But
this conic is also maximal with respect to the other two base lines, contradiction.
If $F$ is in $cr(Q_1)$, the segment $[FD]$ cuts $A_1A_2$. If moreover $F$ is exterior to $cr(O_3)$, the conic $ABFDC$ is maximal with respect to $A_1A_2$. But
this conic is also maximal with respect to the other two base lines, contradiction. 
Hence, $F$ is in $cr(Q_2)$; interior to $cr(O_3)$ in $cr(Q_1)$; or interior to $cr(O_1)$ in $cr(Q_3)$. 

Let $F = F_5$, the conic $ABDFC$ cuts $\mathcal{O}$ four times and is maximal with respect to $A_3A_2$ and $A_1A_3$. The arc $DFC$ doesn't cut $A_3A_2$, and cuts $A_1A_3$ once, see Figure~\ref{acebf}. If $F$ is in $cr(Q_3)$ or in $cr(Q_1)$ exterior to $cr(O_3)$, then $DFC$ is maximal with respect to $A_1A_2$, contradiction. Hence $F$ is in $cr(Q_2)$; or interior to $cr(O_3)$ in $cr(Q_1)$. 

\begin{figure}[htbp]
\centering
\includegraphics{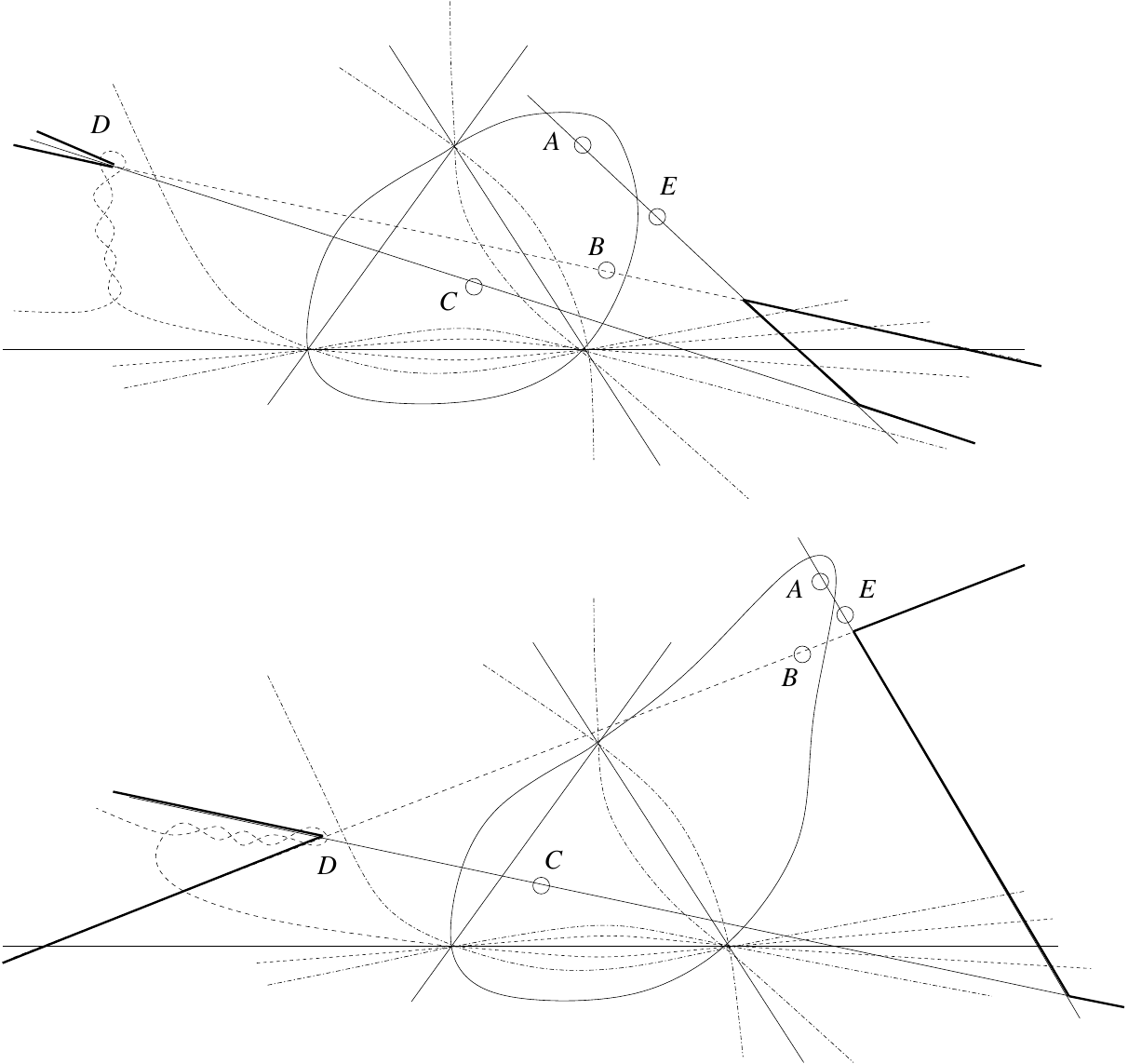}
\caption{\label{dfk} $O_1$ left, two positions of the chain in $cr(O_1)$: exterior to $\mathcal{S}$ and in the zone $F_4$}
\end{figure}

Note that if $D$ is exterior to $cr(O_1)$ ($O_1$ is right or separating), then the zone $F_5$ doesn't intersect the interior of $cr(O_1)$.
Points 2. and 3. follow immediately. To get 1., note that by Bezout's theorem with $cr(O_1)$, the line $DC$ meets successively: $D$, $A_1A_3$, $C$, $A_1A_2$, $A_2A_3$. The zone $F_4$ splits in two parts, one in $cr(Q_3)$, the other interior to $cr(O_1)$ in $cr(Q_2)$, see Figure~\ref{dfk} where the boundary of the zone $F_4$ is drawn in bold, and the vanishing cycles connecting the images of the ovals interior to $O_1$ have been replaced by real double points. 
$\Box$

\begin{lemma}
Consider five ovals $A, E, B, F, C$ satisfying with $D$ the conditions of Lemma 3.
\begin{enumerate}
\item
Let $O_1$ be left. Then $F = F_5$, $F$ is positive exterior in $cr(Q_2)$.
\item
Let $O_1$ be separating or right. If $F = F_4$, then $F$ is negative exterior in $cr(Q_2)$. If $F = F_5$, then $F$ is positve exterior in $cr(Q_2)$. 
\end{enumerate}
\end{lemma}

\begin{figure}[htbp]
\centering
\includegraphics{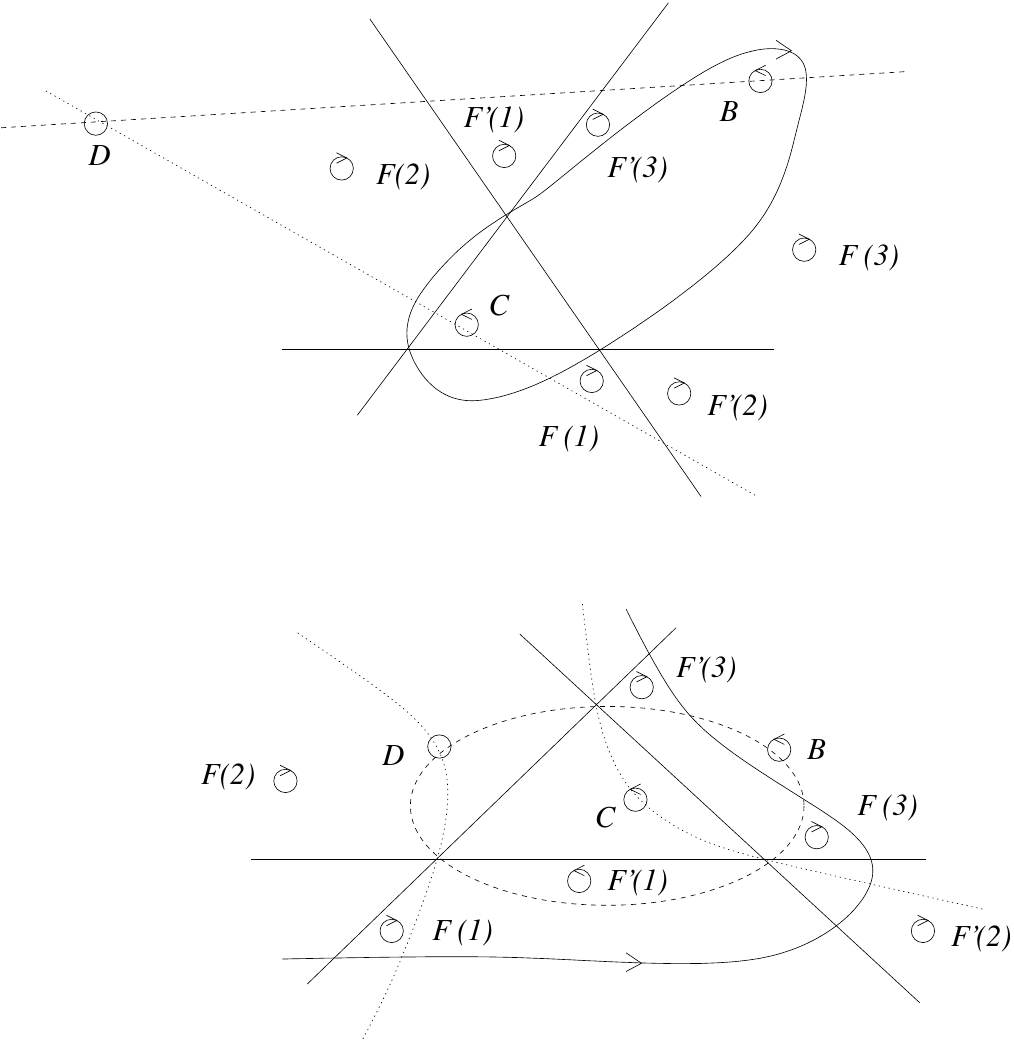}
\caption{\label{fq3} $F(1)$, $F(2)$, $F(3)$, $F'(3)$ are positive, $F'(1)$, $F'(2)$ are negative}
\end{figure}

{\em Proof:\/}
Consider the intersections of the sector $\mathcal{S}$ with the images of the quadrangles. We see that $\mathcal{S} \cap cr(Q_2)$ has two connected components, $\mathcal{S} \cap cr(Q_1)$ and $\mathcal{S} \cap cr(Q_3)$ have each at most two connected components, see Figure~\ref{fq3}, where we have represented the set of objects under consideration in the upper part, and their preimages in the lower part. If $F$ is in $cr(Q_1)$ or $cr(Q_2)$, it is positive ($F(1)$, $F(2)$) or negative ($F'(1)$, $F'(2)$) depending on whether it lies in one connected component or in the other. If $F$ is in $cr(Q_3)$, it is positive for both choices ($F(3)$, $F'(3)$) of its position. 
The case $F = F_4$ corresponds to the three positions $F(1)$, $F(3)$, $F'(2)$, so if $F$ is in $cr(Q_1) \cup cr(Q_3)$, $F$ is positive, if $F$ is in $cr(Q_2)$, $F$ is negative. The case $F = F_5$ corresponds to the other three positions $F'(1)$, $F(2)$, $F'(3)$, so if $F$ is in $cr(Q_2) \cup cr(Q_3)$, $F$ is positive, if
$F$ is in $cr(Q_1)$, $F$ is negative.
The ovals interior to $cr(O_3)$ form an even chain with respect to $D$, so $F$ is not interior to $cr(O_3)$. 
The ovals interior to $O_i, i = 1$ or $2$ bring a contribution of $0$ or $-1$ to
$\lambda_3$. So, if $F$ is in $cr(Q_3)$ or $F$ is in $cr(Q_2)$ with the position $F(2)$, then $F$ is exterior. $\Box$

\begin{lemma}
Let $O_1$ be right or separating.
There exist a base oval $D$ and five ovals $A, E, B, F, C$ (seven ovals $A, E, B, F, P, Q, R$ for the third type) satisfying the conditions of Lemma~3, such that 
the zone $F_4$ is empty or contains only the even chain interior to $cr(O_3)$ (so, $F = F_5$).
\end{lemma}

\begin{figure}[htbp]
\centering
\includegraphics{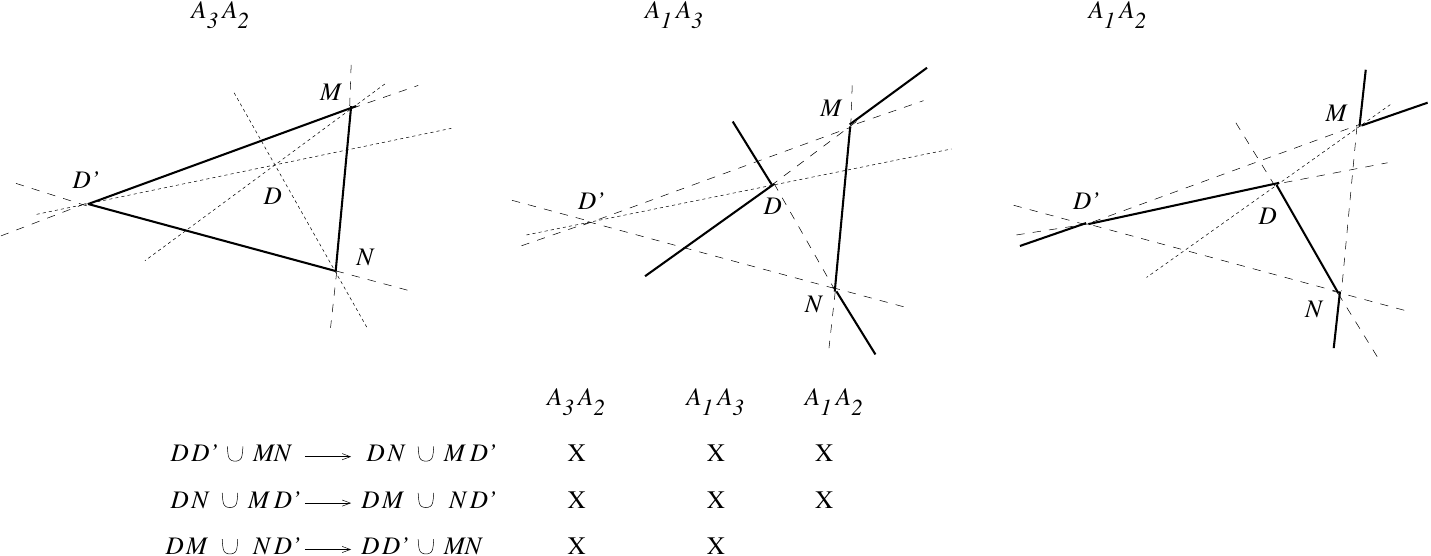}
\caption{\label{emptyf4} The pencil of conics $\mathcal{F}_{DD'MN}$}
\end{figure}

{\em Proof:\/}
Start with a set of five (seven) ovals satisfying the conditions of Lemma~3. 
Assume there exists $D'$ exterior in $cr(Q_2)$, lying in the zone $F_4$. 
Let $\mathcal{F}_{D'}$ be the pencil based at $D'$ and sweeping out $\mathcal{O}$
from the first triangular oval to the last one, in such a way that it meets successively $A_1, A_2, A_3$. 
The pencil $\mathcal{F}_{D'}$ meets successively $A, E, B, \{ D, F \}, C$. 
Replacing the vanishing cycles in the chains by double real points, we see that $\mathcal{F}_D$ and $\mathcal{F}_{D'}$ give rise to the same sets of
triangular chains, and of quadrangular chains that do not contain $D$ nor $D'$.
Using conics, we can get more information.
Let $M$ and $N$ be two triangular oval such that $\mathcal{F}_D$ meets successively $M, D', N$, then $\mathcal{F}_{D'}$ meets successively $M, D, N$.
The pencil of conics $\mathcal{F}_{DD'MN}$ has only one non-totally real portion, see Figure~\ref{emptyf4} . 
Any supplementary oval $S$ lies on a conic: $MSD'DN$, $SMD'DN$, $MD'SDN$, $MD'DSN$. The pencils $\mathcal{F}_D$ and $\mathcal{F}_{D'}$ sweep out the 
set of three ovals $\{ M, N, S \}$ with the same ordering.
Let $M_0$ and $N_0$ be the first and the last (triangular) oval met by $\mathcal{F}_D$, they are also the first and the last oval met by $\mathcal{F}_{D'}$. The two pencils sweep out the same set of ovals. Let $X, Y$ be two ovals met successively by $\mathcal{F}_D$, one of them triangular, then $X, Y$ are met successively by $\mathcal{F}_{D'}$ (if $X$ is triangular, make $(M, S, N) = (X, Y, N_0)$, if $Y$ is triangular, make $(M, S, N) = (M_0, X, Y)$). 
In particular, the orderings with which the triangular chains are swept out by $\mathcal{F}_D$ and $\mathcal{F}_{D'}$ coincide. Let $c_1, c_2$ be two such chains, $\mathcal{F}_D$ and $\mathcal{F}_{D'}$ sweep out the same set of quadrangular ovals between $c_1$ and $c_2$.

\begin{figure}[htbp]
\centering
\includegraphics{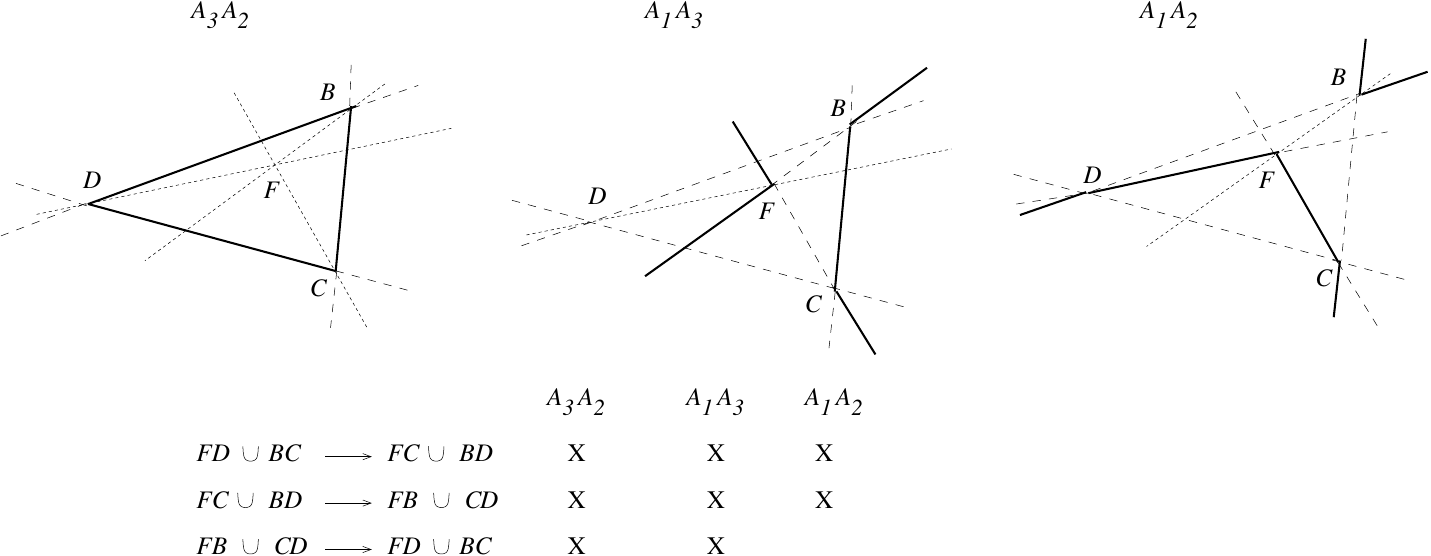}
\caption{\label{pendfcb} The pencil of conics $\mathcal{F}_{DFCB}$}
\end{figure}

We rename the ovals so that $D'$ is now called $D$, the new zone $F_4$ is contained in the old one. Repeat the procedure as many times as necessary, in the end we have a new base oval $D$ such that $A, E, B, C$ ($A, E, B, P, Q, R$) are swept out successively by $\mathcal{F}_D$, these ovals are extremities of odd chains. The ovals in $cr(Q_2)$ swept out between $B$ and $C$ are all in the new zone $F_5$. Finally, the pencil meets between $B$ and $C$ in $cr(T_0) \cup cr(T_3)$ only even chains and disjoint pairs of inessential odd chains.  We may choose now a new oval $F$ to fulfill the conditions of Lemma~3. The pencil of conics $\mathcal{F}_{DFCB}$ is shown in Figure~\ref{pendfcb}.

\begin{figure}[htbp]
\centering
\includegraphics{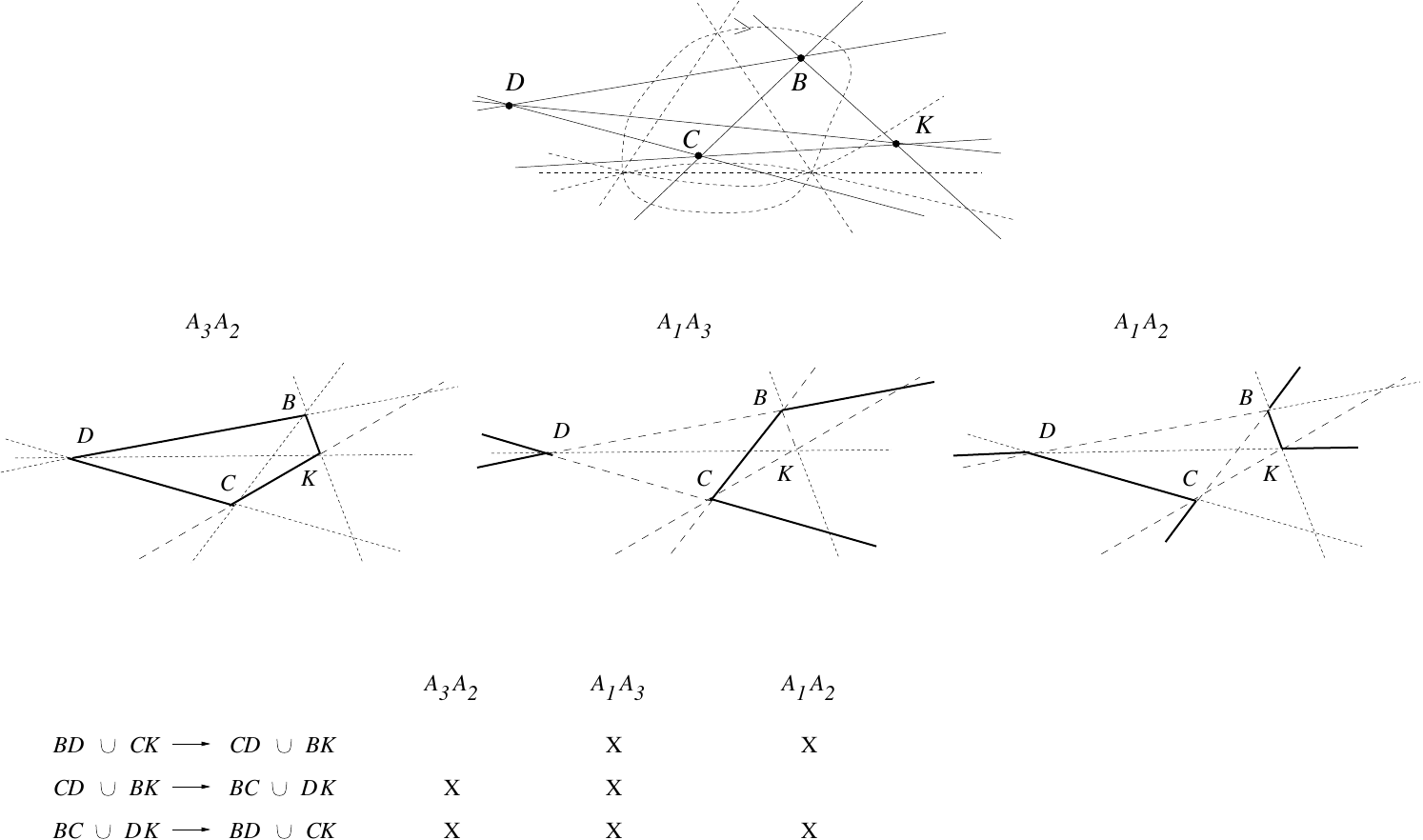}
\caption{\label{dbmkc} $K$ interior to $cr(O_1)$ in $cr(Q_3)$ swept out between $B$ and $C$}
\end{figure}

Let us consider now the ovals lying in the new zone $F_4$.
By Lemma~7, an oval $K$ in this zone may be: exterior or interior to $cr(O_3)$ in $cr(Q_2)$, interior to $cr(O_3)$ in $cr(Q_1)$ or interior to $cr(O_1)$ in $cr(Q_3)$
(in this latter case, $O_1$ is right).
Assume there exists $K$ interior to $cr(O_1)$ in $cr(Q_3)$, then $F_4$ contains the whole of the odd (negative) chain interior to $cr(O_1)$. Let $K$ be extremity of this chain. Consider the pencil of conics $\mathcal{F}_{DBKC}$, the portion $BD \cup CK \to CD \cup BK$ is totally real, see Figure~\ref{dbmkc}.
Let $T$ be an oval swept out between $B$ and $K$. If $T$ is quadrangular, it must be in the zone $F_5$, hence it is on a conic $DTBCK$. But this conic is in the totally real portion, contradiction.
If $T$ is in $cr(T_0)$, $T$ lies on a conic $DTBCK$ or $DCBTK$, in both cases, the conic is maximal with respect to all three base lines and cuts $\mathcal{O}$ at four points, contradiction. Hence, $T$ is in $cr(T_3)$. All of the ovals swept out between $B$ and $K$ are in $cr(T_3)$. The orientations of $B$ and $K$ coincide, hence there should be an odd number of ovals between them. This is a contradiction, as $B$ is by hypothesis extremity of an odd chain. $\Box$

Let us recap the new conditions on our set of ovals.
If $O_1$ is right or separating, the ovals $D, A, E, B, F, C$ ($D$, $A$, $E$, $B$, $F$, $P$, $Q$, $R$) are as in Lemma~9. If $O_1$ is left, $D$ is the image of the negative extreme oval in the chain of $O_1$, and $D, A, E, B, F, C$ satisfy the conditions of Lemma~3. In all cases, the zone $F_4$ contains no exterior ovals.


\begin{lemma}
Let $C_9$ have the third complex type. Then, $P = C$.
\end{lemma}

{\em Proof:\/} 
Recall that the third complex type may realize the three cases shown in Figure~\ref{abpqr}. One has $P = C$ in the first two cases, and $R = C$ in the third case.
Assume that $R = C$, $F$ may be chosen in such a way that it is swept out after $Q$ by $\mathcal{F}_D$ (see proof of Lemma~3). The pencil of conics $\mathcal{F}_{FBCD}$ is divided in three portions by the double lines, in two of these portions, the conics have all $36$ intersection points with the main part of $C_{18}$ and the four base ovals, these portions are totally real, see Figure~\ref{pendfcb}. The remaining ovals are swept out by the third portion $FD \cup BC \to FC \cup BD$. 
As $Q$ is swept out between $B$ and $F$ by $\mathcal{F}_D$, it is on a conic $DFCBQ$.  
The arc $CB$ cuts $\mathcal{O}$ twice, each of the arcs $QD$ and $DF$ cuts $\mathcal{O}$ once, hence $BQ$ must be entirely inside of $\mathcal{O}$. The line $A_1A_2$ separates $B$ from $Q$ in $\mathcal{O}$, as moreover $B, Q$ are not in $cr(O_3)$, the conic $DFCBQ$ is maximal with respect to $A_1A_2$, contradiction.
$\Box$.

\begin{lemma}
$D$ is negative.
\end{lemma}

{\em Proof:\/}
For $O_1$ left, this is already known. Let $O_1$ be right or separating. By Lemma~8, $F (= F_5)$ is positive exterior in $cr(Q_2)$, it contributes $+1$ to $\lambda_2$. Let $K$ be a quadrangular oval belonging to an odd chain, swept out between $B$ and $C$ by $\mathcal{F}_D$. By Lemma~9, $K$ is in the zone $F_5$.   
The pencil of conics $\mathcal{F}_{FBCD}$ has only one non totally real portion, see Figure~\ref{pendfcb}. Thus, $K$ is on a conic $DKFCB$ or $DFKCB$. 
Let us choose $F$ to be the last quadrangular oval belonging to an odd chain that is met by $\mathcal{F}_D$ before $C$, the conic is $DKFCB$. Let us call $\mathcal{Z}_F$ the triangular zone, bounded by the lines $DF$, $FC$, $BD$ containing $K$, this zone lies in $cr(Q_2) \cup cr(Q_1)$ (see positions of the base lines given in Figure~\ref{pendfcb}). The zone $\mathcal{Z}_F$ contains either only exterior ovals in $cr(Q_2)$, or both exterior ovals in $cr(Q_2)$ and the even chain interior to $cr(O_3)$. In the latter case, the chain interior to $cr(O_3)$ lies either in $cr(Q_1)$ or in $cr(Q_2)$.
The pencil $\mathcal{F}_D: B \to F$ sweeps out the ovals of $\mathcal{Z}_F$ plus (only if $O_1$ is left) an even chain interior to $cr(O_1)$ (formed by the images of the non-extreme ovals interior to $O_1$), 
Fiedler's theorem with $\mathcal{F}_D: B \to F$ implies thus that the ovals in $\mathcal{Z}_F \cap cr(Q_2)$ ($F$ included, $D$ not) contribute $+1$ to $\lambda_2$.
The contribution of the ovals in $\mathcal{Z}_F$ ($D$ now included) to $\lambda_2$ is $0$ or $+2$. 
Consider now the maximal pencil $\mathcal{F}_B: D \to F$ sweeping out $\mathcal{Z}_F$. Any oval $K$ met by this pencil is on a conic $DKFCB$, so $K$ is in $\mathcal{Z}_F$. The ovals $D$ and $F$ are connected via this pencil by an even number of ovals, their orientations are opposite. $\Box$

\begin{lemma}
There exists an oval $J$ exterior positive in $cr(Q_2)$, lying on conics $DFCJA$, $DFCJB$, $ABDCJ$.
\end{lemma}

\begin{figure}[htbp]
\centering
\includegraphics{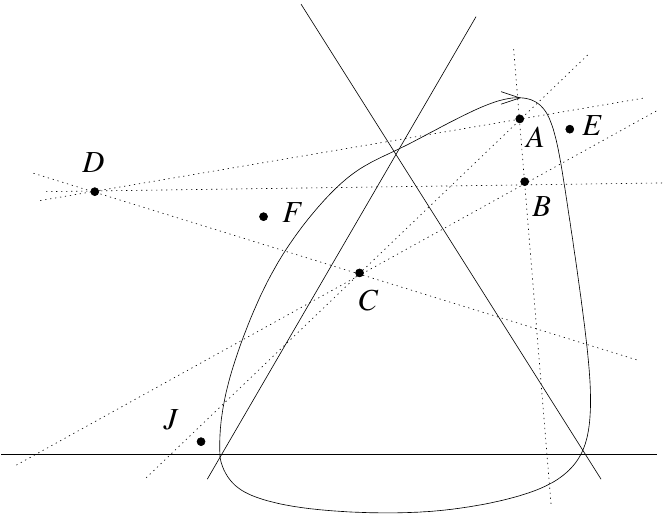}
\caption{\label{abdcj} Position of $J$}
\end{figure}

{\em Proof:\/}
The ovals of $cr(Q_2)$ in $\mathcal{Z}_F$ (including now $D$ and $F$) bring a contribution $0$ to $\lambda_2$. As $\lambda_2 = 1$ (first complex type) or $2$ (third complex type), there must be some supplementary positive oval $J$ extremity of an odd chain in $cr(Q_2)$, that is not swept out by the pencil $\mathcal{F}_D: B \to C$, the conic through $D, B, C, F, J$ is thus $DFCJB$. For $O_1$ left, the ovals interior to $cr(O_1)$ other than $D$ form an even chain, so we may assume that $J$ is not interior to $cr(O_1)$.
The conic through $D, A, C, F, B$ is $DFCAB$, see Figure~\ref{acebf} with $F = F_5$. In Figure~\ref{pendfcb} one may replace $B$ by $A$. The conic through $D, A, C, F, J$ is  $DFCJA$ or $DFCAJ$. In the latter case, $J$ is swept out between $A$ and $B$ by $\mathcal{F}_D$, hence $J$ must verify the conditions required for $E$, see Lemma~5. As $J$ is in $cr(Q_2)$, it should be in $cr(O_1)$, contradiction. So $J$ is on a conic $DFCJA$. 
Using the conics $DFCJA$, $DFCJB$, $DFCAB$, we get the orderings 
$C: D, J, A, B$ and $D: C, J, A, B$ for the pencils of lines based at $C$ and $D$ respectively, hence the conic through $A, B, C, D, J$ is either $ABDCJ$ or $ABCDJ$. If the conic is $ABCDJ$, its arc $DJA$ cuts $A_1A_2$ and $A_1A_3$ (see pencil $\mathcal{F}_{ABCD}$ in Figure~\ref{daebc}), as
$D, J$ are in $cr(Q_2)$ and $A$ is in $cr(Q_3)$, the arc $JA$ cuts $A_1A_2$ and $A_1A_3$. The arc $JA$ of $DFCJA$ is homotopic with fixed extremities to the arc $JA$ of $ABCDJ$. The conic $DFCJA$ is maximal with respect to all three base lines
and cuts $\mathcal{O}$ four times, contradiction. Thus, the conic through $A, B, C, D, J$ is $ABDCJ$, see Figure~\ref{abdcj}. $\Box$

\begin{lemma}
Let $G$ be an exterior oval in $cr(Q_1)$. Then, $G$ lies on a conic $DGAEB$.
\end{lemma} 

\begin{figure}[htbp]
\centering
\includegraphics{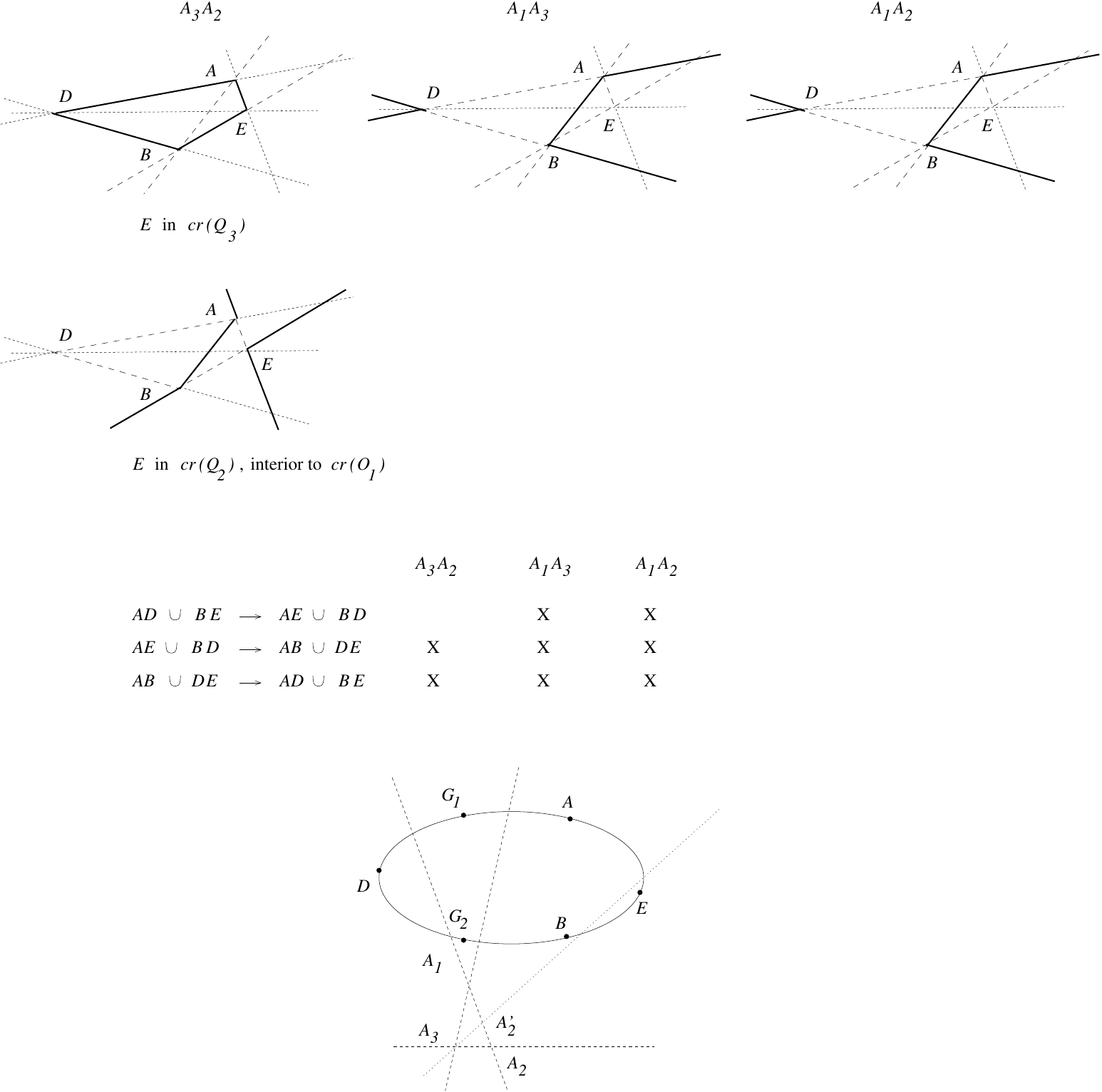}
\caption{\label{dageb} Pencil of conics $\mathcal{F}_{ADBE}$, position of $G$}
\end{figure}

{\em Proof:\/} Consider the pencil of conics $\mathcal{F}_{DAEB}$, all of the conics in it cut $\mathcal{O}$ four times, and two of the three portions determined by the double lines are totally real, see Figure~\ref{dageb}.
In this figure, we have indicated the positions of the base lines, distinguishing for $A_3A_2$ the two cases $E$ in $cr(Q_3)$ and $E$ in $cr(Q_2)$, interior to $cr(O_1)$.
Any remaining oval $H$ will be swept out by the third portion $AD \cup BE \to AE \cup BD$. 
Assume that the conic through $H$ is $DAHEB$ or $DAEHB$, it cannot be maximal with respect to $A_2A_3$ hence $H$ is in $cr(Q_3)$, or in $cr(Q_2)$ and interior to $cr(O_1)$. Let $G$ be an exterior oval in $cr(Q_1)$ and consider the conic of the pencil  $\mathcal{F}_{DAEB}$ passing through $G$, it is $DGAEB$ or $DAEBG$.
In the lower part of the Figure, we have drawn the conic and indicated on it the two admissible positions $G_1$ (conic $DG_1AEB$) and $G_2$ (conic $DAEBG_2$)  of $G$. We have also indicated two positions of the second base point: $A_2$ corresponding to the case $E$ in $cr(Q_3)$ and $A'_2$ corresponding to the case $E$ in $cr(Q_2)$, interior to $cr(O_1)$. The conic through $A, B, D, E, G$ is either $DAEBG$ or $DGAEB$.
The conic through $C, D, F, B, G$ is in the portion $FD \cup CD \to FD \cup BC$ of the pencil $\mathcal{F}_{DCFB}$, it is $DFCGB$ (otherwise, $G$ would be swept out by $\mathcal{F}_D: B \to C$, but the exterior quadrangular ovals met by this pencil
lie all in $cr(Q_2)$, see proof of Lemma~11) 
Similarly, we obtain the conic $DFCGA$. If $G$ lies on a conic $DAEBG$, then the pencil of lines based at $D$ sweeps out successively: $A, E, B, G, A_1, A_2, A_3$ ($cr(T_0)$, hence $C$ is swept out in the portion $A_1 \to A_2 \to A_3$). We have thus the ordering $D: A, E, B, \{ G, F \}, C$, contradiction with the conic $DFCGB$. So, $G$ lies on a conic $DGAEB$. $\Box$

\begin{figure}[htbp]
\centering
\includegraphics{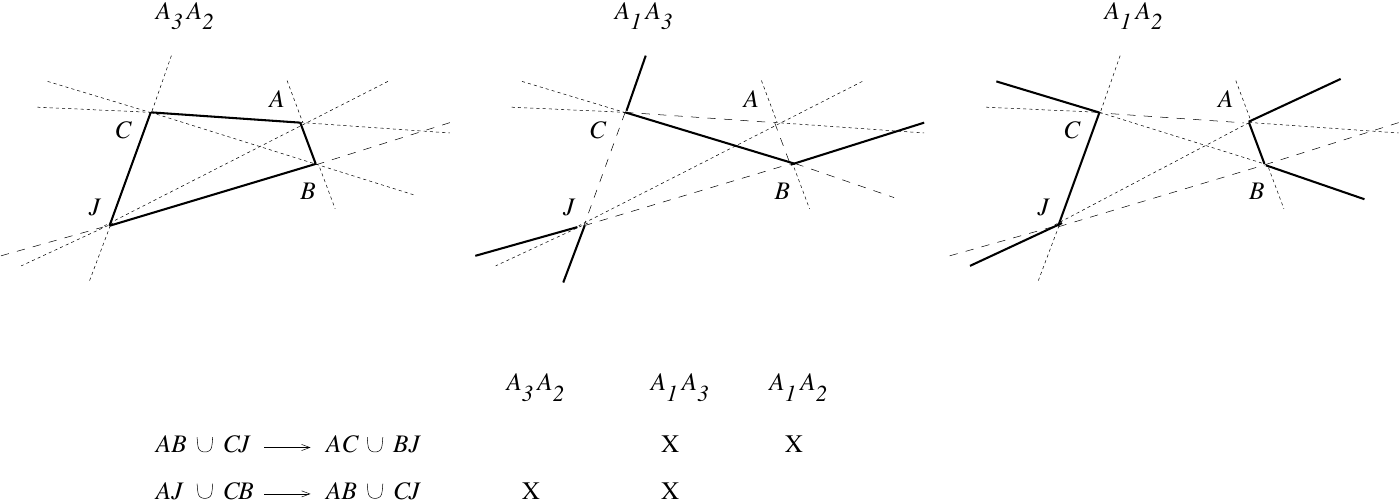}
\caption{\label{penabcj} Pencil of conics $\mathcal{F}_{ABCJ}$}
\end{figure}

We will now find a contradiction using the configuration of eight ovals
$D, A, E, B, F, C, J, G$.
With help of the conics $DGAEB$, $DFCJA$, $DFCAB$, $DAEBC$ and $ABDCJ$, we find the ordering for the pencil of lines based at $A: D, F, C, J, B, E, G$. The conic through $A, G, C, J, B$ can be $ABJCG$, $AGBJC$, $ABGCJ$ or $AJBGC$. With help of $DGAEB$, $DFCJB$, $DFCAB$, $DCJAB$, $DAEBC$, we find the ordering for the pencil based at $B$: $D$, $F$, $C$, $J$, $E$, $A$, $G$. 
This rules out $AGBJC$ and $AJBGC$. The conic is $ABGCJ$ or $ABJCG$. 
Consider the pencil $\mathcal{F}_{ABJC}$, see Figure~\ref{penabcj} where we considered only the two relevant portions.
If the conic is $ABJCG$, the arc $CGA$ doesn't cut $A_1A_3$, and cuts $A_1A_2$ (once), as moreover, $G$ is exterior in $cr(Q_1)$, this arc must be maximal with respect to $A_2A_3$, contradiction. If the conic is $ABGCJ$, the arc $BGC$ cuts both $A_1A_3$ and $A_2A_3$ (once each), and must be maximal with respect to $A_1A_2$, contradiction. 


The first and third complex type are not realizable. 

\subsection{The second complex type}
\begin{figure}[htbp]
\centering
\includegraphics{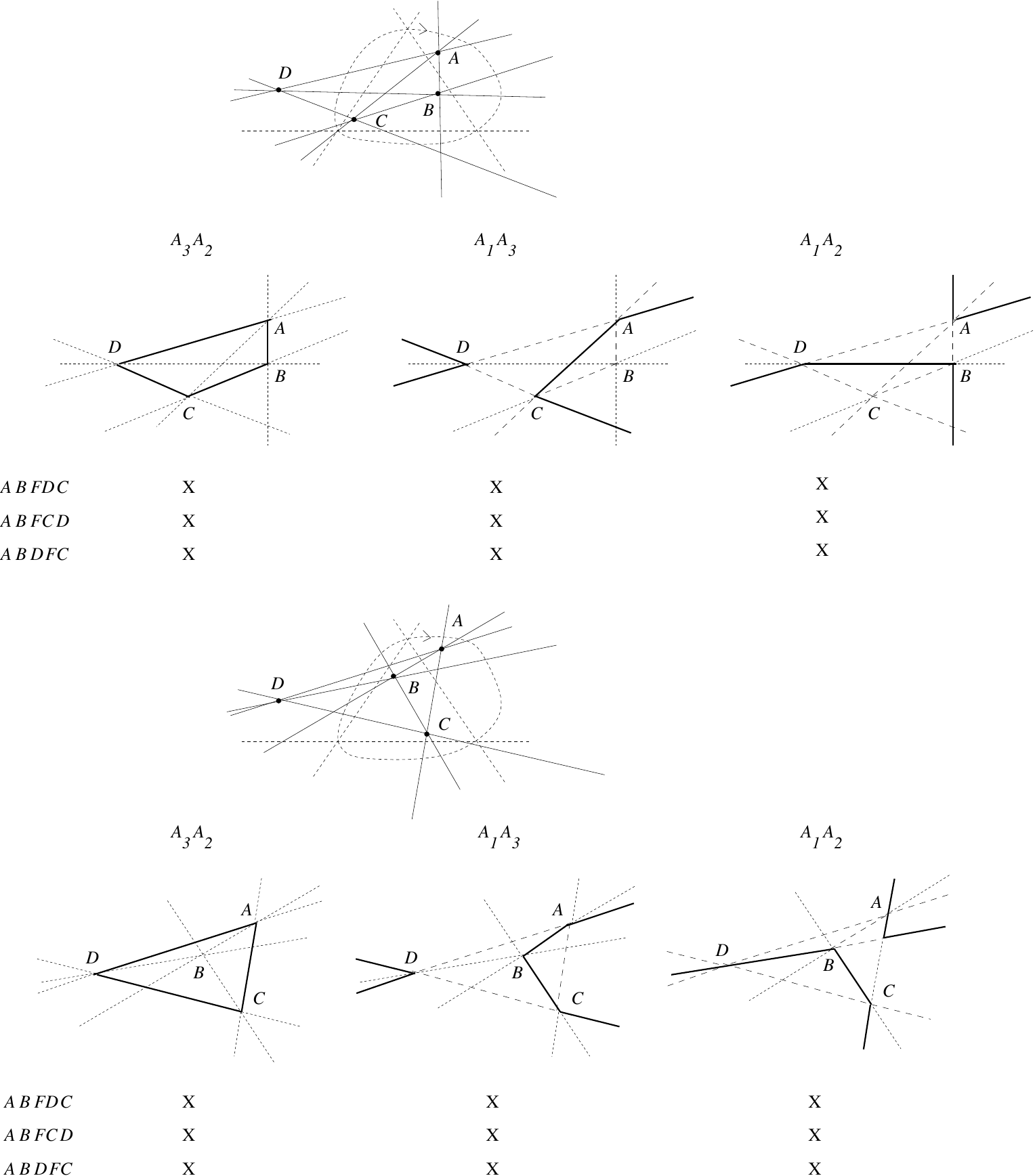}
\caption{\label{ef1} The second complex type}
\end{figure}

Let now $C_9$ realize the second complex type. Let $[DC]$ be the segment of the line $DC$ that doesn't cut $A_1A_3$, and $[DC]'$ be the other segment.
We distinguish two cases, depending on whether the intersection $AB \cap DC$ lies on $[DC]$ or on $[DC]'$, see Figure~\ref{ef1}. 
The oval $F$ is swept out between $B$ and $C$ by $\mathcal{F}_D$, the sector containing $F$ is divided in four triangles by the lines $AB$, $AC$, $BC$, one of them is entirely interior to $\mathcal{O}$ so $F$ is in one of the other three.
In both cases, the conic through $A, B, C, D, F$ is $ABFDC$, $ABFCD$ or $ABDFC$. 
Consider the pencil of conics $\mathcal{F}_{ABCD}$.
In both cases, the portion $AB \cup CD \to BD \cup AC$ is totally real, hence the conic must be $ABFCD$, it is in the portion $AB \cup CD \to AD \cup BC$.
The conics of this portion are maximal with respect to: $A_1A_3$ and $A_1A_2$ in the first case, $A_1A_2$, $A_3A_2$ in the second case.
In the first case, the triangle containing $F$ is divided in two pieces by $A_2A_3$, one in $cr(T_0)$, the other in $cr(T_1) \cup cr(Q_1)$. As $F$ is quadrangular, $F$ is in $cr(Q_1)$, the arc $BFC$ of $ABFCD$ is maximal with respect to $A_2A_3$, and $ABFCD$ cuts $\mathcal{O}$ four times, contradiction. In the second case, the triangle containing $F$ is divided in two pieces by $A_1A_3$, one in $cr(T_0)$, the other in $cr(T_2) \cup cr(Q_2)$. As $F$ is quadrangular, $F$ is in $cr(Q_2)$, the arc $BFC$ of $ABFCD$ is maximal with respect to $A_1A_3$, and $ABFCD$ cuts $\mathcal{O}$ four times, contradiction. The second complex type is not realizable. This finishes the proof of Theorem~1. $\Box$

\section{$M$-curves even, even, odd with jump}
\subsection{Cremona transformation again}
\begin{figure}[htbp]
\centering
\includegraphics{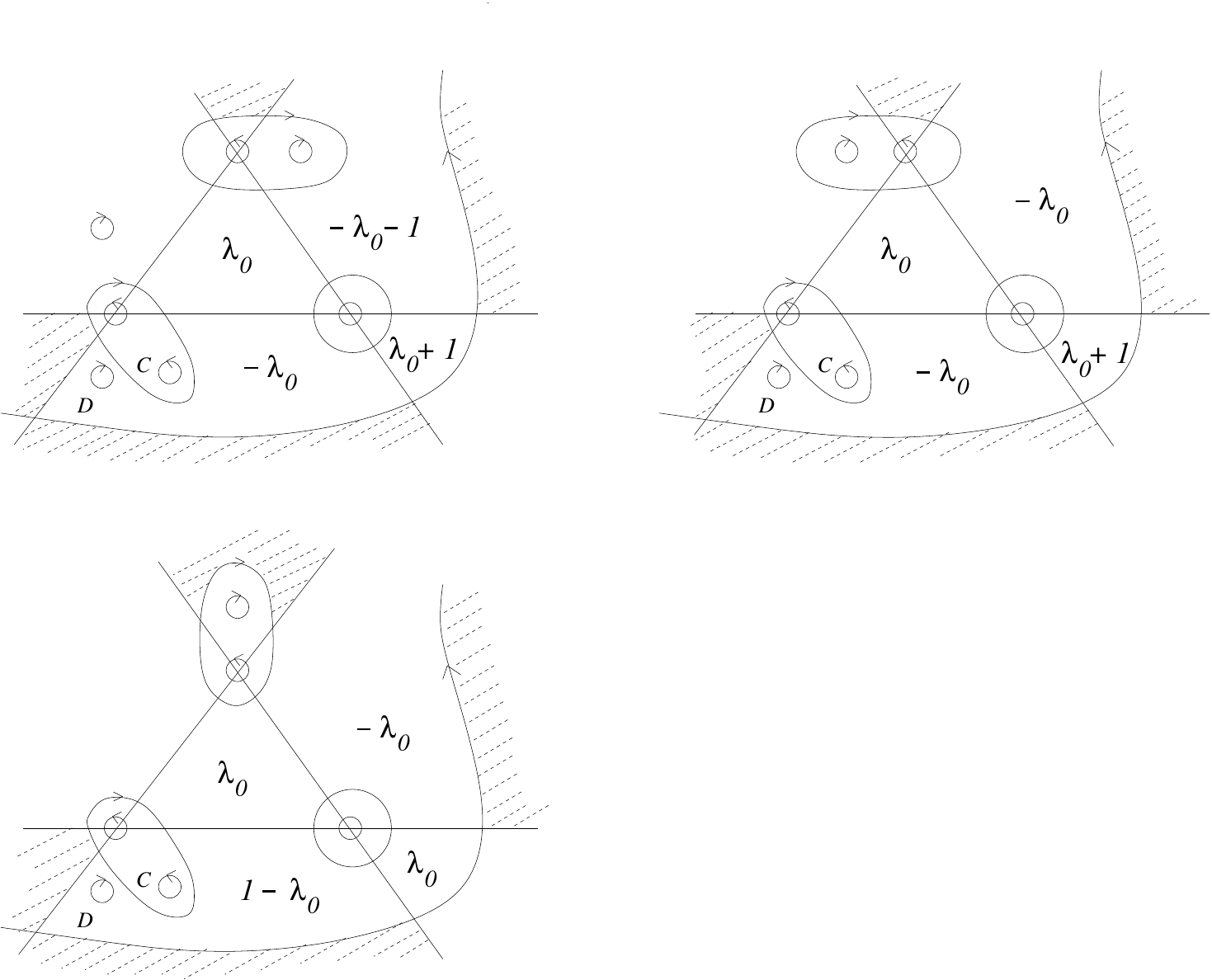}
\caption{\label{xing02} Crossing jump, no exterior ovals in $T_1$ ($O_1$ right, left, up)}
\end{figure}

\begin{figure}[htbp]
\centering
\includegraphics{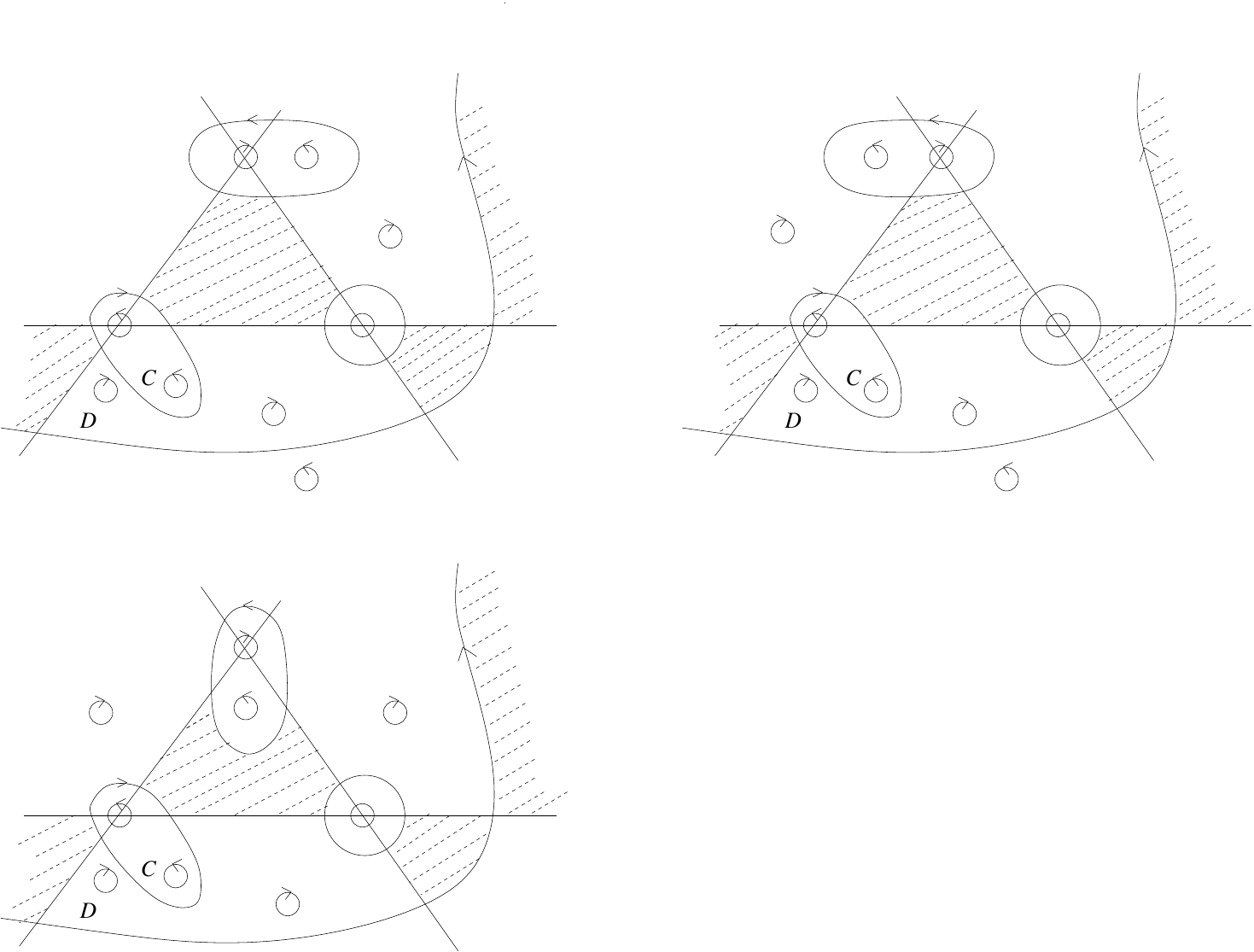}
\caption{\label{crosses} Crossing jump, no exterior ovals in $T_0 \cup T_2$ ($O_1$ right, left, up)}
\end{figure}

\begin{figure}[htbp]
\centering
\includegraphics{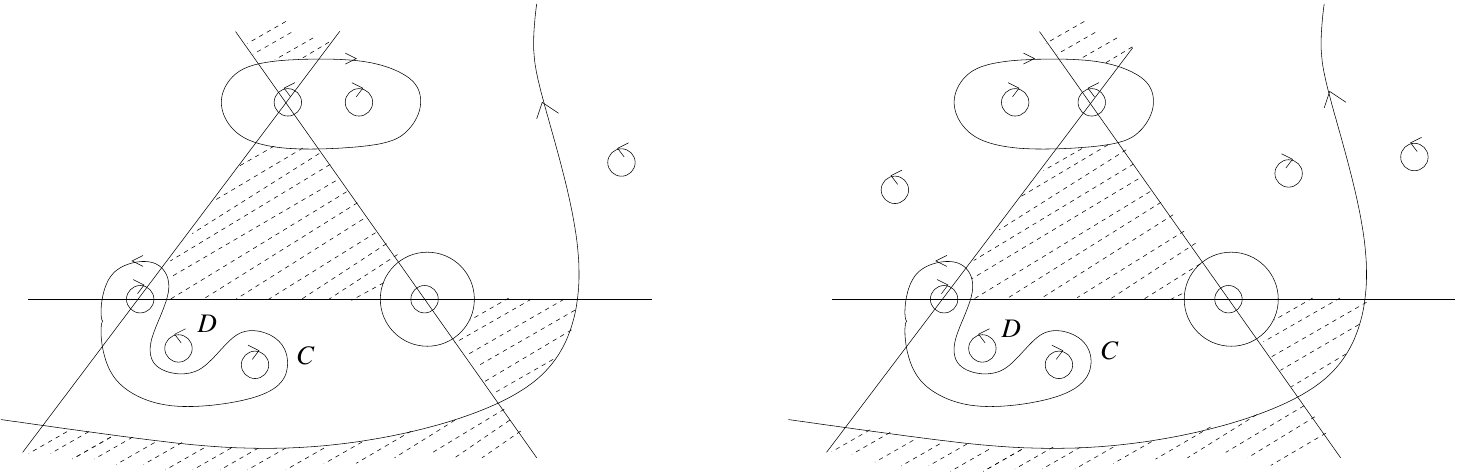}
\caption{\label{crossnot} Non-crossing jump ($O_1$ right, left)}
\end{figure}

Assume there exists an $M$-curve $C_9$ even, even, odd with jump, 
it realizes one of the complex types listed in Table~\ref{orev10}.
Recall that if $O_3$ is positive (negative), $O_3$ is crossing (non-crossing), see section 1.2. Let $\epsilon_3 = +1$ ($-1$) if $O_3$ is positive (negative). The six principal ovals, the base lines, and $\mathcal{J}$ divide the plane in zones. The contributions of the non-principal ovals to $\Lambda_+ - \Lambda_-$ in each zone are indicated in Figures~\ref{xing02}-\ref{crossnot}. 

\begin{figure}[htbp]
\centering
\includegraphics{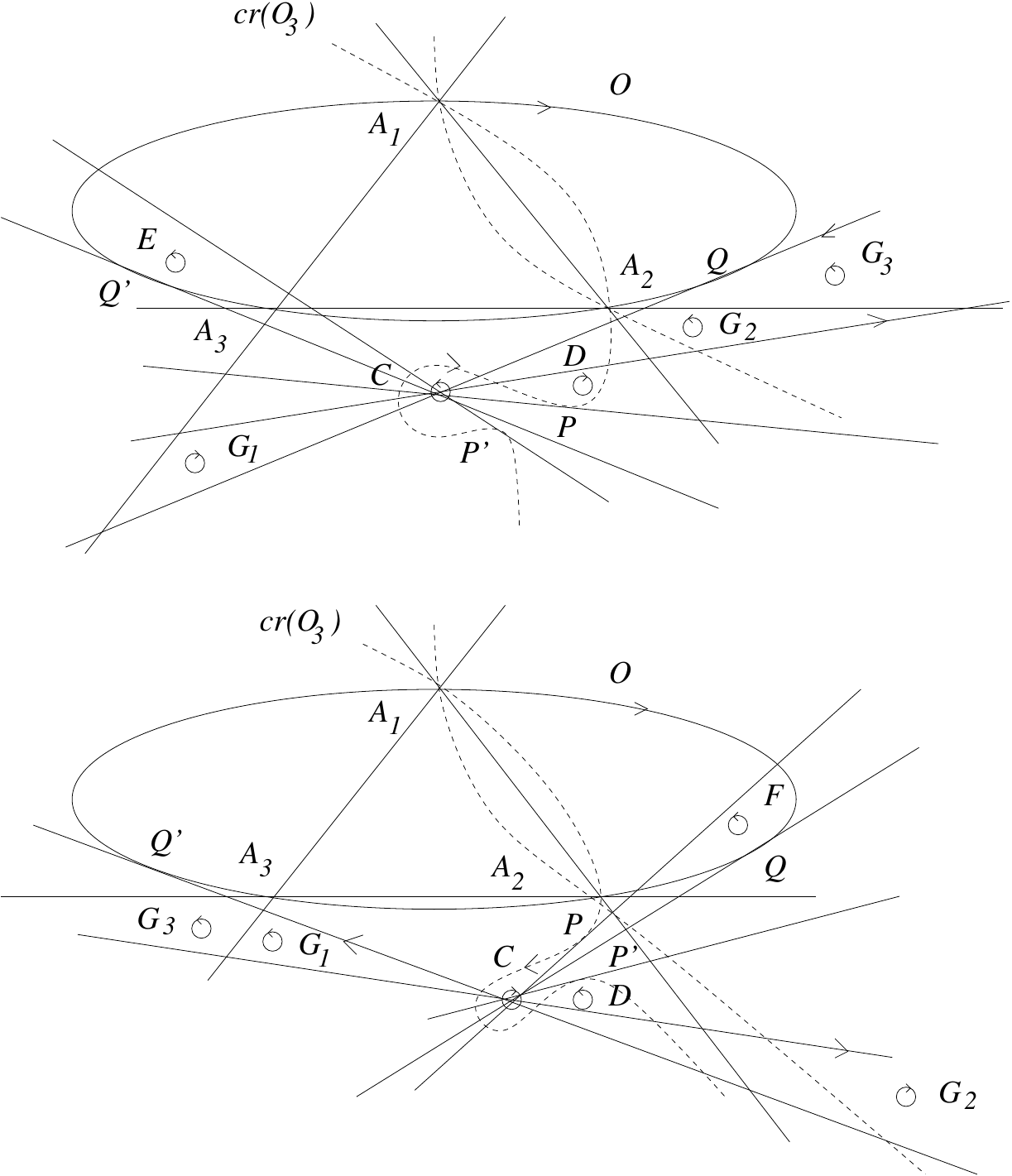}
\caption{\label{lambdaquad} Curve $C_{18}$ and pencil $\mathcal{F}_C$, $O_3$ crossing and $O_3$ non-crossing}
\end{figure}

Let us perform again the Cremona transformation $cr$ based at $A_1, A_2, A_3$, $C_9$ is mapped onto a curve $C_{18}$ with three $9$-fold singularities. 
We use the same conventions as in the case without jump, see section 2.1.
Let us choose two ovals $C$ and $D$ of $C_9$ such that $(A_3, D, C)$ are extremities of the three successive jump chains. 
In Figure~\ref{lambdaquad}, we have represented the image $\mathcal{O}$ of $\mathcal{J}$, plus the images of $O_3$, $C$ and $D$, the upper part corresponds to $O_3$ crossing, the lower part to $O_3$ non-crossing.
The pencil of conics $\mathcal{F}_{A_1A_2A_3C}$ sweeping out $C_9$ is mapped by $cr$ onto a pencil of lines $\mathcal{F}_C$. Let $P, P'$ be the tangency points of $\mathcal{F}_C$ with $cr(O_3)$ and let $Q, Q'$ be the points of tangency of $\mathcal{F}_C$ with $\mathcal{O}$. 
By Lemma~17 from \cite{fi}, the ovals in $\mathcal{O}$ are all met consecutively by $\mathcal{F}_C$, they form a Fiedler chain: {\em triangular ovals\/} $\to P'$ ($O_3$ crossing) or $P \to$ {\em triangular ovals\/} ($O_3$ non-crossing).   
For $O_3$ crossing, one has $\Lambda = \lambda_0 -\lambda_4 -\lambda_5 = -1$. Denote by $E$ the first oval of the triangular Fiedler chain, in the Figure, $E$ has been placed arbitrarily in $cr(T_2)$, it could be as well in $cr(T_1)$ or $cr(T_0)$.
For $O_3$ non-crossing, one has $\Lambda = -\lambda_6 = -1$. Denote by $F$ the last oval of the triangular chain.  
Applying Fiedler's theorem with the pencil $\mathcal{F}_C$, we get the

\begin{lemma}
Let $\mu = \lambda_1 + \lambda_2 - \lambda_3$, one has $\mu = \epsilon_3$.
For $O_3$ crossing, the quadrangular ovals are swept out by $\mathcal{F}_C$ beween $Q$ and $E$. Between $Q$ and $D$, there is one single Fiedler chain, whose contribution to $\mu$ is $1$, between $D$ and $E$, there are two Fiedler chains, starting at $\{ D, P \}$ and ending at $\{ Q', E\}$, they contribute $0$ to $\mu$. 
For $O_3$ non-crossing, the quadrangular ovals are swept out by $\mathcal{F}_C$ between $F$ and $Q'$. Between $F$ and $D$, there are two Fiedler chains starting at $\{ F, Q \}$ and ending at $\{ P', D \}$, they contribute $0$ to $\mu$, between $D$ and $Q'$ there is one single Fiedler chain whose contribution to $\mu$ is $-1$.
\end{lemma}

(In all these Fiedler chains, the contributions of $D$, $E$ and $F$ are not included by convention. As $D$ and $C$ have opposite orientations, they contribute together $0$ to $\mu$.) 
The identity $\mu = \epsilon_3$ may actually be obtained directly from Table~\ref{orev10}.
For $O_3$ crossing, there exists an oval $G$, swept out between $Q$ and $D$, that is positive in $cr(Q_1) \cup cr(Q_2)$ or negative in $cr(Q_3)$, we denote by $G_1, G_2, G_3$ the three admissible choices of $G$.
For $O_3$ non-crossing, there exists an oval $G$ swept out between $D$ and $Q'$, that is negative in $cr(Q_1) \cup cr(Q_2)$ or positive in $cr(Q_3)$, we denote by $G_1, G_2, G_3$ the three admissible choices of $G$. See Figure~\ref{lambdaquad}.

\subsection{Proof of Theorem~2}
The piece of pencil $\mathcal{F}_C$ sweeping out the quadrangular ovals is divided in two portions by the line $CD$, the contribution of each portion to $\mu= \lambda_1 + \lambda_2 - \lambda_3$ is given by Lemma~14. In the next Lemmas, we will find pencils of conics based at $C$, $D$ and two other ovals, having only one non-totally real portion. This allows to  study the distribution of the quadrangular ovals between the two parts of $\mathcal{F}_C$ and sometimes to find a contradiction.
To this purpose, we may ignore the even chains and treat each odd chain as if it consisted of a single oval. So we may assume without loss of generality that $O_2$ contains only the one oval $A_2$, and that each of the three jump chains consists of a single oval $(A_3, D, C)$.  


\begin{lemma}
If $C_9$ has complex type $(+, n), (\pm, \mp), (-, +, +)$, then $O_1$ is right. 
\end{lemma}

\begin{figure}[htbp]
\centering
\includegraphics{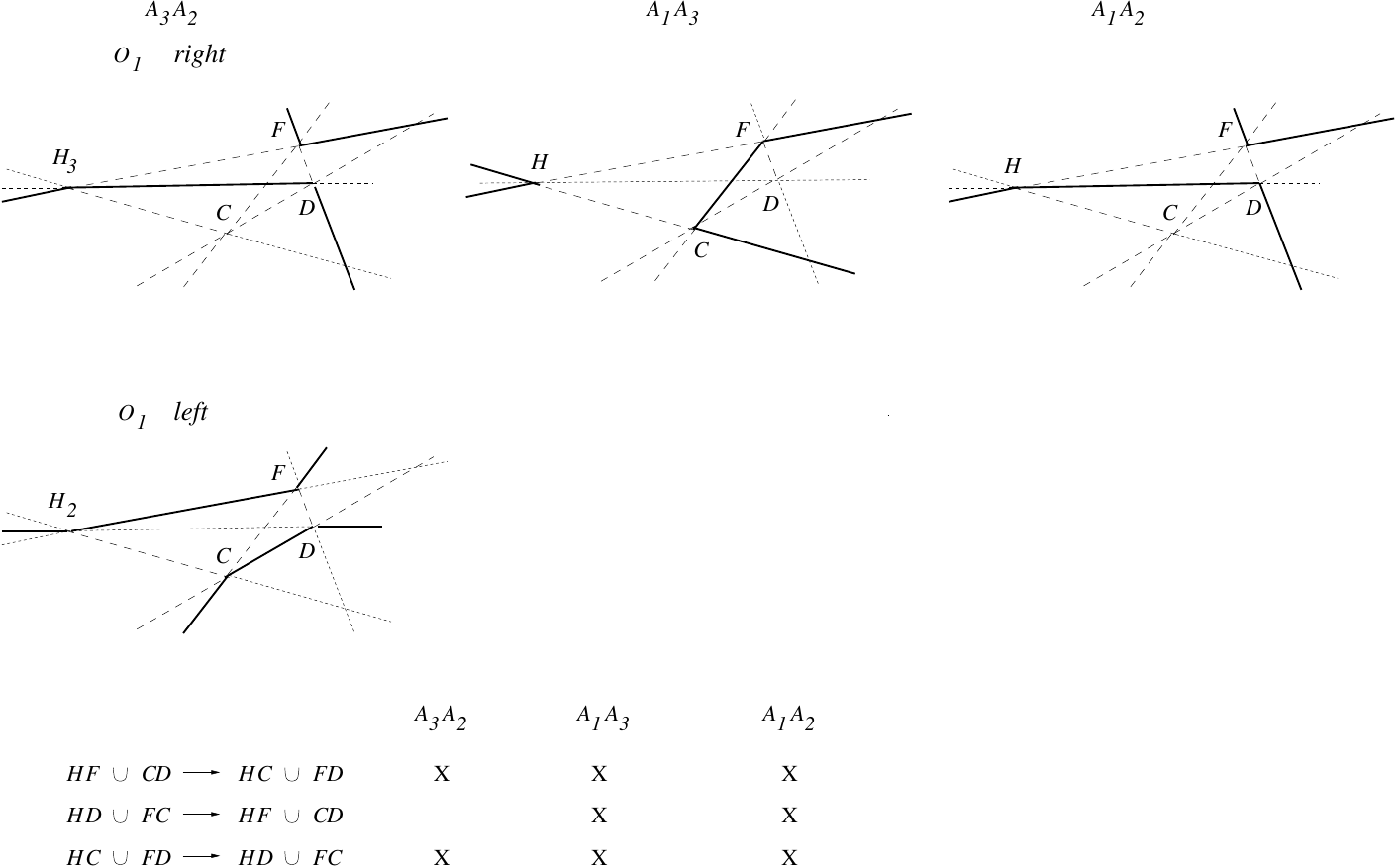}
\caption{\label{hfcd} Non-crossing jump, pencil of conics $\mathcal{F}_{CDFH}$}
\end{figure}

\begin{figure}[htbp]
\centering
\includegraphics{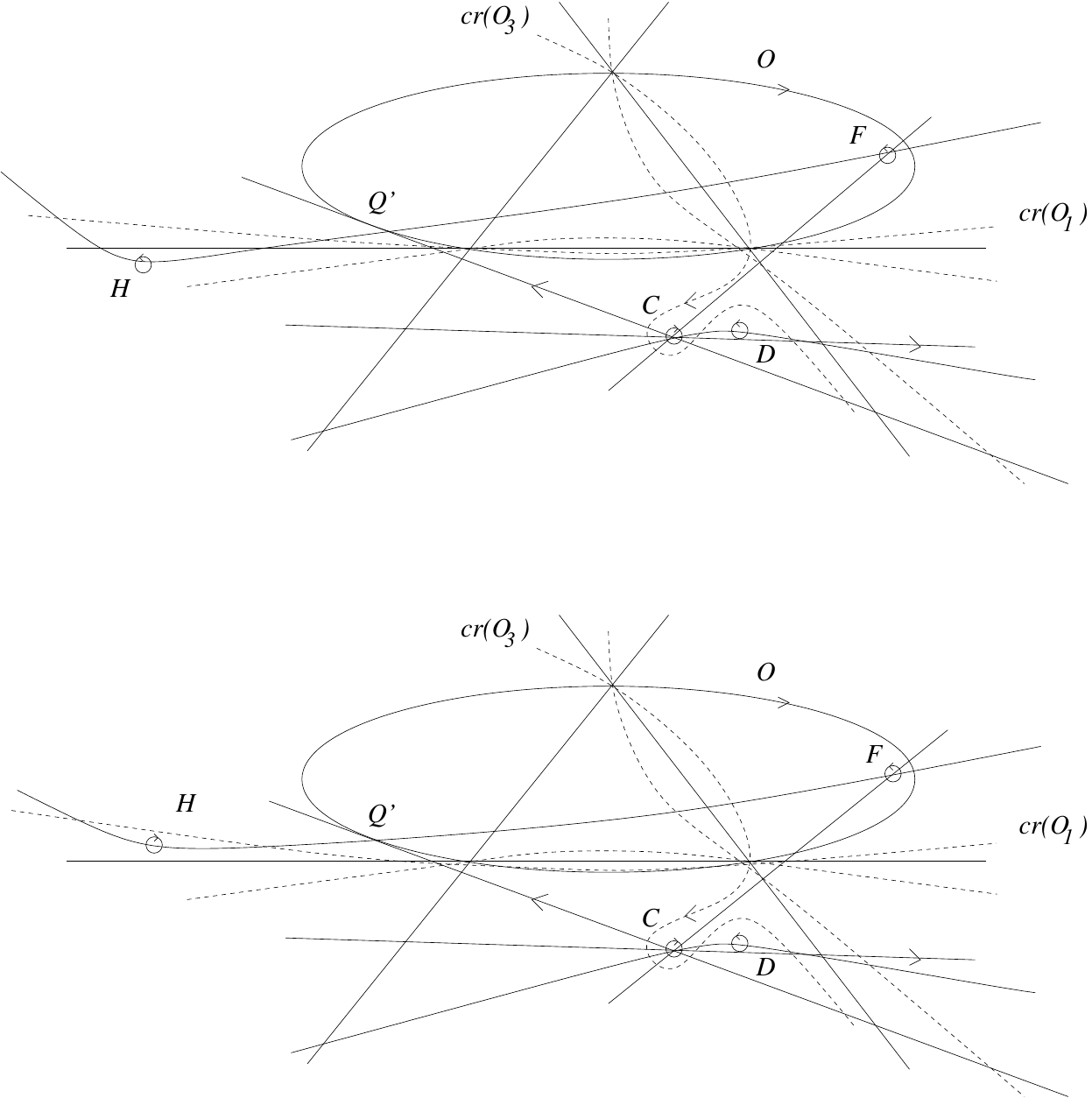}
\caption{\label{noncross} Non-crossing jump, conic through $C, D, F, H$ and a supplementary oval $I$ (non-represented), $O_1$ right and $O_1$ left}
\end{figure}

{\em Proof:\/}
Let $C_9$ have this complex type, it has non-crossing jump, see Figure~\ref{crossnot}. Let $H$ be the positive extremal oval in the chain of $O_1$. If $O_1$ is right, $H = H_3$ is in $Q_3$, if $O_1$ is left, $H = H_2$ is in $Q_2$.
After $cr$, consider the pencil of conics $\mathcal{F}_{CDFH}$. Two portions are maximal with respect to the three base lines, see Figures~\ref{hfcd}.
All conics of the pencil have supplementarily two intersection points with $cr(O_3)$ and two intersection points with $\mathcal{O}$. 
So $\mathcal{F}_{CDFH}$ has only one non-totally real portion, $CD \cup FH \to CH \cup DF$, that sweeps out any empty non-principal oval $I$ of $C_9$ other than $C, D, F, H$, see Figures~\ref{hfcd}-\ref{noncross}. In Figure~\ref{noncross}, we have represented a conic of this portion for either case $O_1$ right and $O_1$ left. Let $I$ be a quadrangular oval. 
If $I$ is in $cr(Q_1)$ or $I$ is exterior to $cr(O_1)$ in $cr(Q_3)$, $I$ is on the arc $FC$ of the conic, hence $I$ is swept out between $F$ and $D$ by $\mathcal{F}_C$. If $I$ is in $cr(Q_2)$ or interior to $cr(O_1)$ in $cr(Q_3)$, $I$ is on the arc $DHF$ of the conic, hence $I$ is swept out between $D$ and $Q'$. For $O_1$ left, one should have $\lambda_3 - \lambda_1 = 0$ and
$\lambda_2 = -1$ (by Lemma~14). But actually $\lambda_1 = \lambda_2 = 0$ and $\lambda_3 = 1$, contradiction. Let  $\lambda^{ext}_3$ and $\lambda_3^{int}$ be the respective contributions to $\lambda_3$ of the exterior ovals and of the ovals interior to $O_1$. 
For $O_1$ right, one must have $\lambda^{ext}_3 - \lambda_1 = 0$ and $\lambda_3^{int} - \lambda_2 = 1$, there is no contradiction. $\Box$


Let now $O_3$ be crossing.
We will say that a quadrangular oval is $Q \to D$ ($D \to E$) if it is swept out between $Q$ and $D$ ($D$ and $E$) by the pencil of lines $\mathcal{F}_C$.

\begin{lemma}
The complex type
$(-, n), (\pm, \mp, n), (+, -, -)$ is not realizable, and the complex type
$(-, u), (\pm, \mp, n), (+, -, -)$ is such that $T_1$ contains no exterior ovals.
\end{lemma}

\begin{figure}[htbp]
\centering
\includegraphics{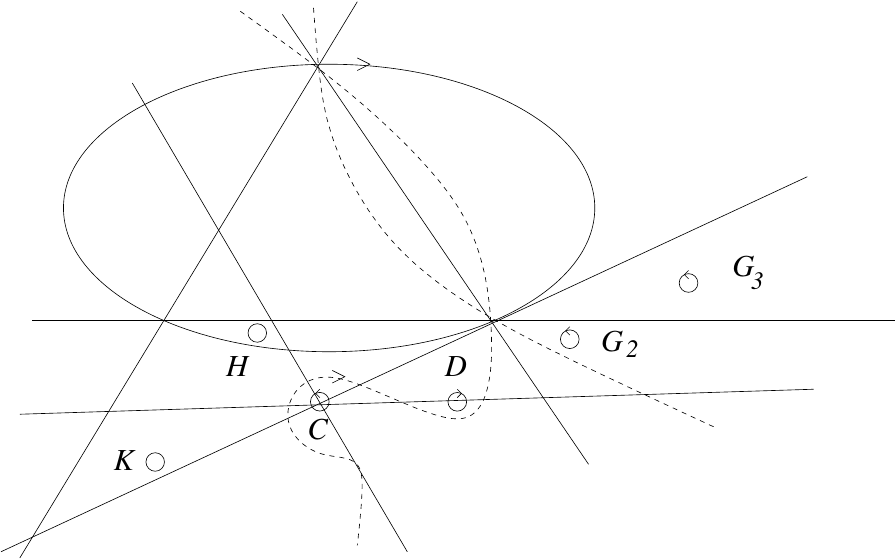}
\caption{\label{ovalinT1} Crossing jump, $O_1$ negative, oval $H$ exterior in $cr(T_1)$}
\end{figure}

\begin{figure}[htbp]
\centering
\includegraphics{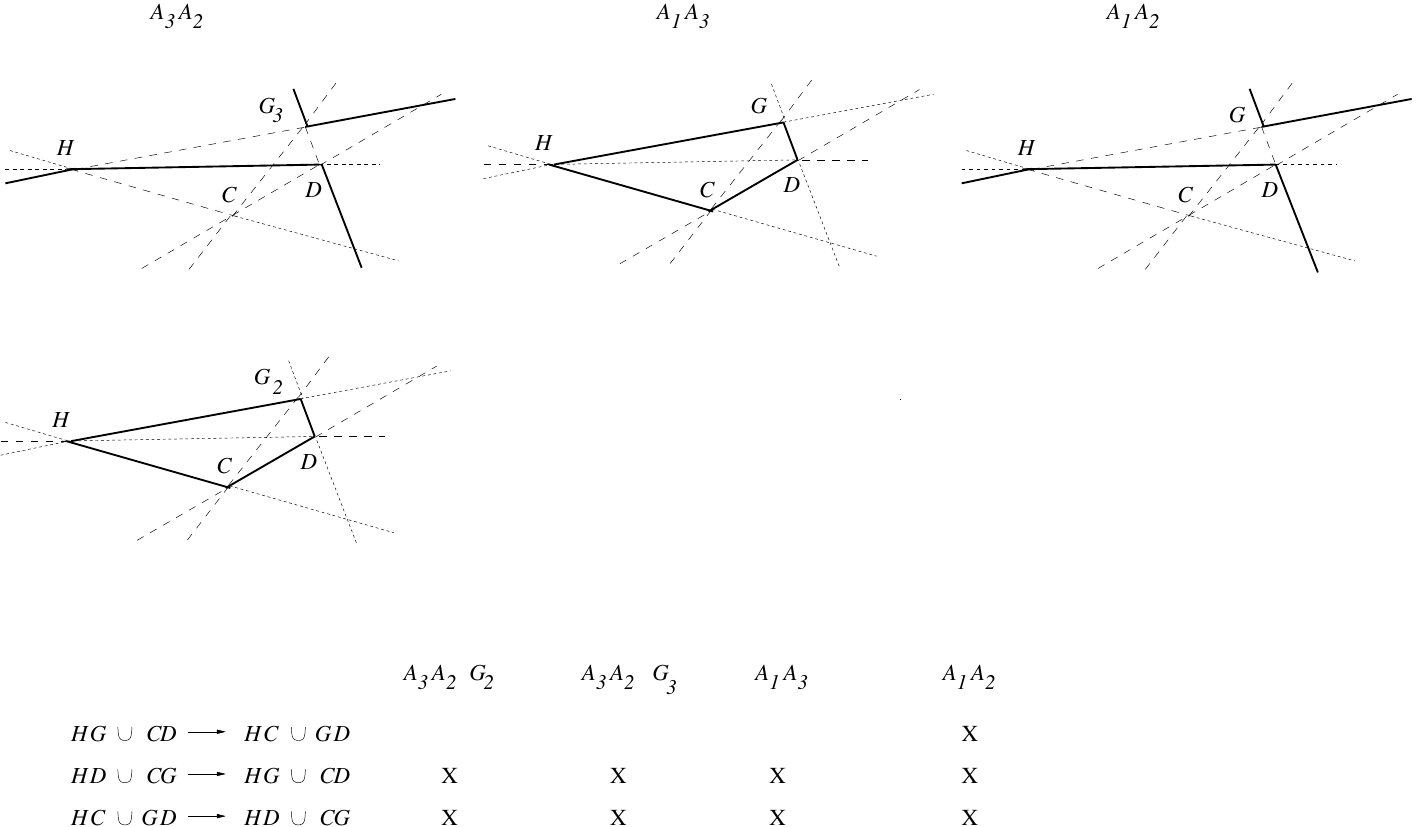}
\caption{\label{hgcd} Crossing jump, $O_1$ negative, pencil of conics $\mathcal{F}_{CDGH}$}
\end{figure}

{\em Proof:\/}
Let $C_9$ realize one of these complex types, $O_1$ may be non-separating,
right or left, see Figure~\ref{crosses}.
Assume there is an exterior oval $H$ in $cr(T_1)$ (this is true for $O_1$ non-separating). Assume there exists an (exterior) oval $K$ that is $Q \to D$ in $cr(Q_1)$, see Figure~\ref{ovalinT1}.
Consider the pencil of conics $\mathcal{F}_{CDHK}$. The three base lines are in the M\"obius band obtained cutting away the affine (in the plane of the figure) triangle $DHK$ from $\mathbb{R}P^2$. All conics of the pencil are maximal with respect to each base line. Moreover, the conics cut $\mathcal{O}$ twice, and $cr(O_3)$ at two supplementary points, hence the pencil is totally real, contradiction. The pencil of lines $\mathcal{F}_C$ meets no oval in $cr(Q_1)$ between $Q$ and $D$, hence there exists a positive oval $G = G_2$ or a negative oval $G = G_3$ that is $Q \to D$ in $cr(Q_2) \cup cr(Q_3)$, see Figure~\ref{ovalinT1}. Consider now the pencil of conics $\mathcal{F}_{CDGH}$. This pencil has only one non totally real portion:
$CD \cup GH \to CH \cup DG$, see Figure~\ref{hgcd}.
(If $G_3$ is not interior to $cr(O_1)$, this portion is maximal also with respect to the line $A_2A_3$.)
Let $I$ be a quadrangular oval. If $I$ is in $cr(Q_2) \cup cr(Q_3)$, $I$ is on the arc $HGD$ of the conic, hence $I$ is $Q \to D$.
If $I$ is in $cr(Q_1)$, $I$ is on the arc $CH$ of the conic, hence $I$ is swept out
between $D$ and $H$ by $\mathcal{F}_C$.
One must have thus $\lambda_2 - \lambda_3 = 1$ (Lemma~14). 
But for all three cases $O_1$ left, right and up, $\lambda_2- \lambda_3 = 0$, 
contradiction. $\Box$

\begin{lemma}
The complex type $(-, u), (\pm, \mp, n), (+, -, -)$ is not realizable.
\end{lemma}

\begin{figure}[htbp]
\centering
\includegraphics{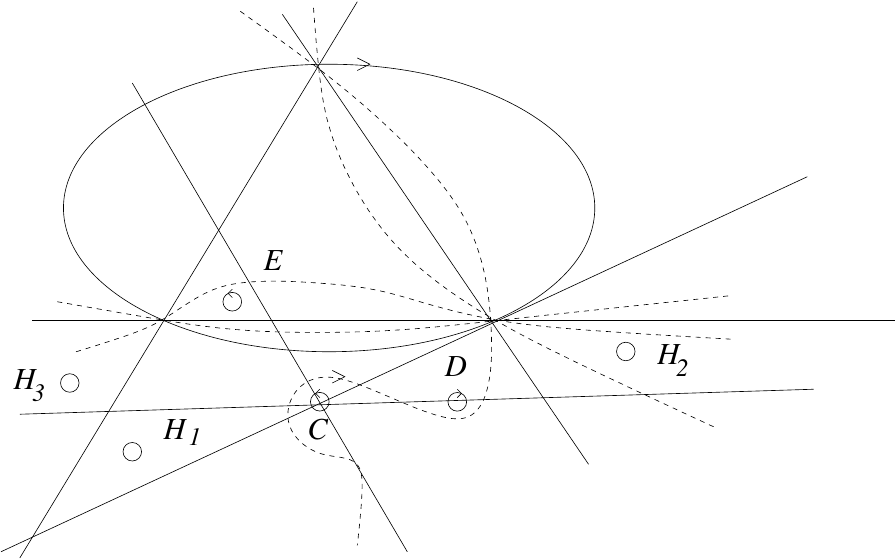}
\caption{\label{ovalinT0} Crossing jump, $O_1$ negative, up, three positions for the oval $H$}
\end{figure}

\begin{figure}[htbp]
\centering
\includegraphics{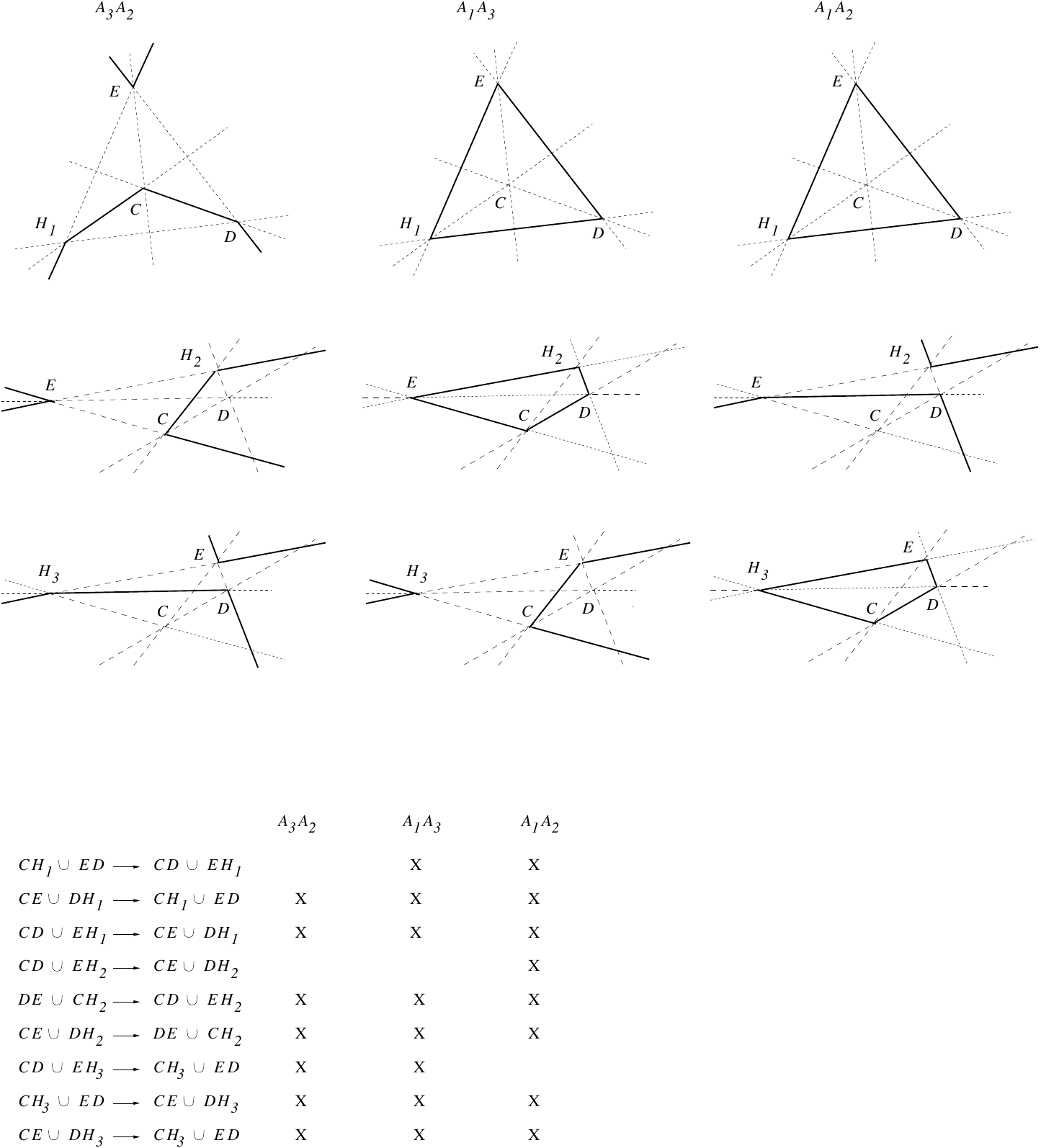}
\caption{\label{ehcd} Crossing jump, $O_1$ negative, up, pencil of conics $\mathcal{F}_{CDEH}$}
\end{figure}

{\em Proof:\/} 
For this complex type, $\lambda_1 = \lambda_2 = \lambda_3 = 1$, see lower part of Figure~\ref{crosses}. 
The triangular ovals are all interior to $cr(O_1)$ in $cr(T_0)$.
Let $H = H_2$ be $Q \to D$ in $cr(Q_2)$, see Figure~\ref{ovalinT0}.
The pencil $\mathcal{F}_{CDEH}$ has only one non-totally real portion: $CD \cup EH
\to CE \cup DH$, see Figure~\ref{ehcd}. Let $I$ be a quadrangular oval.
If $I$ is in $cr(Q_2) \cup cr(Q_3)$, $I$ is on the arc $EHD$ of the conic, hence $I$ is $Q \to D$. If $I$ is in $cr(Q_1)$, $I$ is on the arc $CE$, hence $I$ is $D \to E$.
One should have thus: $\lambda_2 - \lambda_3 = 1$ and 
$\lambda_1 = 0$, contradiction.
Let $H = H_1$ be $Q \to D$ in $cr(Q_1)$. The pencil $\mathcal{F}_{CDEH}$ has only one non-totally real portion: $CD \cup EH \to CH \cup ED$. 
If $I$ is in $cr(Q_1) \cup cr(Q_3)$, $I$ is on the arc $EHC$ of the conic, hence $I$ is $Q \to D$. If $I$ is in $cr(Q_2)$, $I$ is on the arc $DE$, hence $I$ is $D \to E$. One should have thus: $\lambda_1 - \lambda_3 = 1$ and $\lambda_2 = 0$, contradiction.
See Figures~\ref{ovalinT0}-\ref{ehcd}.
All of the ovals between $Q$ and $D$ are in $cr(Q_3)$, they contribute $-1$ to $\lambda_3$, between $D$ and $E$, we find the ovals of $cr(Q_2) \cup cr(Q_1)$,  and the remaining ovals of $cr(Q_3)$, whose contribution to $\lambda_3$ is $+2$.
Let $H_3$ be an oval $D \to E$ in $cr(Q_3)$. The pencil of conics $\mathcal{F}_{CDEH_3}$ has only one non-totally real portion: $CD \cup EH_3 \to CH_3 \cup DE$. If $I$ is in $cr(Q_3)$, $I$ is on the arc $CH_3E$ of the conic, hence $I$ is $D \to E$. There are no ovals in $cr(Q_3)$ between $Q$ and $D$, contradiction. $\Box$

\begin{lemma}
The complex type $(+, n), (\pm, \mp), (+, -, -)$ is not realizable.
For $(+, u), (\pm, \mp)$, $(+, -, -)$, the ovals in $cr(Q_2)$ are $D \to E$, the ovals in $cr(Q_1) \cup cr(Q_3)$ are $Q \to D$. Moreover, $\lambda_0 = 0$ or $1$.
\end{lemma}

{\em Proof:\/}
Hereafter, $I$ and $H$ stand for quadrangular ovals, let us write $H = H_i$ if $H$ is in $cr(Q_i)$. Assume there exists an oval $G = G_0$ that is exterior in $cr(T_0)$ or $G = G_1$ interior to $cr(O_1)$ in $cr(T_1)$.
(Note that if $\lambda_0 \not= 0$ there exists $G = G_0$, and there exists $G = G_1$ if and only if $O_1$ is up).
If all of the ovals in $cr(Q_i)$ are $Q \to D$ ($D \to E$), we say shortly that $cr(Q_i)$ is $Q \to D$ ($D \to E$). 

\begin{figure}[htbp]
\centering
\includegraphics{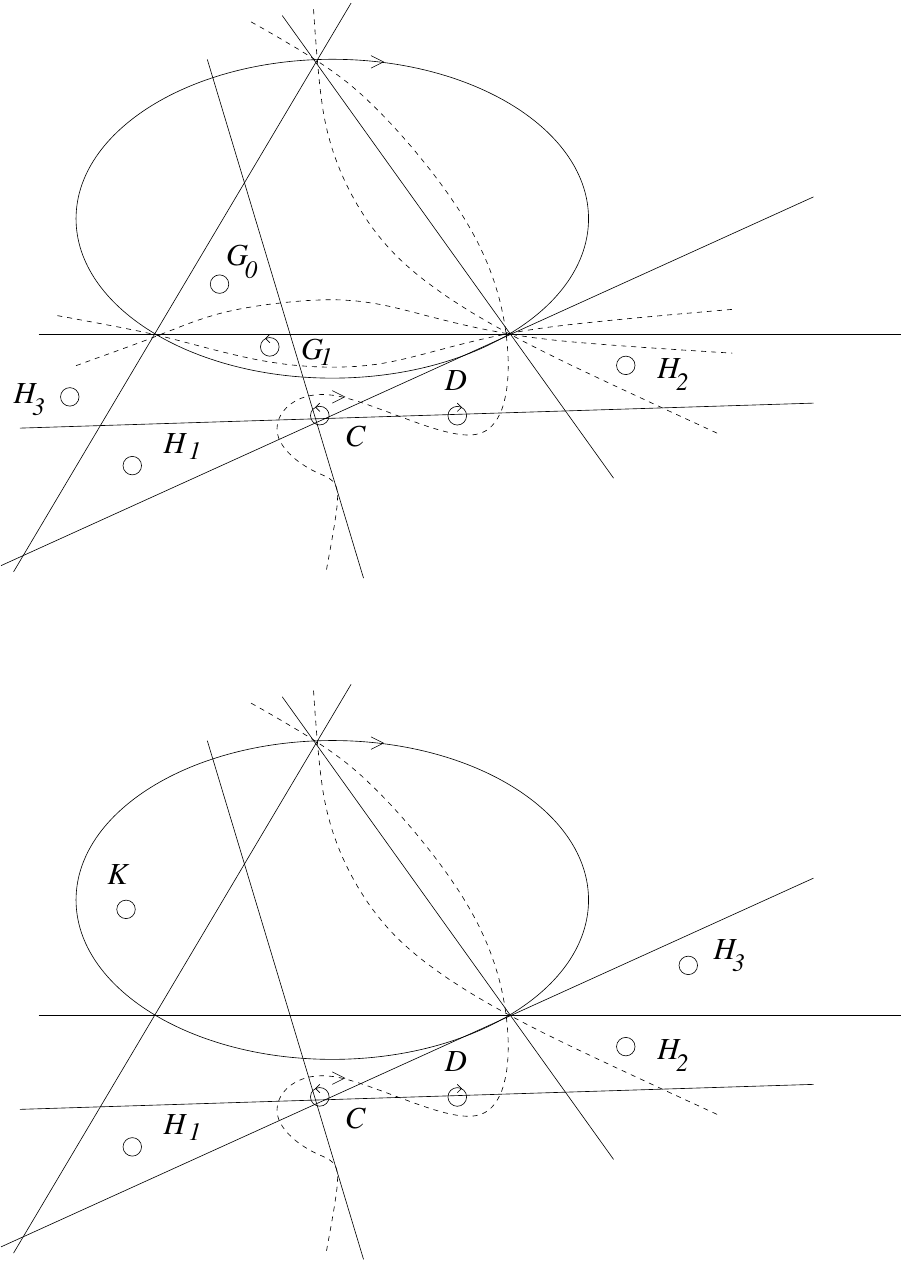}
\caption{\label{xingplus} Crossing jump, $O_1$ positive, $G$ in $cr(T_0)$ or $cr(T_1)$, $K$ in $cr(T_2)$}
\end{figure}

\begin{figure}[htbp]
\centering
\includegraphics{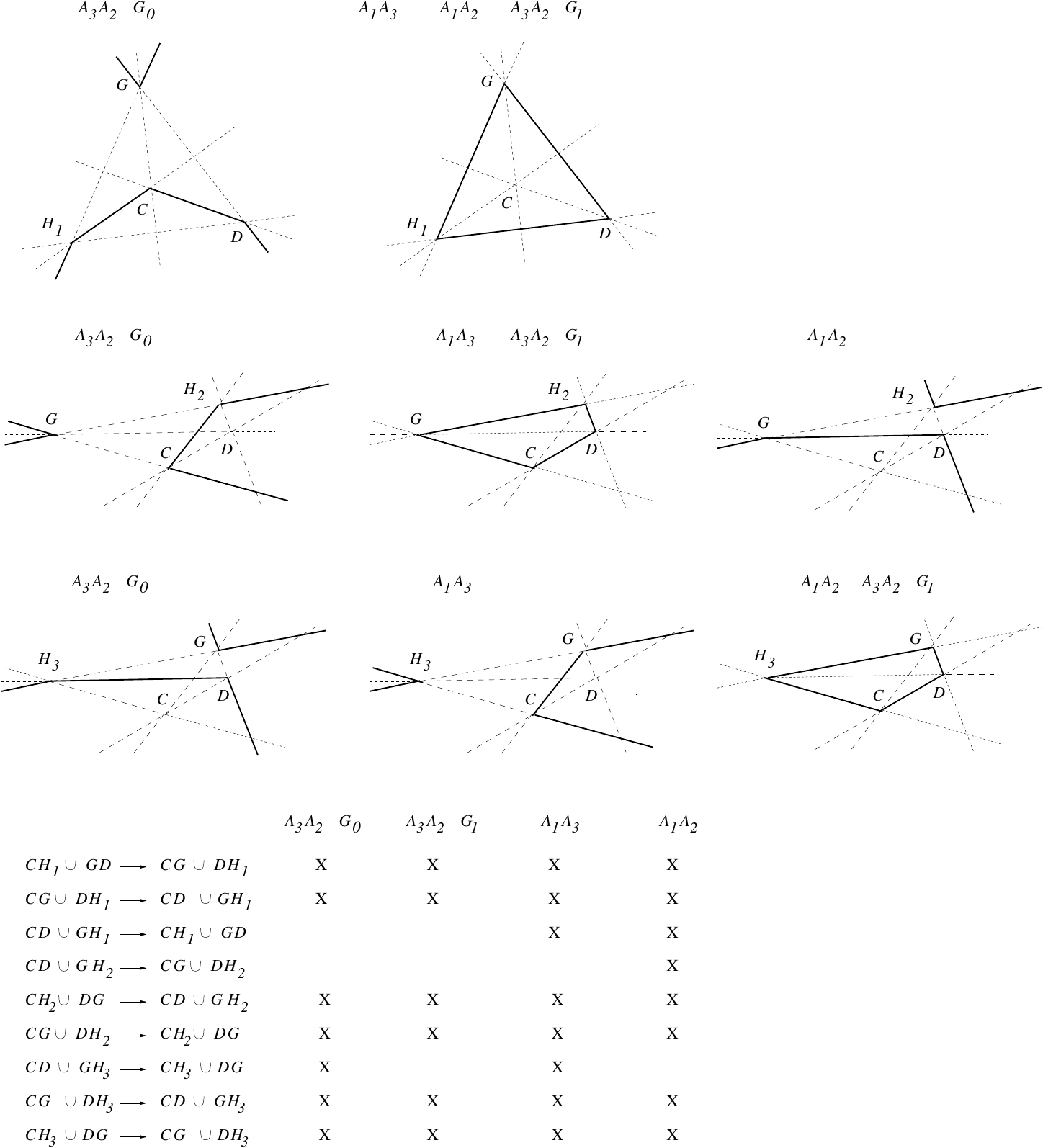}
\caption{\label{gcdh} Crossing jump, $O_1$ positive, pencil of conics $\mathcal{F}_{CDGH}$}
\end{figure}

Assume there exists $H = H_1: Q \to D$, the pencil of conics $\mathcal{F}_{CDGH_1}$ has only one non-totally real portion: $CD \cup GH_1 \to
CH_1 \cup GD$, see upper part of Figure~\ref{xingplus} and Figure~\ref{gcdh}. 
Consider a conic of this portion, passing through a quadrangular oval $I$.
If $I$ is in $cr(Q_1) \cup cr(Q_3)$, then $I$ is on the arc $GH_1C$ of the conic, hence $I$ is $Q \to D$. If $I$ is in $cr(Q_2)$, then $I$ is on the arc $DG$, hence $I$ is $D \to E$. One has thus $\lambda_1 - \lambda_3 = 1$ and $\lambda_2 = 0$. But for $O_1$ non-separating, $\lambda_1 = \lambda_3 = -\lambda_0$ and $\lambda_2 = 1$, contradiction. For $O_1$ up, $\lambda_1 = 1 - \lambda_0$, $\lambda_2 = 0$, $\lambda_3 = - \lambda_0$, there is no contradiction.
In summary:

\begin{quote}
(1) ($O_1$ non-separating and $\exists G_0$)  $\Rightarrow$ ($cr(Q_1): D \to E$).\\
(2) ($O_1$ up and $\exists H_1: Q \to D$) $\Rightarrow$ ($cr(Q_1), cr(Q_3): Q \to D$, $cr(Q_2): D \to E$).
\end{quote}

Assume there exists $H = H_2: Q \to D$, the pencil of conics $\mathcal{F}_{CDGH_2}$ has only one non-totally real portion: $CD \cup GH_2
\to CG \cup DH_2$, see upper part of Figure~\ref{xingplus} and Figure~\ref{gcdh}. 
Consider a conic of this portion, passing through a quadrangular oval $I$.
If $I$ is in $cr(Q_2)$, then $I$ is on the arc $GH_2D$ of the conic, hence $I$ is $Q \to D$. If $I$ is in $cr(Q_1)$, then $I$ is on the arc $CG$ of the conic, hence $I$ is $D \to E$. 
Let now $G = G_0$. If $I$ is in $cr(Q_3)$, then $I$ is on the arc $GH_2$ of the conic, hence $I$ is $Q \to D$. Therefore, $\lambda_2 - \lambda_3 = 1$ and $\lambda_1 = 0$. For $O_1$ non-separating, $\lambda_0 = 0$, and for $O_1$ up, $\lambda_0 = 1$.

\begin{quote}
(3) ($\exists H_2: Q \to D$ and $\exists G_0$) $\Rightarrow$ ($cr(Q_2), cr(Q_3): Q \to D$, $cr(Q_1): D \to E$, $\lambda_0 = 0$ for $O_1$ non-separating, $\lambda_0 = 1$ for $O_1$ up).\\
(4) ($\exists H_2: Q \to D$ and $O_1$ is up) $\rightarrow$ ($cr(Q_2): Q \to D$, $cr(Q_1): D \to E$).
\end{quote}

Assume there exists $H = H_3: D \to E$, the pencil of conics 
$\mathcal{F}_{CDGH_3}$ has only one non-totally real portion: $CD \cup GH_3 \to
CH_3 \cup GD$, see upper part of Figure~\ref{xingplus} and Figure~\ref{gcdh}. 
Consider a conic of this portion, passing through a quadrangular oval $I$.
If $G = G_0$, then $I$ is on the arc $CH_3G$ of the conic, hence $I$ is $D \to E$.
There are no quadrangular ovals $Q \to D$, contradiction.
If $G = G_1$ and $I$ is in $cr(Q_1) \cup cr(Q_3)$, then $I$ is on the arc $CH_3G$ of the conic, hence $I$ is $D \to E$. The zones $cr(Q_1)$ and $cr(Q_3)$ are $D \to E$, hence there exists some oval $H_2: Q \to D$, thus $cr(Q_2)$ is $Q \to D$ (4).
One has thus $\lambda_2 = 1$ and $\lambda_1 - \lambda_3 = 0$. But for $O_1$ up, $\lambda_1 = 1 - \lambda_0$, $\lambda_2 = 0$, $\lambda_3 = -\lambda_0$, contradiction.

\begin{quote}
(5) ($\exists G$) $\Rightarrow$ ($cr(Q_3): Q \to D$).\\
(6) ($\exists G$ and $\exists H_2: Q \to D$) $\Rightarrow$ ($cr(Q_2), cr(Q_3): Q \to D$ and $cr(Q_1): D \to E$, $\lambda_0 = 0$ for $O_1$ non-separating, $\lambda_0 = 1$ for $O_1$ up).\\
\end{quote}

\begin{figure}[htbp]
\centering
\includegraphics{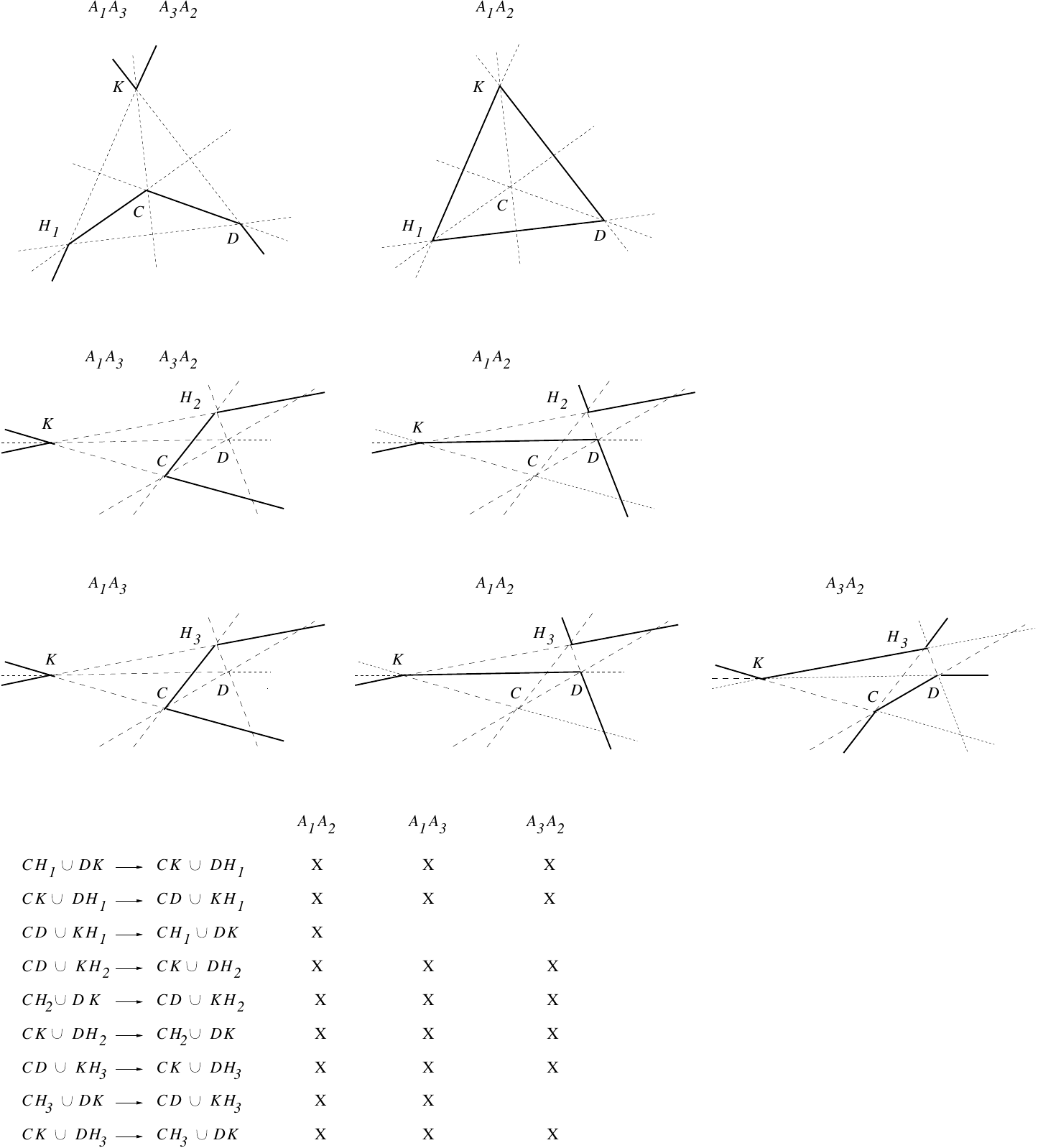}
\caption{\label{cdhk} Crossing jump, $O_1$ positive, pencil of conics $\mathcal{F}_{CDHK}$}
\end{figure}

Assume there exists an oval $K$ in $cr(T_2)$.
Let $H = H_1$ be $Q \to D$, the pencil of conics $\mathcal{F}_{CDH_1K}$ has only one non-totally real portion: $CD \cup KH_1 \to CH_1 \cup KD$, see lower part of Figure~\ref{xingplus} and Figure~\ref{cdhk}. Let $I$ be in $cr(Q_1) \cup cr(Q_3)$, $I$ is on the arc $KH_1C$ of the conic, hence $I$ is $Q \to D$. Let $I$ be in $cr(Q_2)$, $I$ is on the arc $DK$ of the conic, hence $I$ is $D \to E$. One has thus
$\lambda_1 - \lambda_3 = 1$ and $\lambda_2 = 0$, contradiction for $O_1$ non-separating. 

\begin{quote}
(7) ($\exists K$ and $O_1$ is non-separating) $\Rightarrow$ ($cr(Q_1): D \to E$).
\end{quote}

Let $H = H_2$ be $Q \to D$, the pencil of conics 
$\mathcal{F}_{CDKH_2}$ is totally real, contradiction.

\begin{quote}
(8) ($\exists K$) $\Rightarrow$ ($cr(Q_2): D \to E$).
\end{quote}

Let $H = H_3$ be $Q \to D$, the pencil of conics $\mathcal{F}_{CDH_3K}$ has only one non-totally real portion: $CH_3 \cup DK \to CD \cup KH_3$. Let $I$ be in $cr(Q_1) \cup cr(Q_3)$, $I$ is on the arc $KH_3C$ of the conic, hence $I$ is $Q \to D$. Let $I$ be in $cr(Q_2)$, then $I$ is on the arc $DK$ of the conic, hence $I$ is $D \to E$. One has thus $\lambda_1 - \lambda_3 = 1$ and $\lambda_2 = 0$, contradiction for $O_1$ non-separating.

\begin{quote}
(9) ($\exists K$ and $\exists H_3: Q \to D$) $\Rightarrow$ ($cr(Q_1)$, $cr(Q_3): Q \to D$, $cr(Q_2): D \to E$, $O_1$ is up)\\
(10) ($\exists K$ and $O_1$ is non-separating) $\Rightarrow$ ($cr(Q_1)$, $cr(Q_2)$, $cr(Q_3): D \to E$), contradiction.
\end{quote}

For $O_1$ non-separating, the zone $cr(T_2)$ must be empty (10). As $\lambda_5 = 1 + \lambda_0$, one has $\lambda_0 = -1$, so there exists $G = G_0$ in $cr(T_0)$.
The ovals in $cr(Q_1)$ are $D \to E$ (1). As $\lambda_0 \not= 0$, the ovals in $cr(Q_2)$ are also $D \to E$ (3). The ovals in $cr(Q_3)$ are $Q \to D$ (5). One should have thus $-\lambda_3 = \lambda_0 = 1$, this is a contradiction. The complex type with non-separating $O_1$ is not realizable.

For $O_1$ up, the ovals in $cr(Q_1)$ are all $Q \to D$ or $D \to E$ (2). Assume there exists $H = H_2: Q \to D$, then on one hand $\lambda_0 = 1$ (6), on the other hand, the zone $cr(T_2)$ must be empty (8). But
$\lambda_5 = \lambda_0$, contradiction. 
Therefore, $cr(Q_2)$ is $D \to E$.
The zone $cr(Q_3)$ is $Q \to D$ (5). If the ovals in $cr(Q_1)$ are $D \to E$, we have $-\lambda_3 = \lambda_0 = \lambda_5 =1$ hence there exist $H_3$ in $cr(Q_3)$ and $K$ in $cr(T_2)$. 
But the existence of two ovals $H_3: Q \to D$ and $K$ implies that $cr(Q_1)$ is $Q \to D$ (9).

\begin{quote}
(11) $O_1$ is up $\Rightarrow$ $cr(Q_1), cr(Q_3): Q \to D$, $cr(Q_2): D \to E$.
\end{quote}
 
As $\lambda_1 - \lambda_3 = 1$, there exists an oval $H$ that is $Q \to D$ in $cr(Q_1)$ or in $cr(Q_3)$. Let $H =H_1$, recall that for any choice of $G$ in $cr(T_0) \cup cr(T_1)$, the pencil of conics $\mathcal{F}_{CDGH_1}$ has only one non-totally real portion: $CD \cup GH_1 \to GD \cup CH_1$. An oval $K$ in $cr(T_2)$ will be on the arc $DG$ of its conic, see Figure~\ref{gcdh}. 
Let now $H = H_3$ and $K$ be in $cr(T_2)$, recall that the pencil of conics $\mathcal{F}_{CDKH_3}$ has only one non-totally real portion: $CD \cup KH_3 \to KD \cup DH_3$. The oval $G$ is on the arc $KH_3$ of its conic, see Figure~\ref{cdhk}.
In all cases, the pencil $\mathcal{F}_C: D \to A_3 \to A_1 \to A_2$ meets the ovals in $cr(T_2)$ before those in $cr(T_0) \cup cr(T_1)$, the first oval met ($E$) is positive, hence $\lambda_5 = \lambda_0 = 0$ or $+1$. $\Box$ 
  


This finishes the proof of Theorem~2.

\end{document}